\documentclass{article}
\usepackage{Arxiv}

\usepackage{amsmath,amssymb,amsfonts}
\usepackage{rotating}
\usepackage{fancyvrb}
\usepackage{algorithmic}
\usepackage{amsopn}

\usepackage{mathtools}
\usepackage{multirow}
\newtheorem{result}{Main Result}
\usepackage{enumitem,verbatim}
\usepackage{algorithm}
\usepackage{algorithmic}
\usepackage{dsfont}
\usepackage{pifont}
\newcommand{\cmark}{\ding{51}}%
\newcommand{\xmark}{\ding{55}}%
\usepackage{hyperref}
\usepackage[para,online,flushleft]{threeparttable}
\usepackage{todonotes}
\usepackage{lineno}
\usepackage{array}
\usepackage[authoryear]{natbib}
\usepackage{amsthm}

\def\grad{\nabla}

\newcommand{\bq}{\mathbf{q}}
\newcommand{\br}{\mathbf{r}}
\newcommand{\bs}{\mathbf{s}}

\newcommand{\bu}{\mathbf{u}}
\newcommand{\bv}{\mathbf{v}}
\newcommand{\bw}{\mathbf{w}} 
\newcommand{\bx}{\mathbf{x}}
\newcommand{\by}{\mathbf{y}}
\newcommand{\bz}{\mathbf{z}}
\newcommand{\bA}{\mathbf{A}}
\newcommand{\bB}{\mathbf{B}}
\newcommand{\bC}{\mathbf{C}}
\newcommand{\bD}{\mathbf{D}}

\newcommand{\bF}{\mathbf{F}}

\newcommand{\bI}{\mathbf{I}}

\newcommand{\bL}{\mathbf{L}}
\newcommand{\bM}{\mathbf{M}}

\newcommand{\bS}{\mathbf{S}}
\newcommand{\bU}{\mathbf{U}}
\newcommand{\bT}{\mathbf{T}}

\newcommand{\bV}{\mathbf{V}}

\def\ones{\mathbf{1}}

\def\cA{\mathcal{A}}
\def\cB{\mathcal{B}}
\def\cC{\mathcal{C}}
\def\cD{\mathcal{D}}

\def\cF{\mathcal{F}}
\def\cG{\mathcal{G}}
\def\cH{\mathcal{H}}
\def\cI{\mathcal{I}}

\def\cK{\mathcal{K}}
\def\cL{\mathcal{L}}
\def\cM{\mathcal{M}}
\def\cN{\mathcal{N}}
\def\cO{\mathcal{O}}
\def\cP{\mathcal{P}}

\def\cR{\mathcal{R}}
\def\cS{\mathcal{S}}

\def\cU{\mathcal{U}}

\def\cX{\mathcal{X}}
\def\cY{\mathcal{Y}}
\def\cZ{\mathcal{Z}}

\def\mE{\mathbb{E}}

\def\mS{\mathbb{S}}

\def\mZ{\mathbb{Z}}

\def\smskip{\smallskip}

\def\texitem#1{\par\smskip\noindent\hangindent 25pt
               \hbox to 25pt {\hss #1 ~}\ignorespaces}

\def\abs#1{\left|#1\right|}
\def\norm#1{\left\|#1\right\|}

\newcommand{\BEAS}{\begin{eqnarray*}}
\newcommand{\EEAS}{\end{eqnarray*}}
\newcommand{\BEA}{\begin{eqnarray}}
\newcommand{\EEA}{\end{eqnarray}}
\newcommand{\BEQ}{\begin{eqnarray}}
\newcommand{\EEQ}{\end{eqnarray}}
\newcommand{\BIT}{\begin{itemize}}
\newcommand{\EIT}{\end{itemize}}
\newcommand{\BNUM}{\begin{enumerate}}
\newcommand{\ENUM}{\end{enumerate}}

\newcommand{\BA}{\begin{array}}
\newcommand{\EA}{\end{array}}

\newcommand{\reals}{\mathbb{R}}
\newcommand{\integers}{\mathbb{Z}}

\newcommand{\diag}{\mathop{\bf diag}}

\newcommand{\dom}{\mathop{\bf dom}}

\newif\ifpagenumbering
\pagenumberingtrue

\pagenumberingfalse

\newsavebox{\theorembox}
\newsavebox{\lemmabox}
\newsavebox{\defnbox}
\newsavebox{\rsltbox}
\newsavebox{\corollarybox}
\newsavebox{\remarkbox}
\newsavebox{\assbox}
\savebox{\theorembox}{\noindent\bf Theorem}
\savebox{\lemmabox}{\noindent\bf Lemma}
\savebox{\defnbox}{\noindent\bf Definition}
\savebox{\rsltbox}{\noindent\bf Main Result}
\savebox{\corollarybox}{\noindent\bf Corollary}
\savebox{\remarkbox}{\noindent\bf Remark}
\savebox{\assbox}{\noindent\bf Assumption}
\newtheorem{assumption}{\usebox{\assbox}}
\newtheorem{remark}{\usebox{\remarkbox}}
\newtheorem{theorem}{\usebox{\theorembox}}
\newtheorem{lemma}{\usebox{\lemmabox}}
\newtheorem{corollary}{\usebox{\corollarybox}}
\newtheorem{defn}{\usebox{\defnbox}}

\def\fprod#1{\left\langle#1\right\rangle}

\def\ind#1{\mathbb{I}_{#1}}

\def\id{{\bf I}}
\definecolor{darkred}{RGB}{0,128,0}
\def\sa#1{\textcolor{black}{#1}}
\def\nsa#1{\textcolor{black}{#1}}
\def\rev#1{\textcolor{black}{#1}}
\def\us#1{\textcolor{black}{#1}}

\def\aj#1{\textcolor{black}{#1}}
\def\eyh#1{\textcolor{black}{#1}}
\def\missing#1{\textcolor{black}{#1}}
\definecolor{aa}{rgb}{0.54, 0.17, 0.89}
\newcommand{\ey}[1]{\black{#1}}
\def\eyz#1{\textcolor{black}{#1}}
\def\btx{\tilde{\bx}}
\def\bty{\tilde{\by}}
\def\btr{\tilde{\br}}

\def\bdelta{\boldsymbol{\delta}}

\DeclareMathOperator*{\argmin}{argmin}

\def\rv#1{\textcolor{black}{#1}}
\def\eh#1{\textcolor{black}{#1}}
\def\ey#1{\textcolor{black}{#1}}
\def\saa#1{\textcolor{black}{#1}}
\def\erfan#1{\textcolor{black}{#1}}
\def\sp#1{\textcolor{black}{#1}}

\def\rev#1{\textcolor{black}{#1}}

\def\oZ{\overline Z}
\def\oX{\overline X}
\def\oY{\overline Y}
\makeatletter
\newcommand{\thickhline}{%
    \noalign {\ifnum 0=`}\fi \hrule height 1pt
    \futurelet \reserved@a \@xhline
}
\newcolumntype{"}{@{\hskip\tabcolsep\vrule width 1pt\hskip\tabcolsep}}
\makeatother
\allowdisplaybreaks

\title{A Randomized Block-Coordinate Primal-Dual Method for Large-scale Stochastic Saddle Point Problems}
\shorttitle{This is my short title}
\author{
Erfan Yazdandoost Hamedani\\
Systems and Industrial Engineering Department\\
University of Arizona
\\
Tucson, AZ, USA.\\
\texttt{erfany@arizona.edu}
\And
Afrooz Jalilzadeh\\
Systems and Industrial Engineering Department\\
University of Arizona
\\
Tucson, AZ, USA.\\
\texttt{afrooz@arizona.edu}
\And
Necdet Serhat Aybat\\
Department of Industrial and Manufacturing Engineering\\
Pennsylvania State University
\\
University Park, PA,USA.\\
\texttt{nsa10@psu.edu}
}
\begin{document}
\maketitle
\begin{abstract}
    We consider (stochastic) convex-concave saddle point (SP) problems with high-dimensional decision variables, arising in various applications including machine learning problems. To contend with the challenges in computing full gradients, we employ a randomized block-coordinate primal-dual scheme in which randomly selected primal and dual blocks of variables are updated. We consider both deterministic and stochastic settings, where deterministic partial gradients and their randomly sampled estimates are used, respectively, at each iteration. We investigate the convergence of the proposed method under different blocking strategies and provide the corresponding complexity results. While the best-known \rev{computational complexity result for computing a saddle point with $\varepsilon$ primal-dual gap for deterministic primal-dual methods using full gradients is $\mathcal O(\max\{m,n\}^2/\varepsilon)$, where $m$ and $n$ denote the dimensions of primal and dual variables, respectively,} we show that our proposed randomized block-coordinate method achieves an improved complexity of {$\mathcal O(mn/\varepsilon)$ assuming a coordinate-friendly structure on the problem.} Moreover, for the stochastic setting where a mini-batch sample gradient is utilized, we show a \rev{computational complexity of {$\tilde{\mathcal{O}}(m^2n^2/\varepsilon^2)$}} through acceleration. Finally, 
 almost sure convergence of the iterate sequence to a saddle point is established.
 \end{abstract}
\section{Introduction}
\label{sec:intro}
\rev{Saddle point (SP) problems 
are particular minimax optimization problems that arise in a wide range of modern applications due to their ability to model complex interactions and constraints. Notably, they provide a unifying framework for various optimization problems encountered in machine learning, control, game theory, and signal processing. Their versatile formulation encompasses convex optimization problems with nonlinear conic constraints, which subsume classical models such as linear programming (LP), quadratic programming (QP), quadratically constrained quadratic programming (QCQP), second-order cone programming (SOCP), and semidefinite programming, making them fundamental to both theory and practice.} \rev{Next, we formally introduce the problem we consider in this paper.}

Let $(\cX_i,\norm{\cdot}_{\cX_i})$  and $(\cY_j,\norm{\cdot}_{\cY_j})$ be finite-dimensional, normed vector spaces for $i\in\cM\triangleq \{1,2 \break ,\hdots,M\}$ and $j\in\cN\triangleq \{1,2,\hdots,N\}$, respectively, where $\cX_i\triangleq\reals^{m_i}$ for $i\in\cM$ and $\cY_j\triangleq\reals^{n_j}$ for $j\in\cN$. Let 
$\bx=[x_i]_{i\in\cM}\in\Pi_{i\in\cM}\cX_i\triangleq\cX$ and 
$\by=[y_j]_{j\in\cN}\in\Pi_{j\in\cN}\cY_j\triangleq \cY$ with dimensions $m\triangleq \sum_{i\in\cM}m_i$ and $n\triangleq \sum_{j\in\cN}n_i$, respectively; moreover, let $\cZ\triangleq \cX\times \cY$. 
We consider the following 
SP
problem: \vspace*{0mm}
\begin{align}
\label{eq:original-problem}
{\rm (SP)}:\quad \min_{\bx\in\cX}\max_{\by\in\cY}\ \  \cL(\bx,\by)&\triangleq \sa{f(\bx)}+\Phi(\bx,\by)-\sa{h(\by)},
\end{align}
\sa{where \sa{$f(\bx)\triangleq\frac{1}{M}\sum_{i\in\cM}f_i(x_i)$ and $h(\by)\triangleq\frac{1}{N}\sum_{j\in\cN}h_j(y_j)$ such that} $f_i:\cX_i\to\reals\cup\{+\infty\}$ for all $i\in\cM$ and $h_j:\cY_j\to\reals\cup\{+\infty\}$ for all $j\in\cN$ are (possibly nonsmooth) closed convex functions.} 
\saa{Suppose $\Phi(\bx,\by)$
is convex in $\bx$, concave in $\by$; moreover, we assume that we only have access to an \textit{unbiased} stochastic first-order oracle~(SFO) $\grad \Phi(\cdot,\cdot;\xi(\omega))$ for $\Phi$ with a \textit{bounded variance}, where $\xi:\Omega\to\reals$ is a random variable and $(\Omega,\cF,\mathbb{P})$ denotes the associated probability space (see Assumption~\ref{assum:sample}). 
\sp{For instance, when the objective function is an expectation of a stochastic function, since evaluating the gradient of the objective function is expensive, a constant batch-size of \saa{stochastic partial} gradients of $\Phi$ can be used at each iteration \saa{to approximate $\grad_\bx\Phi$ and $\grad_\by\Phi$}.}
Finally, we also assume that} $\grad\Phi(\cdot,\cdot;\xi(\omega))$
satisfies certain differentiability assumptions for all $\omega\in\Omega$
{(see Assumption~\ref{assum-lip})}.

\rev{In many practical applications, SP problems naturally arise in large-scale settings, including high-dimensional computational statistics, machine learning, and distributed computing. For instance, 
{\bf (i)} {\it Robust optimization} {contends} with
{uncertain 
parameters with the aim of minimizing the worst case (maximum) value of the objective function for all realizations of the uncertain parameter that is known to belong to} a given ambiguity set~\citep{ben2009robust}, e.g., robust classification problem \citep{namkoong2016stochastic,shalev2016minimizing} {leads to {an SP} 
problem} -- see Section~\ref{numerics} for details; {\bf (ii)} {\it Distance metric learning} proposed in \citep{xing2003distance} is a convex optimization problem over positive-semidefinite matrices subject to nonlinear convex constraints; {\bf (iii)} {\it Kernel matrix learning} for transduction can be cast as an SDP or a QCQP \citep{lanckriet2004learning,gonen2011multiple}; {\bf (iv)} {\it Training of ellipsoidal kernel machines} \citep{shivaswamy2007ellipsoidal} requires solving nonlinear SDPs. 
These problems often involve high-dimensional variables and {functional constraints, which further complicate their solution through rendering 
the traditional first-order methods building on computing full gradients costly} 
and memory-intensive. As a result, randomized block-coordinate methods have gained popularity due to their scalability and efficiency. By updating only selected blocks of the primal and dual variables at each iteration, these methods significantly reduce per-iteration cost while preserving convergence guarantees under suitable conditions.}


Our algorithm design \us{is inspired by} primal-dual \us{methods} for deterministic bilinear convex-concave 
\nsa{SP} problems proposed 
by~\cite{chambolle2011first}, where 
\sp{the restricted gap metric $G_Z:\cZ\to\reals_+$, i.e., \vspace*{0mm}
\begin{equation}
\label{s-gap}
\begin{aligned}
G_Z(\bar{\bz})\triangleq\sup_{\nsa{\bz\in Z}}\{\cL(\bar{\bx},\by)-\cL(\bx,\bar{\by})\},\quad\forall~\bar\bz\in \cZ,
\end{aligned}
\vspace*{0mm}
\end{equation}
is shown to diminish} at \rev{the non-asymptotic rate of $\mathcal{O}(1/k)$} \sp{for any bounded $Z\subset\cZ$ that contains a saddle point of \eqref{eq:original-problem}}. 
Finite-sum or expectation-valued problems have also been studied by 
~\cite{chen2014optimal,xu2018primal}, 
leading to rates of \rev{$\mathcal{O}({1}/\sqrt{k})$}. We develop a primal-dual (double) randomized block-coordinate framework \us{for} large-scale problems characterized by \sa{(i) a \textit{deterministic} or a \textit{stochastic} objective function; (ii) a \textit{non-bilinear} coupling function $\Phi$; and (iii) a possibly \textit{large} number of primal-dual block coordinates.} Such problems can capture large-scale stochastic optimization problems with a possibly large number of convex constraints. To this end, we propose 
a randomized block-coordinate algorithm for precisely such problems \us{which achieves} an optimal convergence rate in the deterministic setting. Further efficiencies can be garnered through block-specific steplengths and sampling for the stochastic regime.  

\subsection{ Related Work} 
{There has been a 
vast body of work on randomized block-coordinate \sa{(RBC)} 
schemes for \sa{dealing with large-scale} primal optimization problems, \sa{e.g.,} \citep{nesterov2012efficiency,xu2013block,richtarik2014iteration}; but, there are far fewer studies on 
\sa{RBC} algorithms designed for solving SP problems. 
In another line of research, \sa{for minimizing a function in finite-sum or expectation forms,} variance-reduced schemes have been studied for dampening the gradient noise and achieving deterministic rates of convergence {for different scenarios depending on the convexity properties of the objective function:} strongly convex \citep{friedlander2012hybrid,byrd2012sample},
convex \citep{ghadimi2016accelerated,jalilzadeh2018optimal}, and nonconvex optimization \citep{ghadimi15minibatch,lei20asynchronous}. 
Block-based schemes with (randomly generated) block-specific batch sizes are \nsa{employed by \cite{lei20asynchronous} 
achieving the} optimal deterministic rates of convergence. 
Our schemes are based on the accelerated primal-dual algorithm (APD) 
~\citep{hamedani2021primal} which was proposed \saa{for a special case of 
\eqref{eq:original-problem} with} $M=N=1$. In contrast with APD, the methods proposed in this paper can exploit the special block structure whenever $M,N>1$ and they rely on 
{mini-batch samples in the} stochastic setting. Table~\ref{TAB-lit2} provides a review of related 
block-based schemes.}

\begin{table}[!htb]
\renewcommand{\arraystretch}{0.9}
{
	\rev{Compitational complexity} results for primal-dual schemes updating (i) full primal-dual coordinates (PD), (ii) randomized primal (RP),  (iii) randomized dual (RD), and (iv) randomized primal and dual (RPD) block coordinates.\label{TAB-lit2}}
{
  \begin{threeparttable}
  \centering
\newcommand{\tabincell}[2]{\begin{tabular}{@{}#1@{}}#2\end{tabular}}
 \small
  \renewcommand{\arraystretch}{1.28} 
 \begin{tabular}{|c|c|c|c|c|c|c|c|}
        \thickhline
        Type &  References & Setting  & \tabincell{l}{Block\\based\\step}&  \tabincell{l}{\hspace*{-1mm}Non-\\\hspace*{-1mm}bilinear}& \tabincell{l} {Unbounded\\domain}&Deterministic & Stochastic     \\ \thickhline
   \multirow{3}{*}{\tabincell{l}{ \\ \\ \\ \sa{PD}}} & \citep{chen2014optimal} & Convex-Concave& \xmark & \xmark & \cmark&$\cO(1/\varepsilon)$  
										 & $\widetilde\cO(1/\varepsilon^2)$
\\
 \cline{2-8} &
\citep{chambolle2016ergodic} 
& Convex-Concave & \xmark & \xmark  & \xmark&$\cO(1/\varepsilon)$  & --\\
\cline{2-8}
 & \citep{he2016accelerated}
    & Convex-Concave  & \xmark & \xmark & \cmark&$\cO(1/\varepsilon)$ & -- \\
			\cline{2-8}
 & 
 \citep{kolossoski2017accelerated} 
 & Convex-Concave& \xmark & \sa{\cmark} & \xmark&$\cO(1/\varepsilon)$ & --\\
            \cline{2-8}
 & \citep{hamedani2021primal} & Convex-Concave & \xmark   & \sa{\cmark} & \cmark&$\cO(1/\varepsilon)$ &
			-- \\
			\cline{2-8}
 & {\citep{yu2017primal}} & \tabincell{l}{CO}  & \cmark & \rev{\cmark$^*$} &\xmark&$\cO(1/\varepsilon)$ &--\\
\cline{2-8} &   \citep{zhao2022accelerated}  &  Convex-Concave & \xmark & \cmark&\xmark & $\cO(1/\varepsilon)$ & $\widetilde\cO(1/\varepsilon^2)$\\ \thickhline
\multirow{4}{*}{\tabincell{l}{\\ \\ \\ \sa{RP} \\
	\sa{RD}}}  
 &   \citep{chambolle2017stochastic}  &  Convex-Concave& \cmark  & \xmark & \xmark&$\cO(1/\varepsilon)$ & --
					\\
	\cline{2-8} 
 & \citep{xu2017first} & {CO many constraints} & \cmark & \cmark &\xmark& $\cO(1/\varepsilon)$ & -- \\ 
	\cline{2-8} & \citep{xu2018primal} & \tabincell{l}{ CO many constraints} & \xmark  & \cmark &\xmark& -- & $\widetilde\cO(1/\varepsilon^2)$\\ 
 \cline{2-8} & \citep{alacaoglu2022complexity} &  Convex-Concave & \cmark & \xmark &\cmark & $\cO(1/\varepsilon)$ & --\\
 \cline{2-8} & \citep{alacaoglu2022convergence} &  Convex-Concave & \cmark & \xmark &\cmark & $\cO(1/\varepsilon)$ & --\\
 \cline{2-8} & \citep{hamedani23a} &  Convex-Concave & \cmark & \cmark & \xmark & $\cO(1/\varepsilon)$ & --\\ 
 \cline{2-8} &\citep{zhang2023primal} &  Convex-Concave & \xmark & \cmark&\xmark & $\cO(1/\varepsilon)$ & --\\ 
\thickhline
\multirow{1}{*}{\tabincell{l}{
	\sa{RPD}}}  
 & \tabincell{l}{ Our method} & Convex-Concave & \cmark & \cmark &\cmark & $\cO(1/\varepsilon)$  & $\widetilde\cO(1/\varepsilon^2)$  \\
 \thickhline
 \end{tabular}
  \end{threeparttable}}
  {The settings include convex-concave, 
  and constrained optimization (CO) problems. \rev{$^*$\citep{yu2017primal} considers CO problems; therefore, the coupling function is linear in the dual variable.}}
\end{table}
\vspace*{0mm}
\subsection{Discussion on Computational Complexity of Different Partitioning Strategies} 
\label{sec:discussion}

\rev{In this paper, similar to \cite{nesterov2012efficiency,peng2016coordinate} we 
	make the following structural assumption.
\begin{assumption}\label{assum:comp}
For any sample $\xi=(\xi^x,\xi^y)$, computing $\grad_{x_i}\Phi(\cdot,\cdot;\xi^x)$ requires one {primal oracle} call {with computational complexity $\cC_p$} for any $i\in\cM$, and computing $\grad_{y_j}{\Phi}(\cdot,\cdot;\xi^y)$
requires one {dual oracle} call {with computational complexity $\cC_d$} for any $j\in\cN$.
{Moreover, the complexity for computing the full partial gradients $\grad_\bx\Phi(\bx,\by;\xi^x)$ and $\grad_\by\Phi(\bx,\by;\xi^y)$ are $\cO(M\cC_p)$ and $\cO(N\cC_d)$, respectively, for any sample $\xi$.} 
{The same coordinate-friendly structure also holds for the deterministic setting.}
\end{assumption}
To simplify the discussion throughout the manuscript, suppose that the cost of primal-oracles are the same for all $i\in\cM$ and the cost of dual-oracles are the same for all $j\in\cN$.
Based on the above assumption, we formally define \textit{computational complexity} as the number of primal and dual partial (sample) gradients used by an algorithm for finding the $\epsilon$-gap solution. As an example, consider $\min_{\bx\in\cX}\max_{\by\in\Delta}\Phi(\bx,\by)$ where $\Phi(\bx,\by)=\sum_{\ell=1}^Q y_\ell {h_\ell(x_\ell)}$ for some $Q\in\integers_+$, $h_\ell(\cdot)$ is $L_\ell$-smooth convex for $\ell\in\{1,\cdots, Q\}$, $\Delta\subset\cY=\reals^Q$ is a $Q$-dimensional unit simplex; thus, $\Phi$ is a convex-concave function. Considering $M=N=Q$, 
note $\grad_{x_\ell}\Phi(\bx,\by)=y_\ell\grad {h_\ell(x_\ell)}$ and $\grad_{y_\ell}\Phi(\bx,\by)= {h_\ell(x_\ell)}$ for any $\ell$; thus, computing $\grad_{x_\ell}\Phi$ and $\grad_{y_\ell}\Phi$ require $1/Q$ computational effort compared to evaluations of 
$\grad_\bx\Phi$ and 
$\grad_\by\Phi$,
respectively. Thus, this example has a coordinate-friendly structure satisfying Assumption~\ref{assum:comp}.}

\rev{We discuss how the computational complexity of 
our method depends on primal and dual dimensions $m,n\gg1$. 
We {consider} four partitioning strategies for both deterministic 
and stochastic 
settings:\\  
\textit{(i)} {\bf (Block $\bx$, $\by$)} partitions both $\bx$ and $\by$, {i.e., $\cM=\{1,\ldots, M\}$ and $\cN=\{1,\ldots,N\}$ for some $M,N>1$. For simplicity, let $M=m$, $N=n$, i.e., $m_i=1$ for $i\in\cM$ and $n_j=1$ for $j\in\cN$.
 {We assume that the complexity for computing partial gradients $\grad_{x_i}\Phi$ and $\grad_{y_j}\Phi$ (and their unbiased stochastic estimates) are $\cC_p$ and $\cC_d$ for all $i\in\cM$ and $j\in\cN$ and that all coordinate-wise Lipschitz constants are $\cO(1)$ --see Assumption~\ref{assum-lip}. Moreover, for the stochastic case, suppose the variance of these block-wise partial gradient estimates, i.e., $\grad_{x_i}\Phi(\cdot,\cdot;\xi^x)$ and $\grad_{y_j}{\Phi}(\cdot,\cdot;\xi^y)$, are uniformly bounded by $\delta^2\geq 0$ for all $i\in\cM$ and $j\in\cN$.}\\
\textit{(ii)} {\bf (Block $\bx$)} partitions $\bx$ and uses full (stochastic) gradient with respect to $\by$, {i.e., let $M=m$ and $N=1$ leading to $\cM=\{1,\ldots,m\}$ and $\cN=\{1\}$ with $m_i=1$ for $i\in\cM$ and $n_1=n$; thus, the cost of computing a full (stochastic) gradient with respect to $\by$ is $\cO(n\cC_d)$ with a variance bound of $n\delta^2$.}\\ 
\textit{(iii)} {\bf (Block $\by$)} partitions $\by$ and uses full (stochastic) gradient with respect to $\bx$, {i.e.,  let $M=1$ and $N=n$ leading to $\cM=\{1\}$ and $\cN=\{1,\ldots,n\}$ with $m_1=m$ and $n_j=1$ for $j\in\cN$}; thus, the cost of computing a full (stochastic) gradient with respect to $\bx$ is $\cO(m\cC_p)$ with a variance bound of $m\delta^2$.}\\ 
\textit{(iv)} {\bf (Full)} uses full primal and dual {(stochastic) gradients, i.e., let $M=N=1$ leading to $\cM=\{1\}$ and $\cN=\{1\}$ with $m_1=m$ and $n_1=n$}; thus, the cost of computing full (stochastic) gradients with respect to $\bx$ and $\by$ are $\cO(m\cC_p)$ and $\cO(n\cC_d)$ with  variance bounds of $m\delta^2$ and $n\delta^2$, respectively.}  

{\sp{A key consideration in comparing these strategies is to account for the effects of the number of blocks and of block dimension on the magnitude of Lipschitz constants (see Remark~\ref{rem:L-constants}) and on the variance bound of the stochastic estimate,
both of which determine the $\cO(1)$ constant of the complexity bound (see Corollaries~\ref{cor:stoch-complexity} and \ref{col:rate}).}
\sp{The computational complexity for each of these strategies are provided in Table \ref{tab:oracle-comp} and their derivations are explained in detail in Section \ref{sec:complexity} of the Appendix, where} 
\rev{we consider the scenario with both dimensions large, i.e., $m,n\gg 1$, and to be able to count the number of primal and dual oracle calls together, we assumed that the costs of primal and dual oracles are comparable, i.e., $\cC_p\approx\cC_d$ for the sake of simplicity.} The key takeaways from this comparison are as follows:} 
\begin{itemize}[noitemsep,topsep=0pt,leftmargin=*,align=left]
    \item {When $m,n\gg 1$ such that $m\neq n$}, partitioning both variables into multiple blocks may lead to a lower computational complexity. We provide two examples below:
    \begin{itemize}
        \item For $\delta=0$, 
        with {$m=n^2$ and $n\gg 1$,} then ({\bf Block $\bx,\by$}) strategy {with $M=m$ and $N=n$ gives $\cO(n^3/\varepsilon)$} while 
    	{({\bf Full}) and ({\bf Block $\bx$}) strategies lead to $\cO(n^4/\varepsilon)$, and ({\bf Block $\by$}) strategy leads to $\cO(n^5/\varepsilon)$.}
        \item For $\delta>0$, when $m=n^a$ for $n\gg 1$ or when $n=m^a$ for $m\gg 1$ {for some $\integers_+\ni a\geq 2$}, 
    ({\bf Block $\bx,\by$}) with $M=m$ and $N=n$ 
    has the lowest complexity.
    For \sp{$a=2$}, while (\textbf{Full}) and (\textbf{Block $\bx,\by$}) have the same complexity, 
    one may not be able to implement (\textbf{Full}) 
    due to its high per-iteration complexity. 
    For \sp{$a>2$}, (\textbf{Block $\bx,\by$}) is strictly better than the other blocking strategies.
    \end{itemize}
    \item For $\delta=0$, when $m=n\gg 1$, 
    {\bf (Block $\bx$, $\by$)} {with $M=m$ and $N=n$}
    has the same complexity as ({\bf Full}) strategy, \sp{and is better than} the other two strategies. 
    Compared to ({\bf Full}), one may still reside with ({\bf Block $\bx,\by$}) 
    as it has a cheaper per-iteration complexity and possibly larger stepsizes reliant on the block coordinate Lipschitz constants \rev{(See Remark~\ref{rem:long-steps}).}
    \item For both $\delta=0$ and $\delta>0$ cases, 
    	{({\bf Block $\bx,\by$}) with $M=1$ reduces to ({\bf Block $\by$)}, and with $N=1$ it reduces to ({\bf Block $\bx$}), having the same computational complexities as these strategies.} 
\end{itemize}
\begin{table}[h!]
    \centering
    \caption{\rev{Computational complexity of our method under different block partitioning strategies for deterministic and stochastic settings where $\mathfrak{Q}\triangleq \max\{m,n\}$ and $\mathfrak{R}\triangleq \max\{\frac{m}{n},\frac{n}{m}\}$}.}
    \small
    \resizebox{0.99\textwidth}{!}{
    \renewcommand{\arraystretch}{1.5}
    \begin{tabular}{|c|c|c|c|c|}
    \hline
                &  {\bf Block $\bx$,$\by$:} {\small $M=m,~N=n$} & {\bf Block $\bx$:} {$M=m,~N=1$} & {\bf Block $\by$:} {$M=1,~N=n$}  & {\bf Full:} {$M=N=1$} \\ \hline
         $\delta=0$ & $\cO(mn{(\cC_p+\cC_d)}/\varepsilon)$ & $\cO(mn{(\cC_p+n\cC_d)}/\varepsilon)$ & $\cO(mn{(m\cC_p+\cC_d)}/\varepsilon)$ & {$\cO((m\cC_p+n\cC_d)\mathfrak{Q}/\varepsilon)$} \\ \hline
         $\delta>0$ &$\eh{\tilde\cO(}\sp{(mn+\mathfrak{R}\delta^2)^2}{(\cC_p+\cC_d)}/\varepsilon^2)$ & $\tilde\cO({(n+\delta^2)^2m^2{(\cC_p+n\cC_d)}}/\varepsilon^2)$ & $\tilde\cO(\sp{(m+\delta^2)^2n^2{(m\cC_p+\cC_d)}}/\varepsilon^2)$ & {$\tilde\cO({(m\cC_p+n\cC_d)(\mathfrak{Q}+\delta^2)^2}/\varepsilon^2)$}  
         \\\hline
    \end{tabular}}
    \label{tab:oracle-comp}
    \vspace*{2mm}
    \begin{flushleft}
     {\footnotesize
     \textbf{Note:} {Here the primal $(m)$ and dual $(n)$ dimensions are fixed, and \rev{$\delta^2$ denotes the variance bound for each coordinate of the stochastic oracle, i.e.,
     the bound for $\grad_\bx\Phi(\cdot,\cdot;\xi^x)$ and $\grad_\by\Phi(\cdot,\cdot;\xi^y)$ are $m\delta^2$ and $n\delta^2$.} The table shows the complexity bounds for $4$ different blocking strategies on the same problem 
     \rev{depending on how primal and dual block numbers, $M$ and $N$, are chosen.}}}
    \end{flushleft}
    \vspace*{0mm}
\end{table}
\sp{For the case $\delta=0$, \cite{hamedani2021primal} and \cite{hamedani23a} are special cases of RB-PDA corresponding to (\textbf{Full}) and ({\bf Block $\bx$}) strategies. Indeed, the computational complexity bounds in \citet[Theorem 2.2]{hamedani2021primal} and \citet[Theorem 2]{hamedani23a} are $\cO(\mathfrak{Q}^2/\varepsilon)$ and $\cO(mn^2/\varepsilon)$, respectively.}
\subsection{Contributions} 
To efficiently handle 
large-scale SP problems, we {propose} 
a randomized block-coordinate primal-dual algorithm (RB-PDA) when both primal and dual \sa{variables} are partitioned into blocks -- this setting subsumes \sa{using a single block for primal, \nsa{i.e., $M=1$,} and/or dual} variables, \nsa{i.e., $N=1$,} as special cases. 
\saa{For problems with high-dimensional variables, it may be \emph{infeasible} in practice to compute 
$\grad_\bx\Phi(\bx,\by;\xi)$ and $\grad_\by\Phi(\bx,\by;\xi)$ for any sample $\xi$; hence,} primal and/or dual variables 
{can} be partitioned into blocks 
{in case the computational effort for evaluating block partial gradients is significantly cheaper when compared to \sa{using} full \saa{partial} gradients}, \rev{i.e., when the problem has a coordinate-friendly structure as in Assumption~\ref{assum:comp}}. \rev{Let $Z^*$ be the set of saddle point of \eqref{eq:original-problem}.}\\
$(i)$ Our scheme (RB-PDA) can 
\sa{deal} with the objective functions that are \emph{neither} bilinear \emph{nor} {separable}. \\
$(ii)$ In {the} deterministic setting, we show that the iterate sequence generated by 
{RB-PDA} \saa{converges} \sp{almost surely} to a saddle point $\bz^*=(\bx^*,\by^*)$, \rev{i.e., $\mathbb P(\omega\in\Omega:\ \bz^*(\omega)\in Z^*\}=1$,} with {a non-asymptotic} rate of ${\cO}(1/k)$ \sp{in terms of expected gap metric, where the source of randomness is the random blocks chosen within the runtime of RB-PDA, \rev{determining the sample path $\omega\in\Omega$.}} An obvious implication of this rate \sa{is} 
\rev{$\cO(1/\varepsilon)$ computational complexity} to obtain an $\varepsilon$-gap solution \sp{in expectation}.\\
$(iii)$ In {the} stochastic setting, by utilizing a {\em mini-batch} of 
\saa{stochastic partial} gradients in each iteration, 
\rev{it is shown that the complexity for computing an $\varepsilon$-gap solution \sp{in expectation} 
is $\widetilde\cO(1/\varepsilon^2)$.}\\ 
$(iv)$ We investigate the computational complexity of 
\sp{RB-PDA, compare it with {the} existing results, and study 
how the number of blocks affects the $\cO(1)$} constant of the \saa{complexity} result. In particular, we show that partitioning both primal and dual variables into \rev{$M=m$ and $N=n$ blocks, respectively, leads to computational complexities of $\cO(mn(\cC_p+\cC_d)/\varepsilon)$ and $\tilde\cO((mn+\delta^2)^2(\cC_p+\cC_d)/\varepsilon^2)$} for deterministic and stochastic settings, respectively, while the 
best-known computational complexity bounds for the existing methods with full-block updates, e.g., see \citep{hamedani2021primal,zhao2022accelerated}, are \rev{$\cO(\mathfrak{Q}(m\cC_p+n\cC_d)/\varepsilon)$} and \rev{$\cO((\mathfrak{Q}+\delta^2)^2(m\cC_p+n\cC_d)/\varepsilon^2)$,} 
respectively, where $\mathfrak{Q}\triangleq \max\{m,n\}$ and \rev{$\delta^2$ is a variance bound on the unbiased stochastic partial gradient for each coordinate}. 

Notably, our complexity results {for {\bf (Block $\bx$, $\by$)} {with $M=m$ and $N=n$} strategy} significantly improve these bounds for the deterministic setting when \sp{$m\gg n\gg 1$ or $n\gg m\gg 1$}, and for the stochastic setting when either $m=n^a$ for $n\gg 1$, or $n=m^a$ for $m\gg 1$ for some $a\geq 2$. An alternative scenario in the stochastic setting for which {\bf (Block $\bx$, $\by$)} with $M=m$ and $N=n$ strategy does well is when $M=m\gg 1$ and $N=1$ in case $n$ is small, or when $N=n\gg 1$ and $M=1$ in case $m$ is small.\\ 
$(v)$ The block-specific 
step sizes are reliant on the \textit{average} of coordinate-wise Lipschitz constants, allowing for larger steps (leading to improved empirical behavior) when compared to the step sizes for full gradient methods, which are determined by the largest coordinate-wise Lipschitz constant \rev{(see Remark~\ref{rem:long-steps})}.

Before proceeding to the main content of the paper, we would like to mention a recent work by \cite{zhang2023primal} addressing the SP problem \eqref{eq:original-problem} with affine constraints. In this study, the authors proposed multi-block ADMM-based methods in which one block of primal and dual variables are updated at each iteration. In the deterministic setting, when the coupling function is convex-concave and smooth, an iteration complexity of $\cO(1/\varepsilon)$ is achieved. While the focus of their work is to handle the affine constraint, they did not explore the block-specific step-size; hence, the computational complexity of their methods is $\cO(\mathfrak{Q}^2/\varepsilon)$.

{\bf Organization of paper.} 
In Section~\ref{sec:prelim}, we present essential definitions and preliminary results. Section~\ref{sec:method} introduces our assumptions and the proposed RB-PDA algorithm, followed by the main theoretical results in Section~\ref{sec:main}. The convergence analysis of RB-PDA is provided in Section~\ref{sec:convergence-analysis}, based on the fundamental results established in Section~\ref{sec:support}. Finally, numerical experiments demonstrating the practical effectiveness of RB-PDA on robust binary classification are reported in Section~\ref{numerics}, and concluding remarks are given in Section~\ref{sec:con}. All proofs are deferred to the appendix.


{\bf Notation.} Given $i^*\in\cM=\{1,\ldots,M\}$, define $\mathds{1}_{\{i^*\}}:\cM\to \{0,1\}$ such that $\mathds{1}_{\{i^*\}}(i)=1$ if $i=i^*$, and is equal to $0$ otherwise. 
Let $\mathbb{S}^n_{++}$ ($\mathbb{S}^n_+$) be the set of $n\times n$ symmetric positive (semi-) definite matrices, and $\bI_n$ {denotes} the $n\times n$ identity matrix. Given a normed vector space $(\cX,\norm{\cdot}_{\cX})$, $\cX^*$ denotes the dual space with a dual norm $\norm{\cdot}_{\cX^*}$ such that $\norm{r}_{\cX^*}\triangleq\max_{x\in\cX}\{\fprod{r,x}:\ \norm{x}_{\cX}\leq 1\}$. $\mE[\cdot]$ denotes the expectation operation. 
\sa{$\integers_+$ \eyz{($\integers_{++}$)} denotes the set of \eyz{non-negative (positive)} integers.} We also use w.p. to abbreviate ``\textit{with probability}." 
\section{Preliminaries} \label{sec:prelim}
We first define Bregman distances 
for each $\cX_i$ and $\cY_j$, which 
\saa{generalize} the Euclidean distance.
\begin{defn}
\label{def:bregman}
For $i\in\cM$ and $j\in\cN$, let $\varphi_{\cX_i}:\cX_i\rightarrow\reals$ and $\varphi_{\cY_j}:\cY_j\rightarrow\reals$ be differentiable functions on open sets containing $\dom f_i$ and $\dom h_j$, respectively. Suppose $\varphi_{\cX_i}$ and $\varphi_{\cY_j}$ have closed domains and are 1-\textit{strongly convex} with respect to $\norm{\cdot}_{\cX_i}$ and $\norm{\cdot}_{\cY_j}$, respectively. \sa{For all $i\in\cM$, let $\bD_{\cX_i}:\cX_i\times\cX_i\rightarrow\reals_+$ be the Bregman distance function 
defined as $\bD_{\cX_i}(x,\bar{x})\triangleq \varphi_{\cX_i}(x)-\varphi_{\cX_i}(\bar{x})-\fprod{\grad \varphi_{\cX_i}(\bar{x}),x-\bar{x}}$ \sa{for $x\in\cX_i$ and $\bar{x}\in\dom \varphi_{\cX_i}$}; and for all  $j\in\cN$, $\bD_{\cY_j}:\cY_j\times\cY_j\rightarrow\reals_+$ corresponding to $\varphi_{\cY_j}$
is defined similarly.}
\end{defn}
\begin{defn}\label{def:bregman-sum}
\sa{Given 
positive numbers $\{a_i\}_{i\in\cM}$ and $\{b_j\}_{j\in\cN}$,} let $\bA=\diag{[a_i\id_{\sa{m_i}}]_{i\in\cM}}$ $\in\reals^{m\times m}$ and $\bB=\diag{[b_j\id_{\sa{n_j}}]_{j\in\cN}}\in\reals^{n\times n}$ be \sa{the corresponding block-diagonal} matrices. Let $\bD_\cX^{\bA}(\bx,\bar{\bx})\triangleq \sum_{i\in\cM}a_i\bD_{\cX_i}(x_i,\bar{x}_i)$ and 
$\bD_\cY^{\bB}(\by,\bar{\by})\triangleq \sum_{j\in\cN}b_j\bD_{\cY_j}(y_j,\bar{y}_j)$; for $\bA=\id_m$ and $\bB=\id_n$, we simplify the notation using $\bD_\cX(\bx,\bar{\bx})$ and $\bD_\cY(\by,\bar{\by})$, respectively. Similarly, $\norm{\bx}_{\bA}^2\triangleq\sum_{i\in\cM}a_i \norm{x_i}_{\cX_i}^2$ for $\bx\in\cX$ and $\norm{\by}_{\bB}^2\triangleq\sum_{j\in\cN}b_j \norm{y_j}_{\cY_j}^2$ for $\by\in\cY$; \nsa{and similarly $\norm{\sp{\bx'}}_{*,\bA}^2\triangleq\sum_{i\in\cM}a_i \norm{\sp{x_i'}}_{\cX_i^*}^2$ for $\sp{\bx'}\in\cX^*$ and $\norm{\by'}_{*,\bB}^2\triangleq\sum_{j\in\cN}b_j \norm{y_j'}_{\cY_j^*}^2$ for $\by'\in\cY^*$.}
\end{defn}

\begin{defn}\label{def:func-breg}
{Define $U_i\in\reals^{m\times m_i}$ for $i\in\cM$ and $V_j\in\reals^{n\times \sa{n_j}}$ for $j\in\cN$ such that $\bI_m=[U_1,\hdots,U_M]$ and $\bI_n=[V_1,\hdots,V_N]$}. Let $f(\bx)\triangleq \tfrac{1}{M}\sum_{i\in\cM}f_i(x_i)$ 
and $h(\bx)\triangleq \tfrac{1}{N}\sum_{j\in\cN}h_j(y_j)$.
\end{defn}

Next, we state some fundamental lemmas that will facilitate the proof of the proposed algorithms.
\begin{lemma}[\cite{Tseng08_1J}]
\label{lem_app:prox} Let $\cX$ be a finite-dimensional normed vector space with norm $\norm{.}_\cX$, $f:\cX\rightarrow\reals\cup\{+\infty\}$ be a closed convex function with convexity modulus $\mu\geq 0$ w.r.t. $\norm{.}_\cX$, and $\bD:\cX\times\cX\rightarrow\reals_+$ be a Bregman distance function corresponding to a strictly convex function $\phi:\cX\rightarrow\reals$ that is differentiable on an open set containing $\dom f$. Given $\bar{x}\in\dom f$ and $t>0$, let
$x^+=\argmin_{x\in\cX} \big\{f(x)+t \bD(x,\bar{x})\big\}$, 
then 
${f(x)+t\bD(x,\bar{x})\geq} f(x^+) + t\bD(x^+,\bar{x})+t \bD(x,x^+)+\frac{\mu}{2}\norm{x-x^+}_\cX^2$ holds for all $x\in\cX$. 
\end{lemma}

\begin{lemma}\label{lem:inner-w}
\nsa{Given diagonal $\cA=\diag([a_i\id_{m_i}]_{i\in\cM})\in\mathbb{S}^m_{++}$ and an arbitrary 
$\{\bdelta^k\}_{k\geq 0}\subset\sp{\cX^*}$,} let $\{\bv^k\}_{k\geq 0}$ be a sequence such that $\bv^0\in\cX$ and $\bv^{k+1}\triangleq\argmin_{\bx\in\cX}$ $\{-\fprod{\bdelta^k,\bx}+\bD_\cX^{\nsa{\cA}}(\bx,\bv^k)\}$. Then for all $k\geq 0$ and $\bx\in\cX$, it holds that $\fprod{\bdelta^k,\bx-\bv^k}\leq \bD_\cX^{\cA}(\bx,\bv^k)-\nsa{\bD_\cX^{\cA}(\bx,\bv^{k+1})}+\frac{1}{2}\norm{\bdelta^k}_{\nsa{*,\cA^{-1}}}^2$.
\end{lemma}
\begin{proof} See Section \ref{sec:proof-prelim} in the Appendix.\begin{proof}

\missing{Below we provide the celebrated super-martingale convergence lemma, which 
 is subsequently 
 used in deriving the almost sure (a.s.) convergence statements in 
 our main results.}
 \begin{lemma}[\cite{Robbins71}]\label{lem:supermartingale}
 {Let $(\Omega, \cF,\mathbb{P})$ be a probability space, and for each $k\geq 0$ suppose $a^k$, $b^k$, {$c^k$, and $\sa{\zeta^k}$} are finite, nonnegative $\cF^k$-measurable random variables where $\{\cF^k\}_{k\geq 0}$ is a sequence sub-$\sigma$-algebras of $\cF$ such that $\cF^k\subset\cF^{k+1}$ for $k\geq 0$. If $\mE[a^{k+1}|\cF^k]\leq (1+\sa{\zeta^k})a^k-b^k+ c^k$ {and $\sum_{k=0}^\infty c^k<\infty$ and $\sum_{k=0}^\infty \sa{\zeta^k}<\infty$}, then $a=\lim_{k\rightarrow \infty}a^k$ exists almost surely, and $\sum_{k=0}^\infty b^k<\infty$.}
 \end{lemma}

\section{Proposed Method}\label{sec:method}
In this section, we state some technical assumptions and propose our method for solving 
\eqref{eq:original-problem}. 
\begin{assumption}
\label{assum:bound}
\sp{The set of SP solutions $Z^*\subset\cZ$ of {the problem in} \eqref{eq:original-problem} is nonempty.}
\end{assumption}
\begin{assumption}
\label{assum:sample}
\saa{There exist $\delta_x,\delta_y\geq 0$ such that for any $\bx,\by$, the SFO $\grad\Phi(\cdot,\cdot;\xi)$ satisfies}
\begin{subequations}\label{eq:variance}
{\small
\begin{align}
&\mE\left[\grad_\bx\Phi(\bx,\by;\xi^x)\right]=\grad_\bx\Phi(\bx,\by), \quad \sp{\mE\left[\norm{\grad_{x_i}\Phi(\bx,\by;\xi^x)-\grad_{x_i}\Phi(\bx,\by)}^2_{\cX_i^*}\right]\leq \delta_x^2},\quad i\in\cM,\\
&\mE\left[\grad_\by\Phi(\bx,\by;\xi^y)\right]=\grad_\by\Phi(\bx,\by), \quad \sp{\mE\left[\norm{\grad_{y_j}\Phi(\bx,\by;\xi^y)-\grad_{y_j}\Phi(\bx,\by)}^2_{\cY_j^*}\right]\leq \delta_y^2},\quad j\in\cN.
\end{align}}%
\end{subequations}
\end{assumption}
\begin{assumption}
\label{assum-lip}
\sa{\rev{Suppose $\{f_i\}_{i\in\cM}$ and $\{h_j\}_{j\in\cN}$ are closed convex functions,} and $\Phi(\bx,\by)$ is convex in $\bx$ and concave in $\by$ on $\dom f\times \dom h$. Moreover, for all $\by\in\dom h$, $\Phi$ is differentiable in $\bx$ on an open set containing $\dom f$, and for all ${\bx}\in\dom f$, $\Phi$ is 
differentiable in $\by$ on an open set containing $\dom h$.} We further assume that 
\begin{enumerate}[label=\bf{(\roman*)}]
\item 
for all $i,\ell\in\cM$ and $j\in\cN$, there exist \sa{$L_{x_i x_{\ell}}\geq 0$} and {$L_{x_i y_j}> 0$}
 such that for any \sa{$v\in\cX_\ell$} and $u\in\cY_j$, 
\begin{subequations}
\label{eq:smooth-grad-x}
{\small
\begin{align}
&\erfan{\mathbb E[\norm{\grad_{x_i} \Phi({\bx}+U_{\sa{\ell}}v,\by;\xi^x)-\grad_{x_i} \Phi({\bx},\by;\xi^x)}_{\cX_i^*}^2]\leq L_{x_i x_{\ell}}^2\norm{v}_{\cX_{\sa{\ell}}}^2},\quad \forall~(\bx,\by)\in\dom f\times \dom h,
\label{eq:Lxx}\\
&\erfan{\mathbb E[\norm{\grad_{x_i} \Phi({\bx},\by+V_j u;\xi^x)-\grad_{x_i} \Phi({\bx},\by;\xi^x)}_{\cX_i^*}^2]\leq {L_{x_i y_j}^2}\norm{u}_{\cY_j}^2}, \quad \forall~(\bx,\by)\in\dom f\times \dom h,
\label{eq:Lxy}
\end{align}}%
\end{subequations}
\item for all $i\in\cM$ and $j,\ell \in\cN$, there exist \sa{$L_{y_j y_\ell}\geq 0$} and $L_{y_j x_i}> 0$, such that for any $v\in\cX_i$ and \sa{$u\in\cY_\ell$}, 
\begin{subequations}
\label{eq:smooth-grad-y}
{\small
\begin{align}
&\erfan{\mathbb E[\norm{\grad_{y_j} \Phi({\bx},\by+\sa{V_\ell} u;\xi^y)-\grad_{y_j} \Phi({\bx},\by;\xi^y)}^2_{\cY_j^*}]\leq L_{y_j y_\ell}^2\norm{u}_{\sa{\cY_\ell}}^2},\quad \forall~(\bx,\by)\in\dom f\times \dom h,
\label{eq:Lyy}\\
&\erfan{\mathbb E[\norm{\grad_{y_j} \Phi({\bx}+U_iv,\by;\xi^y)-\grad_{y_j} \Phi({\bx},\by;\xi^y)}_{\cY_j^*}^2]\leq {L_{y_j x_i}^2}\norm{v}_{\cX_i}^2},\quad \forall~(\bx,\by)\in\dom f\times \dom h,
\label{eq:Lyx}
\end{align}}%
\end{subequations}
\item for $i\in\cM$, $\argmin_{x_i \in \cX_i} \big\{ tf_i(x_i)+\fprod{r,x_i}+\bD_{\cX_i}(x_i,\bar{x}_i) \big\}$ \us{may be efficiently computed} for any $\bar{x}_i\in\dom f_i$, \sa{$r\in\cX_i^*$}, and $t>0$. Similarly, for $j\in\cN$, $\argmin_{y_j\in\cY_j}\{th_j(y_j)+\fprod{s,y_j}+\bD_{\cY_j}(y_j,\bar{y}_j)\}$ 
\saa{is also efficient} for any $\bar{y}_j\in\dom h_j$, $s_j\in\cY_j^*$ and $t>0$.
\end{enumerate}
\end{assumption}
\begin{remark}
    \sp{Due to Jensen's inequality, one can lower bound 
    \eqref{eq:smooth-grad-x} and \eqref{eq:smooth-grad-y} by switching $\mathbb E[\cdot]$ and $\norm{\cdot}^2$; hence, from Assumption~\ref{assum:sample}, 
    deterministic $\grad \Phi(\cdot,\cdot)$ is smooth with the same block Lipschitz constants.}
\end{remark}

\begin{remark}
\label{rem:L-constants}
\rev{According to \citet[Lemma~2]{nesterov2012efficiency}, {for all $i\in\cM$, if \eqref{eq:Lxx} holds only for $\ell=i$, and for all $j \in \cN$ if \eqref{eq:Lyy} holds only for $\ell = j$,} 
then the partial gradients $\grad_\bx\Phi(\cdot,\by)$ and $\grad_\by\Phi(\bx,\cdot)$ are Lipschitz continuous, with constants $L_{\bx\bx}$ and $L_{\by\by}$, respectively.}
Nonetheless, we introduce {$L_{x_i x_\ell}$ and $L_{y_j y_\ell}$} 
in Assumption~\ref{assum-lip} 
as they are generally much smaller than $L_{\bx\bx}$ and $L_{\by\by}$, respectively.
 \sp{Indeed, consider $\Phi(\bx,\by)=\frac{1}{2}\bx^\top A_1\bx+\by^\top A_2\bx-\frac{1}{2}\by^\top A_3\by$ for $A_1=\mathbf{1}_{m\times m}$, $A_2=\mathbf{1}_{n\times m}$ and $A_3=\mathbf{1}_{n\times n}$ are matrices with all entries being $1$. For this case, $L_{x_i x_\ell}=L_{y_j y_\ell}=L_{x_i y_j}=L_{y_j x_i}=1$ for all $i,j,\ell$; however, $L_{\bx\bx}=m$, $L_{\by\by}=n$, and $L_{\by\bx}=L_{\bx\by}=\sqrt{mn}$, where $L_{\by\bx}$ and $L_{\bx\by}$ are the Lipschitz constants of $\grad_\by\Phi(\cdot,\by)$ and $\grad_\bx\Phi(\bx,\cdot)$.}
\end{remark}

{Based on} \eyh{Assumption \ref{assum-lip}, for $i\in\cM$ and $j\in\cN$}, we define the following constants: 
 \begin{subequations}
 {\small
 \begin{align*}
 &{L_{\by x_i}}\triangleq\Big({\tfrac{1}{N}\sum_{j\in\cN}L_{y_j x_i}^2}\Big)^{1/2},\quad  {C_{x_i}\triangleq \Big(\tfrac{1}{M}\sum_{\ell\in\cM}L_{x_\ell x_i}^2\Big)^{1/2}},\quad
 {L_{\bx y_j}}\triangleq\Big(\tfrac{1}{M}\sum_{i\in\cM}L_{x_iy_j }^2\Big)^{1/2},\quad {C_{y_j}\triangleq \Big(\tfrac{1}{N}\sum_{\ell\in\cN}L_{y_\ell y_j}^2\Big)^{1/2}}, 
 \end{align*}}%
 \end{subequations}
 Moreover, \vspace*{0mm}
 \begin{equation}
 \label{def L}
 {\small
 \begin{aligned}
 &\bL_{\bx\bx}\triangleq \diag([L_{x_ix_i}\id_{m_i}]_{i\in\cM}),\quad {\bL_{\by \bx}}\triangleq \diag([L_{\by x_i}\id_{m_i}]_{i\in\cM}),\quad \bL_{\bx \by}\triangleq \diag([L_{\bx y_j}\id_{n_j}]_{j\in\cN}),\\
 &\bL_{\by\by}\triangleq \diag([L_{y_jy_j}\id_{n_j}]_{j\in\cN}),\quad {\bC_\bx\triangleq \diag([C_{x_i}\id_{m_i}]_{i\in\cM})},\quad {\bC_\by\triangleq \diag([C_{y_j}\id_{n_j}]_{j\in\cN})}.
 \end{aligned}}%
 \vspace*{-1mm}
 \end{equation}
 {
 If 
 $\Phi$ is {\it separable} in $\bx$-variable, i.e., $\Phi(\bx,\by)=
 \sum_{i\in\cM}\Phi_{i}(x_i,\by)$, then \missing{$C_{x_i}=L_{x_ix_i}/\sqrt{M}$}; \sa{
 separability in $\by$-variable, i.e., $\Phi(\bx,\by)=\sum_{j\in\cN}\Phi_{j}(\bx,y_j)$, implies $C_{y_j}=L_{y_jy_j}/\sqrt{N}$.} In general, the constants $L_{\by x_i}, C_{x_i}, L_{\bx y_j}, C_{y_j}$
 are always less than the global Lipschitz constant for $\grad \Phi$. For instance, $\Phi(\bx,\by)=\sum_{\ell=1}^N y_\ell\bx^\top Q_\ell \bx$ implies 
 {$C_{x_i}=\sum_{j=1}^N \abs{y_j}\sqrt{\frac{1}{M}\sum_{\ell\in\cM}\|Q_{j}(\ell,i)\|^2}$}
 while a global Lipschitz \us{constant} of $\nabla_{\bx}\Phi$ is $\sum_{j=1}^N \abs{y_j}\|Q_{j}\|\nsa{\gg C_{x_i}}$ \nsa{when $M\gg 1$}, \rev{where 
$Q_j(\ell,i)\in\reals^{m_\ell\times m_i}$ denotes the submatrix corresponding to $\ell$-th block-row and $i$-th block-column of the block matrix $Q_j\in\reals^{m\times m}$, and $\|\cdot\|$ denotes the spectral-norm.}}

\sp{Next, we {define the notation} 
{related to} {the} {random} sampling of the {partial} gradients.} 
\begin{defn}
\label{def:error}
For any $(\bx,\by)\in\cX\times \cY$, $u\in\mZ_{++}$, and mini-batch sample $\cB=\{\saa{\xi^x_\ell}\}_{\ell=1}^u$, we define $\grad_\bx\Phi_\cB(\bx,\by)\triangleq \frac{1}{u}\sum_{\ell=1}^u\grad_\bx\Phi(\bx,\by;\xi_\ell^x)$. Similarly, for a mini-batch sample $\cS=\{\saa{\xi^y_\ell}\}_{\ell=1}^v$ and $v\in\mZ_{++}$ we define $\grad_\by\Phi_\cS(\bx,\by)\triangleq \frac{1}{v}\sum_{\ell=1}^v\grad_\by\Phi(\bx,\by;\xi^y_\ell)$.
\end{defn}

We propose 
the Randomized Block coordinate Primal-Dual Algorithm (RB-PDA), displayed in Algorithm \ref{alg:RB-PDA}, to solve 
\eqref{eq:original-problem}. In this algorithm, we use an ascent Bregman-proximal step with an additive momentum to update variable $\by$ followed by a descent Bregman-proximal gradient step with an additive momentum to update variable $\bx$. In fact, we use two momentum terms based on the gradients of the objective function, which are crucial to guarantee the convergence rates proposed for this method.
At each iteration $k\in\integers_+$, we randomly select a single dual and a single primal block \us{coordinate}, \us{denoted by} $j_k\in\cN$ and $i_k\in\cM$, \us{respectively}. Next, we evaluate the partial gradient map $\grad_{y_{j_k}}\Phi$ at $(\bx^k,\by^k)$ and $(\bx^{k-1},\by^{k-1})$ and
$\grad_{x_{i_k}}\Phi$ at $(\bx^k,\by^{k})$, $(\bx^k,\by^{k+1})$ and $(\bx^{k-1},\by^{k-1})$ for the deterministic setting. \rev{In the stochastic setting, we replace the gradients with their mini-batch sample partial gradients $\grad_{y_{j_k}}\Phi_{\cB^k}$ and $\grad_{x_{i_k}}\Phi_{\cS^k}$ using mini-batch samples $\cB^k$ and $\cS^k$ at the $k$-th iteration.}

\begin{algorithm}[htb]
\caption{RB-PDA ({$\{{[\tau_i^k]_{i\in\cM}},[\sigma_j^k]_{j\in\cN},\theta^k\}_{k\geq 0}$}, $(\bx_0,\by_0)\in\cX\times\cY$}
 \label{alg:RB-PDA}
 {\small
\begin{algorithmic}[1]
   \STATE {$(\bx_{-1},\by_{-1})\gets(\bx_0,\by_0)$},~ 
   \FOR{$k\geq 0$}
   \STATE {Sample} a coordinate $j_k\in\cN$ uniformly at {random};
   \STATE \ey{Randomly generate \rev{i.i.d.} mini-batch samples $\sp{\cB^k=\{\xi_{\ell}^{y}\}_{\ell=1}^{u^k}}$ with \sp{$\abs{\cB^k}=u^k$};}
			\STATE \label{algeq:s} \sp{$s_{j_k}^k\gets N\grad_{y_{j_k}}\Phi_{\cB^k}(\bx^k,\by^k)+N M \theta^k\left(\grad_{y_{j_k}}\Phi_{\cB^k}(\bx^k,\by^k)-\grad_{y_{j_k}}\Phi_{\cB^k}(\bx^{k-1},\by^{k-1})\right)$;}
   		\STATE \label{algeq:y} 
   		{$\by^{k+1}\gets\by^k$,}\quad
   		{$y_{j_k}^{k+1}\gets\argmin\limits_{y\in\cY_{j_k}}\left[h_{j_k}(y)-\fprod{s_{j_k}^k,y}+\frac{1}{\sigma_{j_k}^k}\bD_{\cY_{j_k}}(y,y_{j_k}^k)\right]$};
   		\STATE {Sample} $i_k\in\cM$ uniformly at {random}; 
		\STATE  \sp{Randomly generate \rev{i.i.d.} mini-batch samples $\cS^k=\{\xi_{\ell}^{x}\}_{\ell=1}^{v^k}$ with $\abs{\cS^k}=v^k$;}
   		\STATE  \sp{$r_{i_k}^k\gets M\grad_{x_{i_k}}\Phi_{\cS^k}(\bx^k,\by^{k+1})+(N-1)M\theta^k(\grad_{x_{i_k}}\Phi_{\cS^k}(\bx^k,\by^k)-\grad_{x_{i_k}}\Phi_{\cS^k}(\bx^{k-1},\by^{k-1}))$};\label{algeq:rik}
   \STATE \label{algeq:x} 
   {$\bx^{k+1}\gets\bx^k$,}\quad
   {$x_{i_k}^{k+1}\gets\argmin\limits_{x\in\cX_{i_k}}\left[f_{i_k}(x)+\fprod{r_{i_k}^k,x} +\frac{1}{\tau_{i_k}^k}\bD_{\cX_{i_k}}(x,x_{i_k}^k)\right]$};
   \ENDFOR
\end{algorithmic}}%
\end{algorithm}
\subsection{Main Results}
\label{sec:main}
\rev{In this subsection, we first formally introduce the gap function and then state the convergence guarantees of the proposed algorithm in stochastic and deterministic settings.}
\rev{
\begin{defn}
    Consider the SP problem in \eqref{eq:original-problem}. Given some bounded set $Z
    \subset\cZ$, $G_Z:\cZ\to\reals_+$ defined in \eqref{s-gap} is called the \textit{restricted gap function}. Moreover, $G:\cZ\to\reals_+$ such that $G(\cdot)=G_Z(\cdot)$ for $Z=\dom f\times\dom h$ is called the \textit{standard gap function}.
\end{defn}}

\rev{Next, we discuss the properties of such restricted gap function. 
\begin{lemma}
\label{lem:restricted_gap}
    Let $\oX\subseteq \cX$ and $\oY\subset\cY$ be some bounded sets, and define $\oZ\triangleq \oX\times \oY$. Let $\bar{Z}^*\subset\cZ$ be the set of saddle points of the following restricted SP problem:
    \begin{equation}\label{eq:restricted-problem}
    {\rm (RSP)}:\quad \min_{\bx\in\oX}\max_{\by\in\oY}~\cL(\bx,\by).
\end{equation}
    Suppose $\overline Z$ contains a saddle point of \eqref{eq:original-problem}, i.e., $\oZ\cap Z^*\neq\emptyset$. Then $G_{\oZ}(\bz)\geq 0$ for all
    $\bz\in\cZ=\cX\times\cY$. Moreover, if $G_{\oZ}(\bar{\bz})=0$ for some $\bar\bz\in \oZ$, then $\bar\bz$ is a saddle-point of (RSP), i.e., $\bar\bz\in \bar Z^*$. Finally, it holds that $\oZ\cap Z^*\subset \bar Z^*$; conversely, if there exists $\bar\bz\in\bar Z^*$ such that $\bar\bz\in \operatorname{int}\oZ$, then $\bar\bz\in Z^*$ as well.
\end{lemma}}
\begin{proof}
\rev{Let $\bz^*=(\bx^*,\by^*)\in \oZ\cap Z^*$. Given any $\bar\bz\in\cZ$, following the definition of the restricted gap function, we have $G_{\oZ}(\bar \bz)=
  \sup_{\bz=(\bx,\by)\in \oZ}\{\mathcal L(\bar \bx,\by)-\mathcal L(\bx,\bar \by)\}\geq \mathcal L(\bar \bx,\by^*)-\mathcal L(\bx^*,\bar \by)\geq 0$ where the final inequality is due to $\bz^*\in Z^*$.} 
  
  \rev{Now suppose $G_{\oZ}(\bar \bx,\bar \by)=0$ for some $(\bar \bx,\bar \by)\in \oZ$, which implies that $\sup_{\by\in \oY}\mathcal L(\bar \bx,\by) = \inf_{\bx\in \oX}\mathcal L(\bx,\bar \by)$. It also trivially holds that $\sup_{\by\in\oY}\mathcal L(\bar \bx,\by)\ge \mathcal L(\bar \bx,\bar \by)
\ge \inf_{\bx\in\oX}\mathcal L(\bx,\bar \by)$; hence, above equality implies that 
\begin{equation*}
  \mathcal L(\bar \bx,\by) \leq 
  \mathcal L(\bar \bx,\bar \by) \leq
  \mathcal L(\bx,\bar \by),\quad \forall (\bx,\by)\in\oX\times \oY,
\end{equation*}
which implies that $\bar\bz=(\bar \bx,\bar\by)$ is an SP solution of (RSP), i.e., $\bar\bz\in\bar Z^*$.}

\rev{Next, we argue that $\oZ\cap Z^*\subset\bar Z^*$. Recall the standard gap function definition $G(\bar \bx,\bar \by)= \sup_{\bz\in\mathcal Z}\{\mathcal L(\bar \bx,\by)-\mathcal L(\bx,\bar \by)\}$. Clearly $G(\bar \bx,\bar \by)\geq G_{\oZ}(\bar \bx,\bar \by)$ for all $(\bar \bx,\bar \by)\in\mathcal Z$. Let $\bz^*\in \oZ\cap Z^*$. Then we have $G(\bz^*)=0$; thus, $0 = G(\bz^*) \geq G_{\oZ}(\bz^*) \geq 0$,
implying that $G_{\oZ}(\bz^\ast)=0$. Since $\bz^*\in\oZ$, we can conclude that $\bz^*\in\bar Z^*$; consequently, $\oZ \cap Z^\ast \ \subseteq\ \bar Z^\ast$. Moreover, for the converse argument, let $(\bar \bx,\bar \by)\in\bar Z^*$ such that $(\bar \bx,\bar \by)\in \operatorname{int}\oX \times \operatorname{int}\oY$. Since $\mathcal L(\cdot,\by)$ and $-\mathcal L(\bx,\cdot)$ are closed convex functions for every $(\bx,\by)\in\cZ$, it must hold that $(\bar \bx,\bar \by)\in Z^*$.}
\end{proof}
\rev{Indeed, Lemma~\ref{lem:restricted_gap} shows that if we can identify a bounded set $\oZ$ such that $\oZ\cap Z^*\neq \emptyset$, then one can consider generating a sequence $\{\bz^k\}_k\subset\oZ$ such that $\bz^k\to\bar\bz\in\oZ$ with $G_{\oZ}(\bar\bz)=0$. Since $\oZ$ is bounded, $\{G_{\oZ}(\bz^k)\}_k\subset\reals_+$, i.e., a real-valued sequence, such that $G_{\oZ}(\bz^k)\to 0$ as $k\to\infty$; on the other hand, when the domain of the SP problem in \eqref{eq:original-problem} is unbounded, $G(\bz^k)$ is not necessarily finite for all $k\geq 0$, and one may not be able to quantify the rate of convergence to $\bar\bz\in\bar Z^*$, i.e., to a saddle point of the restricted problem. That being said, there still remain certain pitfalls concerning the definition of the restricted gap function, which we discuss in the remark below.
\begin{remark}
\label{rem:pathological-cases}
There are 
some issues with the use of the restricted gap function as a convergence metric.
    \begin{enumerate}
    \item It is possible that $\bar Z^*\setminus Z^*\neq\emptyset$. In particular, for some $(\bar \bx,\bar \by)\notin\mathrm{int} \oZ$, one may have $(\bar \bx,\bar \by)\in\bar Z^*$ and yet $(\bar \bx,\bar \by)\notin Z^*$. For example, consider $\min_{\bx\in\reals}\max_{\by\in\reals} \bx\by$ with $\cX=\reals$ and $\cY=\reals$, which has a unique saddle point, i.e., $Z^*=\{(0,0)\}$. For $\oX=[0,1]$ and $\oY=[0,1]$, the point $(\bar \bx,\bar \by)=(0,1)$ satisfies $(\bar \bx,\bar \by)\in \bar Z^*$ and $(\bar \bx,\bar \by)\notin Z^*$.
    \item 
    It is possible that for some $\bar\bz=(\bar \bx,\bar\by)\notin \oZ$, one has $G_{\oZ}(\bar \bz)=0$ while $\bar \bz\notin Z^*$. Again consider $\min_{\bx\in\reals}\max_{\by\in\reals} \bx\by$; for $\oZ=\{(0,0)\}$ and $(\bar\bx,\bar\by)=(1,1)$, it holds that $G_{\oZ}(1,1)=0$ but $(1,1)\notin  Z^*$.
\end{enumerate}
\end{remark}
Now combining the results in Lemma~\ref{lem:restricted_gap} with our observations in Remark~\ref{rem:pathological-cases}, we can conclude with the following key result.
\begin{lemma}\label{lem:restricted-gap-solution}
    Consider the SP problem in \eqref{eq:original-problem}. Let $\{\bz^k\}_{k\geq 0}\subset\cZ$ be a deterministic sequence such that
    $\{\bz^k\}_{k\geq 0}\subset\oZ$ for some compact set $\oZ\subset\dom f\times\dom h$ such that $\oZ\cap Z^*\neq \emptyset$. Let $Z$ be an open set that contains $\oZ$. Suppose $\bz^k\to\bar\bz\in\bar Z$ such that $G_Z(\bz^k)\to 0$. Then, it holds that $\bar\bz\in Z^*$.
\end{lemma}}

\rev{The above discussion does not involve any stochastic component; however, there are two sources of randomness within RB-PDA while generating $\{\bz^k\}$, i.e., one is due to random sampling of block-coordinates and the other is due to use of stochastic gradients to estimate $\grad\Phi$. Still from the arguments above, it is clear that for a restricted gap function $G_Z(\cdot)$ to be a well-defined convergence metric, one needs to ensure that the sequence generated by the algorithm remains bounded so that $Z$ can be chosen appropriately. While this boundedness 
trivially holds when $\dom f\times \dom h$ is compact, when the problem domain $\dom f\times \dom h$ is unbounded, even almost sure convergence of the iterates does not guarantee uniform boundedness of the stochastic sequence, i.e., the norm bound $\Delta(\omega)$ on the sequence $\{\bz^k(\omega)\}_{k\geq 0}$ may possibly depend on the sample path $\omega\in\Omega$, and one may potentially have $\sup_{\omega\in\Omega}\Delta(\omega)=+\infty$. Instead, we show that, with high probability, the generated sequence remains within a bounded set $Z$ --indeed, for any given $p\in(0,1)$, we argue that $\mathbb{P}(\sup_{k\geq 0}\norm{\bz^k}\leq \sqrt{\Delta_p})\geq 1-p$ for some $\Delta_p=\cO(1/p)$. Consequently, we demonstrate that the convergence rate can be expressed in terms of the conditional expectation of the restricted gap function over the iterates in $Z$, providing a well-defined performance metric, e.g., see Corollary \ref{cor:gap-result-stoch}.}


{To provide rate statements, 
we adopt a \textit{weighted restricted 
gap} function $\sp{\cG_Z}:\cZ\times\reals_+\to\reals$ as follows:
\begin{align}\label{w-gap}
    \sp{\cG_Z}(\bar{\bz},T)\triangleq\sup_{\nsa{\bz\in Z}}\big\{(1+\tfrac{M-1}{T})(\cL(\bar{\bx},\by)-\cL(\bx,\by))+(1+\tfrac{N-1}{T})(\nsa{\cL(\bx,\by)-\cL(\bx,\bar{\by})})\big\},
\end{align}
where \rev{$Z\subset \cZ$ is some bounded set that will be defined later on.} We analyze how fast \nsa{$\cG_Z(\bar{\bz}^K,T_K)$ diminishes to zero as $K\to\infty$ for some particular $\{T_K\}$ such that $T_K\to\infty$ as $K\to\infty$}, where $\{\bar{\bz}^k\}_{k\geq 1}\sp{\subset Z}$ is an averaged RB-PDA iterate sequence. Before proceeding to \nsa{our main results, we first discuss how to translate the results obtained in terms of \sp{$\cG_Z$} to guarantees in} 
$G_Z$. 
Clearly 
$\cG_Z$ defined in \eqref{w-gap} reduces to the restricted gap $G_Z$ in~\eqref{s-gap} when $M=N$, i.e., \sp{$\cG_Z(\bar \bz,T)=(1+\frac{M}{T})G_Z(\bar \bz)$.} Therefore, whenever $M\neq N$ one can consider implementing RB-PDA 
with dummy blocks \nsa{so that 
$M=N$; moreover, whenever a dummy block is sampled, one does not compute anything, i.e.,} no update occurs on the added blocks. 
\begin{remark}\label{rem:gap}
\nsa{
\sp{For the case $M\neq N$, even 
without 
dummy blocks,} a rate result on the weighted gap function in~\eqref{w-gap} immediately implies a rate result in terms of the restricted gap in~\eqref{s-gap} through} 
$\cG_Z(\bar \bz,T)\geq \left(1+\frac{\max\{M,N\}-1}{T}\right)G_Z(\bar \bz)-\frac{c}{T}\cH$,
where $c\triangleq\max\{M,N\}-\min\{M,N\}$ and $\cH\triangleq \sup_{\bz\in\nsa{Z}}\{\cL(\bar\bx,\by)-\cL(\bx,\by)\}<+\infty$ if $N\geq M$, and $\cH\triangleq \sup_{\bz\in\nsa{Z}}\{\cL(\bx,\by)-\cL(\bx,\bar\by)\}<+\infty$ otherwise. 
\end{remark}}

Next, we state our main results establishing the complexity bounds for RB-PDA. We first derive the iteration and sample complexity 
for the stochastic setting, i.e., \sp{$\bar\delta\triangleq\max\{\delta_x,\delta_y\}>0$}, where we use fixed batch sizes $u^k=u$ and $v^k=v$ at each step $k\geq 0$ of RB-PDA, for some $u,v\in\integers_+$, to estimate 
the block partial gradients using a stochastic first-order oracle.
Later, in the deterministic setting, \sp{where we assume that the first-order oracle for $\Phi$ is deterministic, i.e., $\bar\delta=0$,}
we analyze the iteration complexity of 
RB-PDA that uses randomly selected block partial gradients of $\Phi$, i.e., $\grad_{x_{i_k}}\Phi$ and $\grad_{y_{j_k}}\Phi$, in each iteration. 
\subsubsection{Stochastic Setting} 
\eyh{Here we consider the situation where} $\grad_{x_i}\Phi$ and $\grad_{y_j}\Phi$ are estimated using stochastic gradients.  
This setting arises in different problems, such as online optimization, where the problem information becomes available over time. 
We show an ergodic convergence rate of 
$\cO(\log(K)/\sqrt{K})$ \sp{assuming that we use constant (mini) batch-sizes \rev{such that $u=\cO(1)$ and $v=\cO(1)$ in terms of dependence on the tolerance $\varepsilon$.}} To control the error of using a mini-batch sample at each iteration, 
the step sizes are chosen to be diminishing. 

\begin{theorem}\label{thm:bounded-stochastic}
Let $\{\bx^k,\by^k\}_{k\geq 0}$ be the sequence generated by 
RB-PDA, {initialized} from arbitrary vectors $\bx^0\in\cX$ and $\by^0\in\cY$. Suppose Assumptions~\ref{assum:bound}-\ref{assum-lip} hold. 
For all $k\geq 0$, let 
$|\cB^k|=u$ and $|\cS^k|=v$ for any $\integers_+\ni u,v\geq 1$, and let
$\tau_{i_k}^k=\tau_i^k$ w.p. $\frac{1}{M}$ for $i\in \cM$ and $\sigma^k_{j_k}=\sigma^k_j$ w.p. $\frac{1}{N}$ for $j\in \cN$, where
\vspace*{0mm}
{\small
\begin{equation}
\label{eq:step-k}
\begin{aligned}
\tau_i^k&={\tfrac{t^k}{\sp{(1+\varsigma)}M}}\Big(L_{x_ix_i}+(N-1)\Big(2({\gamma_1^{-1}+(1+\tfrac{1}{2M})\gamma_2^{-1}})+\eh{(\gamma_1+\tfrac{\erfan{16(N-1)M}}{\sp{\alpha^0}})}C_{x_i}^2\Big)+\eh{(\lambda_2+\tfrac{\erfan{16MN}}{\sp{\tilde\alpha^0}})}N L_{\by x_i}^2
+\tfrac{{\alpha^0}}{M}\Big)^{-1},\\ 
\sigma_j^k&={\tfrac{t^k}{\sp{(1+\varsigma)}N}}\Big(L_{y_jy_j}+ 2M({\lambda_1^{-1}+\lambda_2^{-1}})+ \eh{(\lambda_1+\tfrac{\erfan{16 MN}}{\sp{\tilde\alpha^0}})}MC_{y_j}^2 +(\erfan{\tfrac{M+1}{\erfan{N}}}\gamma_2+\tfrac{\erfan{16M^2(N-1)}}{\alpha^0 N})\erfan{(N-1)} L^2_{\bx y_j}+\tfrac{{\sp{\tilde\alpha^0}}}{N}\Big)^{-1}, 
\end{aligned}
\vspace*{-1mm}
\end{equation}}%
such that 
\sp{$\gamma_1,\gamma_2,\lambda_1,\lambda_2,\alpha^0,\tilde\alpha^0,\rev{\varsigma}>0$ are arbitrary.} 
Moreover, let \erfan{$t^{k}=(\sqrt{k+1}\log(k+3))^{-1}$ and {$\theta^k=t^{k-1}/t^k$}} for $k\geq 1$ such that $t^0=\theta^0=1$. \rev{Then, for any given $p\in (0,1)$ and arbitrary $\bz^*\in Z^*$, 
$\mathbb P\!\left(\sup_{k\ge0}\{\| \bx^k-\bx^*\|^2_{\bT^0}+\| \by^k-\by^*\|^2_{\bS^0}\}\le \Delta_p\right)\ge 1-p$ 
for $\Delta_p\triangleq \frac{2}{\varsigma p}\Big(B(\bz^*)+\frac{4{M^2}\delta_x^2}{v\alpha^0}+\frac{4{N^2}\delta_y^2}{u\tilde\alpha^0}\Big)$, where $\bT^0\triangleq \diag([1/{\tau^0_i}]_{i\in \cM})$, $\bS^0\triangleq \diag(([1/{\sigma^0_j}]_{j\in\cN})$ and $B(\bz^*)\triangleq\Delta(\bz^*)+ C_0(\bz^*)+\cL(\bx^0,\by^*)-\cL(\bx^*,\by^0)$ for $\Delta(\cdot)$ and $C_0(\cdot)$ defined in~\eqref{eq:delta-function} and \eqref{eq:C0-function}, respectively.}
\rev{Thus, for any $r>0$, the bounded open set 
\begin{align}
\label{eq:Z-set}
    Z\triangleq\{(\bx,\by)\in\cX\times\cY:\ \norm{\bx-\bx^*}^2_{\bT^0}+\norm{\by-\by^*}^2_{\bS^0}<(1+r)\Delta_p\}
\end{align}
satisfies $Z\cap Z^*\neq \emptyset$ and $\{\bz^k\}\subset Z$ with probability at least $1-p$.}
\end{theorem}
\begin{proof}
See Section \ref{sec:proof-bounded-thm} in the Appendix.
\end{proof}

\begin{theorem}\label{thm:single-sample}
Under the premise of Theorem \ref{thm:bounded-stochastic}, let $\cA^0\triangleq \alpha^0\id_m$, $\tilde\cA^0\triangleq \tilde\alpha^0\id_n$, and define
\begin{subequations}
\label{eq:Delta-C-def}
{\small
\begin{align}
     \Delta(\bz)&\triangleq \bD_\cX^{\sp{2\bT^0+(N-1)M\bL_{\bx\bx}+\frac{1}{M}\cA^0}}(\bx,\bx^0)+\bD_\cY^{\sp{2\bS^0+NM\bL_{\by\by}+\frac{1}{N}\tilde\cA^0}}(\by,\by^0),\label{eq:delta-function}\\
     C_0(\bz)&\triangleq (M-1)\big(\cL(\bx^0,\by)-\cL(\bx,\by)\big)+(N-1)\big(\cL(\bx,\by)-\cL(\bx,\by^0)\big),\label{eq:C0-function}\\
     \sp{C_1}(\bz)&\triangleq \max\{\tfrac{M-1}{M} L^2_{\varphi_\cX},\tfrac{N-1}{N}L^2_{\varphi_\cY}\} \big(\Delta(\bz)+ C_0(\bz)+\cL(\bx^0,\by)-\cL(\bx,\by^0)\big)/\sp{\varsigma}\label{eq:C1}\\
      B_\delta&\triangleq \sp{\Big(1+\tfrac{1}{\varsigma}\max\{\tfrac{M-1}{M} L^2_{\varphi_\cX},\tfrac{N-1}{N}L^2_{\varphi_\cY}\}\Big) \Big(\frac{4{M^2}\delta_x^2}{v\alpha^0}+\frac{4{N^2}\delta_y^2}{u\tilde\alpha^0}\Big).}
\end{align}}%
\end{subequations}
\rev{If $\dom f\times\dom h$ is not bounded, then let $Z\subset \cZ$ be the bounded open set 
given in \eqref{eq:Z-set} for any $r>0$, $p\in (0,1)$ and $\bz^*\in Z^*$; otherwise, when $\dom f\times\dom h$ is bounded, let $Z=\dom f\times\dom h$.}
 
{\bf (I)} Let \sp{$T_K\triangleq\sum_{k=0}^{K-1} t^{k}\geq \frac{2\sqrt{K+1}-1}{\log(K+2)}$ for $K\geq 1$.} The averaged iterates $\bar\bz^K=(\bar\bx^K,\bar\by^K)$ defined as 
{\small
\begin{equation}
\label{eq:avg-iterate}
    \begin{aligned}
    \bar{\bx}^K&\triangleq(M-1+T_K)^{-1}\left(\sum_{k=0}^{K-1}t^k\Big[1+(M-1)(1-\tfrac{1}{\theta^{k+1}})\Big]\bx^{k+1}+(M-1)t^K\bx^K\right),\\
    \bar{\by}^K&\triangleq(N-1+T_K)^{-1}\left(\sum_{k=0}^{K-1}t^k\Big[1+(N-1)(1-\tfrac{1}{\theta^{k+1}})\Big]\by^{k+1}+(N-1)t^K\by^K\right)
\end{aligned}
\end{equation}}%
satisfy the following bound \sp{for any $K\geq 1$ and \rev{for any $\bz^*\in Z^*$}:} 
{\small
\begin{equation}\label{bound_c_2}
\mE\left[\cG_Z(\bar \bz^K,T_K)\right] \leq\frac{\bar \Delta
+\sp{C_1(\bz^*)+B_\delta}}{T_K}, \quad\mbox{where}\quad \bar \Delta\triangleq \sup_{\bz\in Z}\{\Delta(\bz)+{C_0(\bz)}\}<+\infty.
\vspace*{0mm}
\end{equation}}%

{\bf (II)} 
\rev{The number of primal and dual oracle calls} 
for achieving 
\sp{
$\mE[G_Z(\bar\bz^K)]\leq \varepsilon$} are 
$\cO\left(\frac{v\bar D^2}{\varepsilon^2}\log(1/\varepsilon)\right)$ and $\cO\left(\frac{u\bar D^2}{\varepsilon^2}\log(1/\varepsilon)\right)$, 
where $\bar D\triangleq \bar{\Delta}+\sp{{C}_1(\bz^*)+B_\delta}+c\cH$,  
and $c,\cH\in (0,+\infty)$ are defined in Remark \ref{rem:gap}.
\end{theorem}
\begin{proof}
See Section \ref{sec:stochastic-proof} in the Appendix.
\end{proof}

\begin{corollary}\label{cor:gap-result-stoch}
\rev{Under the premise of Theorem \ref{thm:bounded-stochastic}, suppose $\dom f\times \dom h$ is not bounded and let $Z\subset \cZ$ be the bounded open set defined in \eqref{eq:Z-set} for any given $r>0$, $p\in (0,1)$ and $\bz^*\in Z^*$. Then, for any $K\geq 1$, $\mE[G_Z(\bar\bz^K)\mid \{\bar\bz^k\}_{k\geq 0}\subset Z]\leq \frac{1}{1-p}\frac{\bar{D}}{T_K}$, where $\bar D$ is defined in Theorem \ref{thm:single-sample}.}
\end{corollary}
\begin{proof} 
\rev{
According to Theorem~\ref{thm:bounded-stochastic}, we have $\mathbb P\!\left(\sup_{k\ge0}\{\| \bx^k-\bx^*\|^2_{\bT^0}+\| \by^k-\by^*\|^2_{\bS^0}\}\le \Delta_p\right)\ge 1-p$.
Hence, for the event $E\triangleq \{\{\bar \bz^k\}_{k\ge 0}\subset Z\}$, we have that $\mathbb P(E)\ge 1-p$ since $\{\bz^k\}_{k\ge 0}\subset Z$ implies that $\{\bar\bz^k\}_{k\ge 0}\subset Z$.
Moreover, $\mathbb E\!\left[G_Z\!\big(\bar{\bz}^{K}\big)\right]
= \mathbb E\!\left[G_Z\!\big(\bar{\bz}^{K}\big)\mid E\right]\mathbb P(E)
  +  \mathbb E\!\left[G_Z\!\big(\bar{\bz}^{K}\big)\mid E^{\mathrm c}\right](1-\mathbb P(E))$.
Since $G_{Z}(\cdot)\geq 0$, the second conditional expectation is nonnegative; thus, 
using Theorem \ref{thm:single-sample}-\textbf{(II)} and Remark~\ref{rem:gap} leads to the desired result:
$\mathbb E\!\left[G_Z\!\big(\bar{\bz}^{K}\big)\mid E\right]
 \leq \mathbb E\!\left[G_Z\!\big(\bar{\bz}^{K}\big)\right]/{\mathbb P(E)}
\ \le\ \frac{1}{1-p}\frac{\bar{D}}{T_K}$.}
\qed

Next, we simplify the convergence rate bound obtained in Theorem \ref{thm:single-sample} by assuming that the domains of \rev{$f,h$} are compact and the coordinate-wise Lipschitz constants 
\nsa{are the same for all coordinates.}
\begin{assumption}
\label{assump:simple}
    \sp{Suppose 
    $\cR_x^2\triangleq {\sup_{\bx\in\dom f}\bD_{\cX}(\bx,\bx^0)}<\infty$ and $\cR_y^2\triangleq {\sup_{\by\in\dom h}\bD_{\cY}(\by,\by^0)}<\infty$.} Let $\nsa{L_{x_{i} x_{i'}}}=L_{xx}$, $L_{x_iy_j}=L_{xy}$, $\nsa{L_{y_jx_i}}=L_{yx}$, and $\nsa{L_{y_{j} y_{j'}}}=L_{yy}$, for all \nsa{$i,i'\in\cM$ and $j,j'\in\cN$}, and we define $L_x\triangleq L_{xx}+L_{yx}+
\tfrac{N-1}{N}
\nsa{L_{xy}}$, and $L_y\triangleq L_{yy}+L_{yx}+\tfrac{N-1}{N}
\nsa{L_{xy}}$. 
\sp{Moreover, $\grad\varphi_\cX$, $\grad\varphi_\cY$ and $\cL$ 
are Lipschitz\footnote{\sp{Since $Z=\dom f\times \dom h$ is compact and $\grad \Phi$ is Lipschitz on $Z$, 
if $f$ and $h$ are Lipschitz, 
then $\cL$ is Lipschitz on $Z$ as well.}} on $\dom f$, $\dom h$, and $\dom f\times \dom h$, respectively.} 
\end{assumption}
\begin{corollary}\label{cor:stoch-complexity}
\sp{Suppose Assumption~\ref{assump:simple} holds for some $\cR_x,\cR_y\gg 1$. 
} Under the premise of Theorem \ref{thm:single-sample}, choosing $\gamma_1=1/L_{xx}$, $\lambda_1=1/L_{yy}$, $\gamma_2=1/L_{xy}$, $\lambda_2=1/L_{yx}$, \sp{$\alpha^0=M\max\{NL_{xx},(N-1)L_{xy}\}$ and $\tilde\alpha^0=MN\max\{L_{yy},L_{yx}\}$}
implies that $\bar\bz^K=(\bar\bx^K,\bar\by^K)$ satisfies 
{\small
\begin{equation}\label{eq:bound-O-1-s}
    \mE\left[G(\bar\bz^K)\right]\leq  \widetilde{\cO}\left(\frac{MN(L_x\cR_x^2+L_y\cR_y^2)+\max\{\frac{M}{N}\frac{1}{L_{xx}}\frac{\delta_x^2}{v},~\frac{N}{M}\frac{1}{L_{yy}}\frac{\delta^2_y}{u}\}}{\sqrt{K}+\max\{M,N\}-1}\right),\quad \forall~K\geq 1,
\end{equation}}%
where $G(\cdot)$ denotes the standard gap function.
\end{corollary}
\begin{proof} 
\sp{Let $Z=\dom f\times \dom h$} and 
\sp{note that $C_{x_i}=L_{x x}$, $L_{\by x_i}=L_{y x}$ for $i\in\cM$, and $C_{y_j}=L_{y y}$, $L_{\bx y_j}=L_{x y}$ for $j\in\cN$. Therefore, the step-sizes in \eqref{eq:step-k} satisfy $\tau_i^k=\tau^k$ for $i\in\cM$ and $\sigma_j^k=\sigma^k$ for $j\in\cN$ and $k\geq 0$ such that} $\frac{1}{\tau^k}=\frac{1}{t^k}\cO(L_{xx}MN+L_{xy}M(N-1)+L_{yx}MN)$ and $\frac{1}{\sigma^k}=\frac{1}{t^k}\cO(L_{yy}MN+L_{yx}MN+L_{xy}M(N-1))$. Therefore, 
$\frac{1}{\tau^k}=\sp{\frac{1}{t^k}\cO(MN L_x)}$  
and $\frac{1}{\sigma^k}=\sp{\frac{1}{t^k}\cO(MN L_y)}$. 
Note that $\nsa{\sup_{\bz\in Z}\Delta(\bz)=}\cO((\frac{1}{\tau^0}+(N-1)ML_{xx}+ N L_x)\cR_x^2+(\frac{1}{\sigma^0}+NML_{yy}+M L_y)\cR^2_y)$. Thus, $\sup_{\bz\in Z}\Delta(\bz)=\cO(MN(L_x\cR_x^2+L_y\cR_y^2))$. \nsa{\sp{On the other hand, since $\cL(\cdot)$ is Lipschitz on $Z$, say with constant $L_{\cL}>0$,} it implies that 
$\sup_{\bz\in Z} C_0(\bz)=\cO((M-1)\cR_x \sp{L_\cL}+(N-1)\cR_y\sp{L_\cL})$. Thus, \sp{for $\cR_x,\cR_y\gg 1$ sufficiently large}, we have $\bar{\Delta}=\cO(\sup_{\bz\in Z}\Delta(\bz))=\sp{\cO\Big(MN(L_x\cR_x^2+L_y\cR_y^2)\Big)}$.} Furthermore, since $\cR_x^2\gg 1$ and $\cR_y^2\gg 1$, according to \eqref{eq:C1}, we have \sp{$C_1(\bz^*)= \max\{\tfrac{M-1}{M} L^2_{\varphi_\cX},\tfrac{N-1}{N}L^2_{\varphi_\cY}\} \cO(\sup_{\bz\in Z}C_0(\bz))$. Finally, for our choice of $\alpha^0$ and $\tilde\alpha^0$, we get $B_\delta=\cO(\max\{\frac{M}{N L_{xx}}\frac{\delta_x^2}{v},~\frac{N}{M L_{yy}}\frac{\delta_y^2}{u}\})$}
and
we have the desired result in \eqref{eq:bound-O-1-s} since 
\sp{the term $\cH$ in Remark~\ref{rem:gap} is $\cO(\cR_x L_\cL+\cR_y L_\cL)$}. 
\end{proof}
\begin{remark}\label{rem:parameter-s}
{Corollary~\ref{cor:stoch-complexity} provides a guideline for selecting the design parameters under Assumption~\ref{assump:simple} to obtain a simplified convergence bound. For the more general setting where the coordinate-wise Lipschitz constants are not identical, one can select 
$\gamma_1=M/\sum_{i\in\cM}C_{x_i}$,
$\lambda_1=N/\sum_{j\in\cN}C_{y_j}$,
$\gamma_2=
N/\sum_{j\in\cN} L_{\bx y_j}$,
$\lambda_2=
M/\sum_{i\in\cM} L_{\by x_i}$,
and {$\alpha^0=\cO({M\max\{{N\gamma_1^{-1},(N-1)\gamma_2^{-1}}
\}})$} and $\tilde\alpha^0=\cO({MN\max\{{\lambda_1^{-1},\lambda_2^{-1}}
\}})$.} 
\end{remark}
\subsubsection{Deterministic Setting} 
\sp{
For RB-PDA, stated in Algorithm \ref{alg:RB-PDA}, suppose the mini-batch sizes are set to $|\cB^k|=|\cS^k|=1$, and $\bar\delta=\max\{\delta_x,\delta_y\}=0$; hence, $\grad_\bx\Phi_{\cB^k}(\cdot)=\grad_\bx\Phi(\cdot)$ and $\grad_\by\Phi_{\cS^k}(\cdot)=\grad_\by\Phi(\cdot)$.}

\begin{theorem}\label{thm:variable-sample}
Let Assumptions~\ref{assum:bound}-\sa{\ref{assum-lip}} hold. \sp{Suppose $\grad\varphi_\cX$ and $\grad\varphi_\cY$ are Lipschitz on $\dom f$ and $\dom h$ with constants $L_{\varphi_\cX}$ and $L_{\varphi_\cY}$, respectively.}  Initialized from arbitrary vectors $\bx^0\in\cX$ and $\by^0\in\cY$, let $\{\bx^k,\by^k\}_{k\geq 0}$ be the sequence generated by RB-PDA, stated in Algorithm \ref{alg:RB-PDA} for parameters chosen as follows: 
\sa{$t^k=\theta^k=1$,  
\sp{$\bT^k\triangleq \diag([1/\tau_i]_{i\in \cM})$ and $\bS^k\triangleq \diag(([1/\sigma_j]_{j\in\cN})$} for $k\geq 0$ such that}
{\small
\begin{subequations}
\label{eq:step-size}
\begin{align}
&\tau_i\sa{=}\tfrac{1}{\sp{(1+\varsigma)}M}\Big( L_{x_ix_i}+(N-1){(\gamma_1^{-1}+\gamma_1C_{x_i}^2+\tfrac{M+1}{M}\gamma_2^{-1})} +N {\lambda_2} L_{\by x_i}^2 \Big)^{-1},\label{eq:step-size-tau}\\
&\sigma_j\sa{=}\tfrac{1}{\sp{(1+\varsigma)}N}\Big(L_{y_jy_j}+M({\lambda_1^{-1}+\lambda_2^{-1})+M\sa{\lambda_1}C_{y_j}^2}+\tfrac{N-1}{N}(M+1){\gamma_2} L^2_{\bx y_j}\Big)^{-1}, \label{eq:step-size-s}
\end{align}
\end{subequations}}%
where \sp{$\gamma_1,\gamma_2,\lambda_1,\lambda_2,\rev{\varsigma}>0$ are arbitrary.}
\rev{If $\dom f\times \dom h$ is bounded, then let $Z=\dom f\times \dom h$; otherwise, for any $r>0$, $p\in (0,1)$ and $\bz^*\in Z^*$, let $Z\subset \cZ$ be the bounded open set given in \eqref{eq:Z-set} with $\Delta_p\triangleq \frac{2}{\varsigma p}B(\bz^*)$ such that $\alpha^0=\tilde\alpha^0=0$ chosen within the definition of $B(\bz^*)$ in Theorem~\ref{thm:bounded-stochastic}. Then, for the case $\dom f\times \dom h$ is unbounded, $\{\bz^k\}\subset Z$ with probability at least $1-p$. Moreover, the following results also hold for both bounded and unbounded problem domain scenarios:}

{\bf (I)} 
For any 
$K\geq 1$, the \sa{averaged iterates} $\bar\bz^K\triangleq (\bar\bx^K,\bar\by^K)$ defined as $\bar{\bx}^K\triangleq\frac{1}{K+M-1}(M\bx^K+\sum_{k=1}^{K-1}\bx^{k})$ and $\bar{\by}^K\triangleq\frac{1}{K+N-1}(N\by^K+\sum_{k=1}^{K-1}\by^{k})$ satisfy 
$\erfan{\mE\left[\sp{\cG_Z}(\bar \bz^K,K)\right]}\leq (\bar\Delta\sp{+C_1(\bz^*)})/K,\quad \sp{\forall \bz^*\in Z^*}$,
where $\bar \Delta\triangleq \sup_{\bz\in Z}\{\Delta(\bz)+{C_0(\bz)}\}<+\infty$ and ${C}_1(\bz^*)$ are as in Theorem~\ref{thm:single-sample} with $\cA^0=\mathbf{0}$ and $\tilde\cA^0=\mathbf{0}$.

{\bf (II)} 
The number of primal and dual gradient calls 
for achieving 
$\sp{\mE}[\sp{G_Z}(\bar\bz^K)]\leq \varepsilon$ is 
$\cO(\bar D/\varepsilon)$, where $\bar D\triangleq \bar \Delta+C_1(\bz^*)+c\cH$ and 
$0<c,\cH<+\infty$ are defined in Remark \ref{rem:gap}.

{\bf (III)} 
$\{\bx^k,\by^k\}_{k\geq 0}$ converges to a saddle point \erfan{almost surely}.
\end{theorem}
\begin{proof}
See Section \ref{sec:deterministic-proof} in the Appendix.
\end{proof}
Next, we simplify the rate result 
by considering the case where the coordinate-wise Lipschitz constants \nsa{are the same}, similar to Corollary \eqref{cor:stoch-complexity} 
to illustrate the effect of the number of blocks on the convergence rate.

\begin{corollary}\label{col:rate}
Suppose Assumption~\ref{assump:simple} holds for some $\cR_x,\cR_y\gg 1$. Under the premise of Theorem \ref{thm:variable-sample}, 
\nsa{for all $K\geq 1$}, $\bar\bz^K=(\bar\bx^K,\bar\by^K)$ satisfies 
{\small
\begin{equation}\label{eq:bound-O-1}
    \mE\left[G(\bar\bz^K)\right]\leq 
\cO\left(\frac{MN\left(L_x\cR_x^2+L_y\cR_y^2\right)
    }{K+\max\{M,N\}-1}\right).
\end{equation}}%
\end{corollary}
\begin{proof}
\sp{It follows from the same arguments in the proof of Corollary~\ref{cor:stoch-complexity} by setting $\delta_x=\delta_y=0$.}
\begin{remark}\label{rem:parameter-d}
\rev{When the coordinate-wise Lipschitz constants are not identical, in the deterministic setting, one can select 
$\gamma_1=\frac{M}{\sum_{i\in\cM}C_{x_i}}$,
$\lambda_1=\frac{N}{\sum_{j\in\cN}C_{y_j}}$,
$\gamma_2={\sqrt{\frac{N}{M}}}
\frac{N}{\sum_{j\in\cN} L_{\bx y_j}}$,
and 
$\lambda_2={\sqrt{\frac{M}{N}}}
\frac{M}{\sum_{i\in\cM} L_{\by x_i}}$.}
\end{remark}
\begin{remark}
	\label{rem:long-steps}
	  \rev{Randomized block-coordinate methods employ possibly larger steps in expectation due to their step-sizes being dependent on the \textit{average} of coordinate-wise Lipschitz constants (because of uniform randomization of the coordinates) when compared to the step sizes for full gradient methods, which are determined by the \textit{largest} coordinate-wise Lipschitz constant, e.g., see \cite{nesterov2012efficiency}, where $\min_{\bx} g(\bx)$ with $\bx=[x_i]_{i\in\cM}$ is considered and the 
	  block coordinate method therein uses a constant stepsize of $1/L_i$ {when 
      $i\in\cM$ is sampled with some predetermined probability} while the full gradient method uses a step-size of $1/L$ with $L=\max_{i\in\cM}L_i$, where $L_i$ denotes the block Lipschitz constant of $\grad g$ with respect to $x_i$.}
	
	\rev{As an example for our setting, we will compare the effective primal step-sizes of RB-PDA algorithm with the full gradient variant for the deterministic setting. For instance, similar to the quadratic example in Remark~\ref{rem:L-constants}, consider $\Phi(\bx,\by)=\frac{1}{2}\bx^\top A\bx+\by^\top Q\bx-\frac{1}{2}\by^\top B\by$ where $A\in\reals^{m\times m}$ and $B\in\reals^{n\times n}$ are diagonal matrices such that $A_{ii}=1$ for $i=1,\ldots,m-1$ and $A_{mm}=m+1$, and $B_{jj}=1$ for $j=1,\ldots,n-1$ and $B_{nn}=n+1$; moreover, $Q=\mathbf{1}_{n\times m}$ is a matrix with all entries being $1$.} 
	
	\rev{First, we consider the setting with $M=m$, i.e., $m_i=1$ for $i\in\cM$, and $N=n$, i.e., $n_j=1$ for $j\in\cN$. 
    The coordinate-wise Lipschitz constants satisfy $L_{x_i x_i}=1$ for $i\in\cM\setminus\{M\}$ and $L_{x_M x_M}=M+1$ while $L_{x_i x_\ell}= 0$ for $\ell\in\cM\setminus\{i\}$ and $i\in\cM$. Similarly, $L_{y_j y_j}=1$ for $j\in\cN\setminus\{N\}$ and $L_{y_N y_N}=N+1$ while $L_{y_j y_\ell}= 0$ for $\ell\in\cN\setminus\{j\}$ and $j\in\cN$. Furthermore, $L_{x_i y_j}=L_{y_j x_i}=1$ for all $i\in\cM$ and $j\in\cN$. 
    One can also verify that $L_{\by x_i}=1$ for $i\in\cM$ and $L_{\bx y_j}=1$ for $j\in\cN$; moreover, $C_{x_i}=\frac{1}{\sqrt{M}}$ for $i=1,\ldots,M-1$ and $C_{x_M}=\frac{M+1}{\sqrt{M}}$, and similarly $C_{y_j}=\frac{1}{\sqrt{N}}$ for $j=1,\ldots,N-1$ and $C_{y_N}=\frac{N+1}{\sqrt{N}}$.}
	
	\rev{Note that the primal and dual steps of RB-PDA can be rewritten as follows after a proper scaling (scaled direction vectors $\hat s_{j_k}^k$, $\hat r_{i_k}^k$ and scaled primal-dual step sizes $\tilde\sigma_{j_k}^k$, $\tilde\tau_{i_k}^k$ are introduced so that we can compare multi-block and full gradient settings in a unified manner):
    {\small
		\begin{align*}
			&\hat s_{j_k}^k\gets \grad_{y_{j_k}}\Phi (\bx^k,\by^k)+M\theta^k\Big(\grad_{y_{j_k}}\Phi(\bx^k,\by^k)-\grad_{y_{j_k}}\Phi(\bx^{k-1},\by^{k-1})\Big)\\
			&\by^{k+1}\gets\by^k, \quad y_{j_k}^{k+1}\gets\argmin_{y\in\cY_{j_k}}\left[\frac{1}{N}h_{j_k}(y)-\langle{\hat s_{j_k}^k,y}\rangle +\frac{1}{\tilde\sigma_{j_k}}\bD_{\cY_{j_k}}(y,y_{i_k}^k)\right],\\
			&\hat r_{i_k}^k\gets \grad_{x_{i_k}}\Phi (\bx^k,\by^{k+1})+(N-1)\theta^k\Big(\grad_{x_{i_k}}\Phi(\bx^k,\by^k)-\grad_{x_{i_k}}\Phi(\bx^{k-1},\by^{k-1})\Big)\\
			&\bx^{k+1}\gets\bx^k, \quad x_{i_k}^{k+1}\gets\argmin_{x\in\cX_{i_k}}\left[\frac{1}{M}f_{i_k}(x)+\langle{\hat r_{i_k}^k,x}\rangle +\frac{1}{\tilde\tau_{i_k}}\bD_{\cX_{i_k}}(x,x_{i_k}^k)\right],
		\end{align*}}%
		where $\tilde\tau_i\triangleq M\tau_i$ for some $\tau_i$ satisfying \eqref{eq:step-size-tau} for all $i\in\cM$and $\tilde\sigma_j\triangleq N\sigma_j$  for some $\sigma_j$ satisfying \eqref{eq:step-size-s} for all $j\in\cN$, 
		 and the design parameters $\gamma_1,\lambda_1,\gamma_2,\lambda_2$ are chosen based on Remark~\ref{rem:parameter-d}.
Thus, plugging in the constants for this problem, we get $\gamma_1=\sqrt{M}/2$, $\lambda_1=\sqrt{N}/2$, $\gamma_2=\sqrt{\frac{N}{M}}$, and $\lambda_2=\sqrt{\frac{M}{N}}$, which implies that
{\small
		\begin{align*}
			\tilde\tau_i=\Omega\Big(\frac{1}{\sqrt{NM}}\Big),\ i=1,\ldots,M-1,\quad\tilde\tau_M=\Omega\Big(\frac{1}{NM^{3/2}}\Big),\\
			\tilde\sigma_j=\Omega\Big(\frac{1}{\sqrt{NM}}\Big),\ j=1,\ldots,N-1,\quad\tilde\sigma_N=\Omega\Big(\frac{1}{MN^{3/2}}\Big).
		\end{align*}}%
		Since at each iteration we sample primal and dual coordinates uniformly at random, the expected primal and dual step sizes, $\tilde\tau$ and $\tilde\sigma$, are
        $\tilde\tau:=\sum_{i\in\cM}\tilde\tau_i\frac{1}{M}=\Omega\Big(\frac{1}{\sqrt{NM}}\Big)$, and $\tilde\sigma:=\sum_{j\in\cN}\tilde\sigma_j\frac{1}{N}=\Omega\Big(\frac{1}{\sqrt{NM}}\Big)$.
		As $M=m$, i.e., $m_i=1$ for $i\in\cM$, and $N=n$, i.e., $n_j=1$ for $j\in\cN$, our choices of parameter imply that the expected primal and dual step sizes for RB-APD are 
			$\tilde\tau=\Omega\Big(\frac{1}{\sqrt{nm}}\Big)$
			and $\tilde\sigma=\Omega\Big(\frac{1}{\sqrt{nm}}\Big)$,
		respectively.}
	
	\rev{Next, we discuss the full gradient computation scenario, i.e., $M=N=1$; thus, $\cM=\{1\}$ and $\cN=\{1\}$ implying $m_1=m$ and $n_1=n$. For this case, we have that $L_{\bx\bx}=m$, $L_{\by\by}=n$, and $L_{\by\bx}=L_{\bx\by}=\sqrt{mn}$, where $L_{\bx\bx}$, $L_{\by\by}$, $L_{\by\bx}$ and $L_{\bx\by}$ are the Lipschitz constants corresponding to the full gradients. According to \cite[Remark 2.3]{hamedani2021primal}, the primal and dual step sizes have the form $\tau=(L_{\bx\bx}+L^2_{\by\bx}/\alpha)^{-1}$ and $\sigma=(\alpha+2L_{\by\by})^{-1}$ for any $\alpha>0$, respectively. Therefore, choosing $\alpha=L_{\by\by}$ implies that
    {\small
		\begin{align*}
			\tau=(L_{\bx\bx}+L_{\by\bx})^{-1}=\frac{1}{m+\sqrt{mn}},\qquad  \sigma=(2L_{\by\by}+L_{\by\bx})^{-1}=\frac{1}{2n+\sqrt{mn}}.
		\end{align*}}%
		Indeed, this choice of step sizes can be recovered as a special case by setting $N=M=1$ in \eqref{eq:step-size} and choosing $\lambda_1=\frac{1}{L_{\by\by}}$ and $\lambda_2=1/L_{\by\bx}$ following Remark~\ref{rem:parameter-d}.}
	
	\rev{When either $n\gg m$ or $m\gg n$, the expected step sizes employed by RB-PDA for $M=m$ and $N=n$, i.e., using primal and dual block coordinates, are better than the deterministic step sizes for APD method in \cite{hamedani2021primal} using full gradient computations, i.e., corresponding to RB-APD with $M=1$ and $N=1$. That is to say $\tilde\tau=\Omega(1/\sqrt{mn})>1/(m+\sqrt{mn})=\tau$ and $\tilde\sigma=\Omega(1/\sqrt{mn})>1/(2n+\sqrt{mn})=\sigma$. In conclusion, this example shows that RB-PDA takes larger step sizes on average compared to primal-dual methods employing full gradients when there is a large deviation among different block Lipschitz constants, e.g., when average and max are significantly different as in the example above where the average is $\sum_{i=1}^m L_{x_ix_i}/m=\sum_{j=1}^n L_{y_jy_j}/n=2$ while $\max_{i=1,\ldots,m} L_{x_ix_i}=m+1$ and $\max_{j=1,\ldots,n} L_{y_jy_j}=n+1$.}
\end{remark}
\section{Convergence Analysis}
\label{sec:convergence-analysis}
\rev{In this section, we present a unified convergence analysis of the proposed algorithm, covering both stochastic and deterministic settings. After introducing the necessary definitions and preliminaries, we first establish a one-step bound for the proposed method in Lemma \ref{lem:one-step}. We then consider a class of generic step-size rules in Assumption \ref{assum:step} to derive our main result, a convergence bound on the expected weighted gap function, from which the proofs of Theorem~\ref{thm:single-sample} and Theorem~\ref{thm:variable-sample} in the previous section follow directly.}
\subsection{Supporting Lemmas}\label{sec:support}
 \begin{defn}\label{def:M-matrix}
 {For $k\geq 0$, \nsa{given some positive numbers} 
 $\{\tau_i^k\}_{i\in\cM}$ and $\{\sigma_j^k\}_{j\in\cN}$, let 
 $\bT^k\triangleq \diag([1/\tau_i^k ~\id_{m_i}]_{i\in\cM})\in\mathbb{S}^n_{++}$, and $\bS^k\triangleq \diag([1/\sigma_j^k~\id_{n_j}]_{j\in\cN})\in\mathbb{S}^m_{++}$.}
 Moreover, for any $k\geq 0$, 
 let \sa{$\{\bM_i^k\}_{i=1}^4\subset\mathbb{S}^m$ be block diagonal matrices} defined as follows:
 {\small
 \begin{align*}
 &\bM_1^k\triangleq M\theta^k({\lambda_2}N\bL_{\by\bx}^2+(N-1)\gamma_1\bC_{\bx}^2),\quad \bM_3^k\triangleq \bT^k-M\bL_{\bx\bx}-M(N-1)\theta^k(\gamma_1^{-1}+\gamma_2^{-1})\id_m-{\gamma_2^{-1}}(N-1)\id_m,\\ 
   &\bM_2^k\triangleq M\theta^k({\gamma_2}(N-1)\bL^2_{\bx\by}+\lambda_1 N \bC_{\by}^2),\quad \bM_4^k\triangleq \bS^k-N\bL_{\by\by}- M N\theta^k({\lambda_1^{-1}+\lambda_2^{-1}})\id_n -{\gamma_2}(N-1)\bL_{\bx\by}^2,
 \end{align*}}%
 where {$\gamma_1,\gamma_2,\lambda_1,\lambda_2>0$ are free \nsa{algorithm design} parameters}.
 \end{defn}

 \begin{defn}\label{def:pos-set}
 Let $\cZ\triangleq \cX\times\cY$ and $\bz\triangleq (\bx,\by)$. Define $G_f(\bz,\bar{\bz})\triangleq f(\bar{\bx})-f(\bx)+\Phi(\bar{\bx},\by)-\Phi(\bx,\by)$, $G_h(\bz,\bar{\bz})\triangleq h(\bar{\by})-h(\by)-\Phi(\bar{\bx},\bar{\by})+\Phi(\bx,\by)$, and $G_\Phi(\bz,\bar{\bz})\triangleq \Phi(\bar{\bx},\by)-\Phi(\bar{\bx},\bar{\by})$ for $\bz,\bar{\bz}\in\cZ$.
 For 
 simplicity, let $G_f^k(\bz)\triangleq G_f(\bz,\bz^k)$, $G_h^k(\bz)\triangleq G_h(\bz,\bz^k)$, and $G_\Phi^k(\bz)\triangleq G_\Phi(\bz,\bz^k)$ for any $k\geq 0$.
 \end{defn}

 \begin{defn}
 \label{def:sigma_algebra}
 $\mE^{k}[\cdot]\triangleq\mE\left[\cdot\mid \mathcal F_{k-1}\right]$ denotes the conditional expectation 
 and $\mathcal F_{k-1}\triangleq\sigma(\Psi_{k-1})$ is the $\sigma$-algebra generated by $\Psi_{k-1}\triangleq\{\bx^0,i_0,\hdots,i_{k-1}\}\bigcup\{\by^0,j_0,\hdots,j_{k-1}\}\bigcup\sp{\{\cB^0,\cS^0,\hdots,\cB^{k-1},\cS^{k-1}\}}$ for $k\geq 1$.
 \end{defn}

 \begin{lemma}\label{lem:pos-theta}
 \sa{Given 
 $\bz=(\bx,\by)\in{\cZ}$, 
 for all $k\geq 0$, we have \sp{$G_f^k(\bz)=\cL(\bx^k,\by)-\cL(\bx,\by)$}, and} 
 {\small
 \begin{align}
 \label{eq:G-connection}
 G_h^k(\bz)  \geq \cL(\bx,\by)-\cL(\bx,\by^k)+\fprod{\grad_\bx\Phi(\bx^k,\by^k),\bx-\bx^k}, \quad
 \sp{G_\Phi^k(\bz)  
  \leq \fprod{\grad_\by\Phi({\bx}^k,{\by}^k),\by-{\by}^k}.} 
 \end{align}}%
 \end{lemma}
 \begin{proof}
 \sa{Since $\by^k\in\dom h$, using convexity of $\Phi(\cdot,\by^k)$ gives} the first inequality.
 Moreover, \sa{since $\bx^k\in\dom f$,} from the concavity of $\Phi(\bx^k,\cdot)$ we have 
 $G_\Phi^k(\bz)\leq \fprod{\grad_\by\Phi({\bx}^k,{\by}^k),\by-{\by}^k}$.
 \end{proof}

 \begin{defn}\label{def:ek}
 For $k\geq 0$,  $\mathbf{e}_\bx^k(\bx,\by)\triangleq[{e}_{x_i}^k(\bx,\by)]_{i\in\cM}$ and $\mathbf{e}_\by^k(\bx,\by)\triangleq[{e}_{y_j}^k(\bx,\by)]_{j\in\cN}$ denote the error of sample gradients with respect to $\bx$ and $\by$ respectively, where $e_{x_i}^k$ and $e_{y_j}^k$ are 
  defined as
  {\small
 \begin{align}
 \label{eq:error_seq} 
 &e_{x_i}^k(\bx,\by)\triangleq \sp{\grad_{x_i}\Phi_{\cS^k}(\bx,\by)}-\grad_{x_i}\Phi(\bx,\by),\quad e_{y_j}^k(\bx,\by)\triangleq \sp{\grad_{y_j}\Phi_{\cB^k}(\bx,\by)}-\grad_{y_j}\Phi(\bx,\by).
 \end{align}}%
 Moreover, let $\bw_\by^k=[w_{y_j}^k]_{j\in\cN}$ and $\bw_\bx^k=[w_{x_i}^k]_{i\in\cN}$ denote the error sequence due to using sampled gradients in updating $\by^{k+1}$ and $\bx^{k+1}$ respectively, i.e.,
 {\small
 \begin{align}\label{eq:wk-y}
 &\bw^k_\by\triangleq \mathbf{e}^k_\by(\bx^k,\by^{k})+M\theta^k\left({\mathbf{e}}^k_\by(\bx^k,\by^{k})-{\mathbf{e}}^{k}_\by(\bx^{k-1},\by^{k-1})\right),  \\
 \label{eq:wk}
 &\bw^k_\bx\triangleq \mathbf{e}^k_\bx(\bx^k,\by^{k+1})+(N-1)\theta^k\left({\mathbf{e}}^k_\bx(\bx^k,\by^{k})-{\mathbf{e}}^{k}_\bx(\bx^{k-1},\by^{k-1})\right),   
 \end{align}}%
 for $k\geq 0$, such that $\mathbf{e}^{0}_\bx(\bx^{-1},\by^{-1})=0$ and $\mathbf{e}^{0}_\by(\bx^{-1},\by^{-1})=0$.
 \end{defn}
Next, we 
\sp{consider the progress in one iteration of RB-PDA,} 
which is the main building \us{block for deriving} rate results 
in Theorems \ref{thm:single-sample} and \ref{thm:variable-sample}. 
First, we need to give some definitions.
\sp{Let $\cA^k,\tilde\cA^k\in\mathbb{S}^m_{++}$ be diagonal matrices and we define two auxiliary sequences $\{\bu_\bx^k\}_{k\geq 0}$ and $\{\bu_\by^k\}_{k\geq 0}$ such that $\bu_\bx^0=\bx^0$ and $\bu_\by^0=\by^0$, and
{\small
 \begin{align}
 \label{eq:uk}
 \bu_\bx^{k+1}\triangleq \argmin_{\bx\in\cX}\Big\{-\fprod{\sp{\bw_\bx^k},\bx}+\bD_\cX^{\frac{1}{M}\cA^k}(\bx,\bu_\bx^k)\Big\},\quad
 \bu_\by^{k+1}\triangleq \argmin_{\by\in\cY}\Big\{\fprod{\bw_\by^k,\by}+\bD_\cY^{\frac{1}{N}{\tilde\cA}^k}(\by,\bu_\by^k)\Big\}, 
 \end{align}
 }}%
 for all $k\geq 0$, which is helpful to bound the stochastic error of sampled gradient using Lemma \ref{lem:inner-w}.
 \sp{For $k\geq 0$, let $\btr^k\triangleq M\grad_\bx\Phi(\bx^k,\by^{k+1})+\tilde\bq_\bx^k$ where $\tilde\bq_\bx^k\triangleq M(N-1)\theta^k(\grad_\bx\Phi(\bx^k,\by^{k})-\grad_\bx\Phi(\bx^{k-1},\by^{k-1}))$, and similarly, we define $\tilde\bs^k\triangleq M\grad_\by\Phi(\bx^k,\by^k)+\tilde\bq_\by^k$ where $\tilde\bq_\by^k\triangleq MN\theta^k(\grad_\by\Phi(\bx^k,\by^k)-\grad_\by\Phi(\bx^{k-1},\by^{k-1}))$. Moreover, we define}\vspace*{0mm}
 \sp{
 \begin{subequations}
 {\small
    \begin{align}
    \bs^k\triangleq [s_j^k]_{j\in\cN}\quad\mbox{s.t.}\quad s_j^k=N\grad_{y_j}\Phi_{\cB^k}(\bx^k,\by^k)+q_{y_{j}}^k,\quad\forall~j\in\cN,  
    \end{align}}%
    where $q_{y_{j}}^k\triangleq N M \theta^k(\grad_{y_{j}}\Phi_{\cB^k}(\bx^k,\by^k)-\grad_{y_{j}}\Phi_{\cB^k}(\bx^{k-1},\by^{k-1}))$ and $\bq_\by^k\triangleq [q_{y_j}^k]_{j\in\cN}$. Similarly we define
    {\small
    \begin{align}
    \br^k\triangleq [r_i^k]_{j\in\cM}\quad\mbox{s.t.}\quad r_i^k=M\grad_{x_{i}}\Phi_{S^k}(\bx^k,\by^{k+1})+q_{x_{i}}^k,\quad\forall~i\in\cM,  
    \end{align}}%
    where $q_{x_{i}}^k\triangleq M(N-1)\theta^k(\grad_{x_{i}}\Phi_{S^k}(\bx^k,\by^{k})-\grad_{x_{i}}\Phi_{S^k}(\bx^{k-1},\by^{k-1}))$ and $\bq_\bx^k\triangleq [q_{x_i}^k]_{i\in\cM}$.
\end{subequations}}%
Finally, let $\{\btx^k\}_{k\geq 1}\subseteq\cX$ and $\{\bty^k\}_{k\geq 1}\subseteq\cY$ be some auxiliary sequences such that for all $k\geq 0$,
\begin{subequations}
{\small
 \begin{align}
 \tilde{y}_j^{k+1} \triangleq \argmin_{y_j\in\cY_j} &~ \us{\left[h_j(y_j)-\fprod{s_j^k, ~y_j}+{\frac{1}{\sigma_j^k}}\bD_{\cY_j}(y_j,y_j^k)\right]},\quad j\in\cN,\label{eq:tilde-problem-y} \\
 \tilde{x}_i^{k+1} \triangleq \argmin_{x_i\in\cX_i} &~\us{\left[f_i(x_i)+\fprod{r_i^k, ~x_i}+{\frac{1}{\tau_i^k}}\bD_{\cX_i}(x_i,x_i^k)\right]},\quad i\in\cM,  \label{eq:tilde-problem-x}
 \end{align}}%
 \end{subequations}
\us{Note that} the auxiliary sequence {$\{\btx^k,\bty^k\}_{k\geq 1}$} is never actually computed in the implementation of the algorithm RB-PDA and it is defined \sa{to} 
 facilitate the analysis of convergence of $\{\bx^k\}_{k\geq 1}\subseteq\cX$ and $\{\by^k\}_{k\geq 1}\subseteq\cY$. \sp{Moreover, note that according to \eqref{eq:wk}, $\bw^k_\bx= \frac{1}{M}(\br^k-\btr^k)$ and $\bw_\by^k= \frac{1}{N}(\bs^k-\tilde\bs^k)$.}
 \begin{lemma}\label{lem:one-step}
 (One-step analysis) Let $\{\bx^k,\by^k\}_{k\geq 0}$ be the sequence generated by RB-PDA, stated in Algorithm \ref{alg:RB-PDA} {initialized} from arbitrary $\bx^0\in\cX$ and $\by^0\in\cY$ \rev{such that $\theta^k\geq 1$ for all $k\geq 0$}. Suppose Assumptions \ref{assum-lip} 
 holds.
 \sa{Given $\{\cA^k\}_k$ and $\{\tilde\cA^k\}_k$,  arbitrary sequences of diagonal matrices with positive entries,} 
 \begin{subequations}\label{eq:one-step}
 {\small
 \begin{align}\label{eq:one-step-main}
 &\cL(\bx^{k+1},\by)-\cL(\bx,\by^{k+1})+(\theta^k-1)\big[(M-1)(\cL(\bx^k,\by)-\cL(\bx,\by))+(N-1)(\cL(\bx,\by)-\cL(\bx,\by^k))\big]\\
 &\leq (M-1)(\theta^kG_f^k(\bz) -G_f^{k+1}(\bz))+(N-1)(\theta^kG_h^k(\bz)-G_h^{k+1}(\bz))+ M(G_\Phi^{k+1}(\bz)-\theta^kG_\Phi^k(\bz)),\nonumber\\
 &\quad+ Q^k(\bz)-R^{k+1}(\bz)+M\norm{\bw_\bx^k}^{\sa{2}}_{\nsa{*,(\cA^k)^{-1}}}+N\norm{\bw_\by^k}^2_{{*,(\tilde\cA^k)^{-1}}}+\fprod{\bw_\bx^k,\nsa{\bu_\bx^k-\bx^k}}+\fprod{\bw_\by^k,\by^k-\bu_\by^k}+{\Xi_x^k(\bx)}+{\Xi_y^k(\by)}, \nonumber
 \end{align}}%
 holds for all $k\geq 0$ and $\bz=(\bx,\by)\in{\cZ}$, where $G_f^k$, $G_h^k$, and $G_\Phi^k$ are given in Definition \ref{def:pos-set}, \nsa{$\bw_x^k,\bw_\by^k$ and $\bu_x^k,\bu_\by^k$ are given in~\eqref{eq:wk} and \eqref{eq:uk}, respectively;} and \sa{$Q^k(\cdot)$, $R^{k+1}(\cdot)$, ${\Xi_x^k(\bx)}$ and ${\Xi_y^k(\by)}$ are defined as follows:}
 {\small
 \begin{align}
  & Q^k(\bz)\triangleq \bD_\cX^{\bT^k}(\bx,\bx^k)+ \bD_\cY^{\bS^k}(\by,\by^k) +\tfrac{1}{M}\bD_\cX^{\cA^k}(\bx,\bu_\bx^k)+\ey{\tfrac{1}{N}\bD_\cY^{\tilde\cA^k}(\by,\bu_\by^k)}+\bD_\cX^{\bM_1^k}(\sa{\bx^k,\bx^{k-1}})  \label{eq:Qk-def}\\ 
 &\nonumber\quad+\bD_\cY^{\bM_2^k}(\sa{\by^k,\by^{k-1}})-M\theta^k\fprod{\grad_\by\Phi(\bx^{k-1},\by^{k-1}),~\by^k-\by} -(N-1)\theta^k\fprod{\grad_\bx\Phi(\bx^{k-1},\by^{k-1}),~\bx-\bx^k}, \nonumber\\ 
 \label{eq:Rk-def}
 &R^{k+1}\sa{(\bz)}\triangleq \bD_\cX^{\bT^k}(\bx,\bx^{k+1})+ \bD_\cY^{\bS^k}(\by,\by^{k+1})+\tfrac{1}{M}\bD_\cX^{\cA^k}(\bx,\bu_\bx^{k+1})\ey{+\tfrac{1}{N}\bD_\cY^{\tilde\cA^k}(\by,\bu_\by^{k+1})}+\bD_\cX^{\bM_3^k-{\cA^k}}(\sa{\bx^{k+1},\bx^{k}}) \\&\nonumber\quad+\bD_\cY^{\bM_4^k-\ey{\tilde\cA^k}}(\sa{\by^{k+1},\by^{k}}) -M\fprod{\grad_\by\Phi(\bx^k,\by^k),~\by^{\sa{k+1}}-\by}
 -(N-1)\fprod{\grad_\bx\Phi(\bx^k,\by^{k}),~\bx-\bx^{k+1}}, \\
 \label{eq:error-x}
 &{\Xi_x^k(\bx)} \triangleq Mf(\bx^{k+1})-(M-1)f(\bx^k)-f(\btx^{k+1})+M\Phi(\bx^{k+1},\by^{k+1})-\sum_{i\in\cM}\Phi(\bx^k+U_i(\tilde{x}_i^{k+1}-x_i^k),\by^{k+1})\\ \nonumber
 &  \quad  +\erfan{\tfrac{1}{2}\norm{\bx^k-\bx^{k-1}}^2_{\bU^k}}+\tfrac{1}{M}\Big(\bD_{\cX}^{\bT^k}(\bx,\bx^k) -\bD_{\cX}^{\bT^k}(\bx,\tilde{\bx}^{k+1})-\nsa{\bD_{\cX}^{\bT^k-M\bL_{\bx\bx}-\cA^k}}(\tilde{\bx}^{k+1},\bx^k)\Big)\nonumber\\ \nonumber
 &  \quad -\Big(\bD_{\cX}^{\bT^k}(\bx,\bx^k) -\bD_{\cX}^{\bT^k}(\bx,{\bx}^{k+1})-\bD_{\cX}^{\bT^k-M\bL_{\bx\bx}-\cA^k}({\bx}^{k+1},\bx^k)\Big)+\tfrac{1}{M}\fprod{\erfan{\tilde\bq}_\bx^k,M\bx^{k+1}-(M-1)\bx^k-\btx^{k+1}},\\
 \label{eq:error-y}
 &{\Xi_y^k(\by)}\triangleq  Nh(\by^{k+1})-(N-1)h(\by^k)-h(\bty^{k+1})
 +\sum_{j\in\cN}\Phi(\bx^k,\by^k+V_j(\tilde y_j^{k+1}-y_j^k))-N\Phi(\bx^k,\by^{k+1})\\\nonumber
 &  \quad +\erfan{\tfrac{1}{2}\norm{\by^{k+1}-\by^k}^2_{\bV^k_1}+\tfrac{1}{2}\norm{\by^k-\by^{k-1}}^2_{\bV^k_2}}+\tfrac{1}{N}\Big(\bD_{\cY}^{\bS^k}(\by,\by^k) -\bD_{\cY}^{\bS^k}(\by,\tilde{\by}^{k+1})-\bD_{\cY}^{\bS^k-N\bL_{\by\by}-\ey{\tilde\cA^k}}(\tilde{\by}^{k+1},\by^k)\Big)\nonumber\\ \nonumber
 &  \quad -\Big(\bD_{\cY}^{\bS^k}(\by,\by^k) -\bD_{\cY}^{\bS^k}(\by,{\by}^{k+1})-\bD_{\cY}^{\bS^k-N\bL_{\by\by}-\ey{\tilde\cA^k}}({\by}^{k+1},\by^k)\Big)+\tfrac{1}{N}\fprod{\erfan{\tilde\bq}_\by^k,\bty^{k+1}-N\by^{k+1}+(N-1)\by^k},
 \end{align}}%
 \erfan{where {\small $\bU^k\triangleq M\theta^k[ (N-1)\gamma_1(L^2_{x_{i_k}x_{i_{k-1}}}\bI_m-\bC_\bx^2)+ N\lambda_2(L^2_{y_{j_k}x_{i_{k-1}}}\bI_m-\bL^2_{\by\bx})]$}, {\small $\bV^k_1\triangleq (N-1)\gamma_2(L^2_{x_{i_k}y_{j_k}}\bI_n-\bL^2_{\bx\by})$}, {\small $\bV^k_2\triangleq M\theta^k[ N\lambda_1(L^2_{y_{j_k}y_{j_{k-1}}}\bI_n-\bC^2_{\by})+(N-1)\gamma_2(L^2_{x_{i_k}y_{j_{k-1}}}\bI_n-\bL^2_{\bx\by})]$}, and $\{M_i^k\}_{i=1}^4$ are given in Definition \ref{def:M-matrix}}.
 \end{subequations}
 \end{lemma}%
\begin{proof} See Section \ref{sec:proof-lem-one-step} in the Appendix.\end{proof}
 In Lemma \ref{lem:one-step}, we derived a one-step bound on the Lagrangian metric, 
 where the terms ${M\norm{\bw_\bx^k}^{\sa{2}}_{\nsa{*,(\cA^k)^{-1}}}}$ and $\rev{N}\norm{\bw_\by^k}^{\sa{2}}_{\nsa{*,(\tilde\cA^k)^{-1}}}$ provide us with an explicit error bound due to using sample gradients to approximate $\grad_{x_{i_k}}\Phi$ and $\grad_{y_{j_k}}\Phi$ respectively. \rev{These terms vanish in the deterministic setting since the exact gradients are available; hence, $\bw_\bx^k=\mathbf{0}$ and $\bw_\by^k=\mathbf{0}$ for $k\geq 0$.} In the next lemma, we will further bound \sa{this error \rev{for the stochastic setting}, in terms of the number of samples required,} using Assumption~\ref{assum:sample} and show that \nsa{the last four terms in \eqref{eq:one-step-main} have zero means.} 

 \begin{lemma}\label{error_sampling}
 \ey{Let $\{\bw_\bx^k,\bw_\by^k\}_k$ be the error sequence 
 as in \eqref{eq:wk}.}
 Suppose Assumptions~\ref{assum:sample} and~\ref{assum-lip} hold, and for all $k\geq 0$, \sp{$v^k=v$} and \sp{$u^k=u$}, \ey{ $\cA^k=\alpha^k\id_{m}$ and $\tilde\cA^k=\tilde\alpha^k\id_{n}$ for some 
 $\alpha^k,\tilde\alpha^k>0$}. Then, for $k\geq 0$, 
 {\small
 \begin{subequations}
 \begin{align}
 &{\mE}^k\Big[\norm{\bw_\bx^k}^2_{\nsa{*,(\cA^k)^{-1}}}\Big]\leq 
 \frac{2\sp{M}\delta_x^2}{v\alpha^k}+\frac{8\sp{M(N-1)^2(\theta^k)^2}}{\alpha^k}\left(\norm{\bx^k-\bx^{k-1}}_{\bC_\bx^2}^2+\norm{\by^k-\by^{k-1}}_{\bL_{\bx\by}^2}^2\right);\label{eq:part1}\\
 &{\mE}^k\Big[\norm{\bw_\by^k}^2_{\nsa{*,(\tilde\cA^k)^{-1}}}\Big]\leq \frac{2\sp{N}\delta_y^2}{u\tilde\alpha^k}+\frac{8\sp{NM^2(\theta^k)^2}}{\tilde\alpha^k}\left(\norm{\bx^k-\bx^{k-1}}_{\bL_{\by\bx}^2}^2+\norm{\by^k-\by^{k-1}}_{\bC_{\by}^2}^2\right);\label{eq:part2}\\
 &\mE^k[\fprod{\bw_\bx^k,\nsa{\bu_\bx^k-\bx^k}}]=\mE^k[\fprod{\bw_\by^k,\nsa{\by^k-\bu_\by^k}}]=0,\quad \mE^k[\Xi_x^k(\bx)]=\mE^k[\Xi_y^k(\by)]=0,\quad \forall(\bx,\by)\in\dom f\times \dom h,\label{eq:part3}
 \end{align}
 \end{subequations}}%
where \sp{$\bu_\bx^k$, $\bu_\by^k$,} ${\Xi_x^k(\bx)}$, and ${\Xi_y^k(\by)}$, are defined in \nsa{\eqref{eq:uk}}, \eqref{eq:error-x}, and \eqref{eq:error-y}, respectively. 
 \end{lemma}
 \begin{proof} See Section \ref{sec:proof-lem-error} in  the Appendix for the proof.\end{proof}
\subsection{Generic step-size rule and the meta-result}

In this section, we state 
\sp{a more generic step-size rule, and then we provide a convergence result for RB-PDA under this more general rule as a meta result.}
Indeed, \sp{in the proofs of our main results stated in}~Theorems~\ref{thm:single-sample} and~\ref{thm:variable-sample}, we show that the step-size choices given in these theorems satisfy Assumption~\ref{assum:step}.
 \begin{assumption}\label{assum:step}
 {\bf (Step-size Condition)} \sp{Let $\bar\delta\triangleq\max\{\delta_x,\delta_y\}$,
 $\bT^k$, $\bS^k$ and $\{\bM_i\}_{i=1}^4$ be as 
 in Definition \ref{def:M-matrix}.} {We assume that for $k\geq 0$, there exist \nsa{diagonal} 
 \sp{$\bT^k$, $\bS^k$,}
 $\theta^k\geq 1$, and \sp{$t^k\leq 1$} with $t^0=1$ satisfying}
 \sp{\small
 \begin{subequations}\label{eq:step-size-condition}
 \begin{align}
 &{t^{k}\bT^{k}= t^{k+1}\bT^{k+1},\quad t^{k}\bS^{k}= t^{k+1}\bS^{k+1}},\quad t^k=t^{k+1}{\theta^{k+1}},\quad t^{k}\cA^{k}\succeq t^{k+1}\cA^{k+1},\quad t^{k}\tilde\cA^{k}\succeq t^{k+1}\tilde\cA^{k+1}, 
 \label{eq:step-size-theta}\\
 &
 \sa{t^{k+1}}(\bM_1^{k+1}\erfan{+\cD^{k+1}})+\sp{\varsigma \bT^0}\preceq\sa{t^{k}}(\bM_3^k-\cA^k),\quad \sa{t^{k+1}}(\bM_2^{k+1}+\erfan{\tilde\cD^{k+1}})+\sp{\varsigma \bS^0}\preceq\sa{t^{k}}(\bM_4^k-\tilde\cA^k),\label{eq:step-size-lip} 
 \end{align}
 \end{subequations}}%
 \sp{for some $\varsigma\geq 0$ and $\gamma_1,\gamma_2,\lambda_1,\lambda_2>0$}, \sp{where $\cD^k=\cA^k=\mathbf{0}$, $\tilde\cD^k=\tilde\cA^k=\mathbf{0}$ when $\bar\delta=0$; otherwise, when $\bar\delta>0$, $\cA^k= \alpha^k\id_{m}$, $\tilde\cA^k=\tilde\alpha^k\id_{n}$ for some $\alpha^k,\tilde\alpha^k>0$ such that $\sum_{k=0}^\infty t^k\Big(\frac{1}{\alpha^k}+\frac{1}{\tilde\alpha^k}\Big)<+\infty$, and $\cD^k\triangleq 16(M\theta^k)^2\Big((N-1)^2\bC_\bx^2
 /\sp{\alpha^k}+N^2\bL_{\by\bx}^2
 /\sp{\tilde\alpha^k}\Big)$ and $\tilde\cD^k\triangleq  16(M\theta^k)^2\Big(N^2\bC^2_\by
 /\sp{\tilde\alpha^k}+(N-1)^2\bL_{\bx\by}^2
 /\sp{\alpha^k}\Big)$.}
 \end{assumption}

 \erfan{Before discussing the convergence rate result of the proposed algorithm, \sp{we first demonstrate 
 that the generated sequence by Algorithm \ref{alg:RB-PDA} 
 is bounded \rev{almost surely}, which is necessary for establishing the convergence result as well as for} defining a valid gap function; \rev{moreover, a uniform bound (over sample paths $\omega\in\Omega$) holds with high probability}.}
\begin{lemma}\label{lem:result-summable}
    Let $\{\bx^k,\by^k\}_{k\geq 0}$ be generated by RB-PDA, stated in Algorithm \ref{alg:RB-PDA}, {and initialized} from arbitrary $\bx^0\in\cX$ and $\by^0\in\cY$. Suppose Assumptions~\ref{assum:bound}-
    \ref{assum-lip} 
    and \ref{assum:step} hold with $\varsigma>0$. Then, it holds that
    $\sum_{k=0}^\infty \bD_\cX(\bx^{k+1},\bx^{k})+\bD_\cY(\by^{k+1},\by^{k}) <+\infty$ a.s. and  $\mE[\sum_{k=0}^\infty\bD^{\bT^0}_\cX(\bx^{k+1},\bx^{k})+\bD^{\bS^0}_\cY(\by^{k+1},\by^{k})]=\sum_{k=0}^\infty\mE[\bD^{\bT^0}_\cX(\bx^{k+1},\bx^{k})+\bD_\cY^{\bS^0}(\by^{k+1},\by^{k})]\leq (B(\bz^*)+\sum_{k=0}^\infty B^k_\delta)/\varsigma$ for all $\bz^*\in Z^*$, where $B(\bz)\triangleq \Delta(\bz)+ C_0(\bz)+\cL(\bx^0,\by)-\cL(\bx,\by^0)$, $\sp{B^k_\delta}\triangleq 0$ 
 \sp{when $\bar\delta=0$;} $\sp{B^k_\delta}\triangleq 2\erfan{t^k}\left(\frac{\erfan{M^2}\delta_x^2}{v\alpha^k}+\frac{\erfan{N^2}\delta_y^2}{u\tilde\alpha^k}\right)$ 
 \sp{when $\bar\delta>0$}, and $\Delta(\cdot)$ and $C_0(\cdot)$ are defined in~\eqref{eq:delta-function} and \eqref{eq:C0-function}, respectively. 
 
 \rev{Finally, 
 for any $p\in (0,1)$ and $\bz^*\in Z^*$, let $\Delta_p\triangleq \frac{2}{\varsigma p}\Big(B(\bz^*)+\sum_{k=0}^\infty B_\delta^k\Big)=\cO(1/p)$. Then, it holds that 
 $\mathbb P\!\left(\sup_{k\ge0}\{\| \bx^k-\bx^*\|^2_{\bT^0}+\| \by^k-\by^*\|^2_{\bS^0}\}\le \Delta_p\right)\ge 1-p$.}
\end{lemma}
\begin{proof} See Section \ref{sec:proof-lem-summable} in the Appendix for the proof.\end{proof}
\rev{Next, we provide a bound that will play a key role for establishing the convergence rate results.} 
\begin{lemma}
    \label{lem:sup-bound}
    Under the premise of 
    \sp{Lemma~\ref{lem:result-summable},} the following bound holds:
    {\small
    \begin{align}\label{eq:bound-delta}
 \MoveEqLeft\mE\Big[\sup_{\bz\in Z} \Big\{C_0(\bz)+ \bD_\cX^{\bT^0+(N-1)M\bL_{\bx\bx}+\frac{1}{M}\cA^0}(\bx,\bx^0)+\bD_\cY^{\bS^0+NM\bL_{\by\by}+\ey{\frac{1}{N}\tilde\cA^0}}(\by,\by^0)\nonumber \\
 &+\sum_{k=0}^{K-1}t^k\Xi_x^k(\bx)+\sum_{k=0}^{K-1}t^k\Xi_y^k(\by)\Big\}\Big]\leq  \sup_{\bz\in Z}\Big\{ C_0(\bz)+\Delta(\bz)\Big\}+\sp{ C_1(\bz^*)+ C_2,\quad\forall~\bz^*\in Z^*,} 
 \end{align}}%
 \sp{where $C_1(\cdot)$ is defined in~\eqref{eq:C1} and $C_2\triangleq \tfrac{1}{\varsigma}\max\{\tfrac{M-1}{M} L^2_{\varphi_\cX},\tfrac{N-1}{N}L^2_{\varphi_\cY}\} \sum_{k=0}^{\infty}B_\delta^k$.}
\end{lemma}
\begin{proof} \sp{See Section~\ref{sec:proof-sup-bound} in the Appendix for the proof.}\end{proof}

Based on the boundedness result of Lemma \ref{lem:result-summable}, we can define the restricted gap function as follows.
\begin{defn}
\label{def:restricted-domain}
    \rev{If $\dom f\times \dom h$ is bounded, then let $Z=\dom f\times \dom h$; otherwise, for arbitrary $\bz^*\in Z^*$ and $p\in (0,1)$, let 
    the restricted domain $Z\subset \cZ$ be defined as $Z\triangleq\{(\bx,\by)\in\cX\times\cY:\ \norm{\bx-\bx^*}^2_{\bT^0}+\norm{\by-\by^*}^2_{\bS^0}<(1+r)\Delta_p\}$ for any $r>0$ and $\Delta_p>0$ as in Lemma~\ref{lem:result-summable}.}
\end{defn}

\sp{Next, we state our meta-theorem that establishes the convergence guarantees of RB-PDA under the generic step-size rule stated in Assumption~\ref{assum:step}. The rate and complexity guarantees given in our main results directly follow from the meta-theorem by showing that the particular step-size rules for the stochastic and deterministic scenarios satisfy Assumption~\ref{assum:step}.}
\begin{theorem}\label{thm:meta}
\sp{Suppose Assumptions~\ref{assum:bound}-\ref{assum-lip} hold. Let $\{\bx^k,\by^k\}_{k\geq 0}$ be the 
RB-PDA sequence, {initialized} from arbitrary 
$\bx^0\in\cX$ and $\by^0\in\cY$, corresponding to some algorithm parameters satisfying Assumption~\ref{assum:step} when $|\cB^k|=u$ and $|\cS^k|=v$ for some $u,v\geq 1$ for all $k\geq 0$. \rev{Let $Z\subset\cZ$ be the bounded set given in Definition~\ref{def:restricted-domain}.
If $\dom f\times \dom g$ is unbounded, then $\{\bz^k\}\subset Z$ with probability at least $1-p$; moreover, for both bounded and unbounded problem domain scenarios,} 
$\bar\bz^K=(\bar\bx^K,\bar\by^K)$ defined as in~\eqref{eq:avg-iterate}
satisfies $\mE\left[\cG_Z(\bar \bz^K,T_K)\right] \leq\frac{\bar \Delta
+\sp{C_1(\bz^*)}}{T_K}+\frac{1}{T_K}\Big(1+\frac{1}{\varsigma}\max\{\tfrac{M-1}{M} L^2_{\varphi_\cX},\tfrac{N-1}{N}L^2_{\varphi_\cY}\}\Big)\sum_{k=0}^\infty B_\delta^k$ for all $K\geq 1$, where $T_K\triangleq\sum_{k=0}^{K-1} t^{k}$, $\bar \Delta\triangleq \sup_{\bz\in Z}\{\Delta(\bz)+{C_0(\bz)}\}$, $C_0(\bz),{C}_1(\bz)$ are given in \eqref{eq:Delta-C-def}, and $B_\delta^k$ is defined in Lemma \ref{lem:result-summable}. Finally, in the deterministic case, $\delta_x=\delta_y=0$, 
one has 
$\mE\left[\cG_Z(\bar \bz^K,T_K)\right] \leq\frac{\bar \Delta
+\sp{C_1(\bz^*)}}{T_K}$ for $K\geq 1$.}
\end{theorem}
\begin{proof} \sp{See Section~\ref{sec:meta-proof} in the Appendix for the proof.}\end{proof}
\section{Numerical Experiments}\label{numerics}
In this experiment, we aim to compare the performance of different blocking strategies for our proposed RB-PDA method against traditional primal–dual methods for solving SP problems. 
Consider a set of data points $\{a_j\}_{j\in\cN}\subset \reals^m$ with $\{b_j\}_{j\in\cN}\subset \{-1,+1\}$ labels.
Let $L:\cX\times \reals^m\times \{-1,1\}\to \reals$ \rev{be a convex loss function}, and we define $L_j(\bx)\triangleq L(\bx;a_j,b_j)$. Distributionally robust optimization studies worse case performance under uncertainty to find solutions with some specific confidence level \citep{namkoong2016stochastic,shalev2016minimizing}. This problem can be formulated as follows:
$\min_{\norm{\bx}_\infty\leq R} \max_{\by=[y_j]_{j\in\cN}\in\cU}~ \sum_{j\in\cN} y_j L_j(\bx)$, 
where $\cU$ is an uncertainty set. As suggested in \citep{namkoong2016stochastic}, 
we consider 
the uncertainty set  $\cU=\{\by\in\reals^n_+\mid \sum_{j\in\cN}y_j=1, \frac{1}{2}\norm{n\by-\ones_n}^2\leq \rho\}$, representing the family of discrete probability distributions in a neighborhood of the uniform distribution.
To have a separable structure in maximization variable 
for our method, we relax the constraint in $\by$ by introducing primal variables $w_1$ and $w_2$ for the equality and inequality constraints, respectively, to obtain an equivalent reformulation in the form of \eqref{eq:original-problem} with $f(\bw)=\ind{W}(\bw)$, $h(\by)=\ind{\reals^n_+}(\by)$ and \vspace*{0mm}
{\small
\begin{align*}
 \Phi(\bw,\by)=\sum_{j\in\cN} \by_\ell L_j(\bx)~+w_1(\ones_n^\top \by-1)-w_2(\norm{n\by-\ones_n}^2/2-\rho)/n, \vspace*{0mm}
\end{align*}}%
where $\bw=[\bx^\top, w_1, w_2]^\top$ and $W=[-R,R]^m\times \reals\times\reals_+$. 
In our experiment, we 
set 
$R=10$. We test the performance of our method under various block partitioning strategies in deterministic and stochastic settings and compare it with Stochastic Mirror Descent (SMD) \citep{juditsky2011solving} and Stochastic Mirror-prox (SMP) \citep{nemirovski2009robust} \eh{by fixing the running time of algorithms}.

\noindent\textbf{Dataset.} 
\rev{We use three real datasets: \texttt{w7a} ($24692\times 300$) from \cite{zeng2008fast}, preprocessed using LIBSVM \cite{chang2011libsvm}, \texttt{madelon} ($2000\times 500$), and \texttt{mushroom} ($8124\times 112$). To ensure a divisible dimension in the variable $\mathbf{y}$, we randomly select $95\%$ of the \texttt{w7a} data points, i.e., $A \in \mathbb{R}^{23458\times 300}$.}

\noindent\textbf{Implementation details.} All algorithms are initialized from $\bx^0={\bf 0}$ and $\by^0=\ones_n/n$. Moreover, we fixed the running time of all methods to $t=300$ seconds. For RB-PDA in deterministic setting we set $\sigma=0.01/n$ and $\tau$ is selected from $\{10^{-2},10^{-1},10^{0},10^{1},10^{2}\}$ based on the best performance of the algorithm. In the stochastic setting, for RB-PDA, SMD, and SMP algorithms, we set $\sigma^k=\sigma/\sqrt{k}$ and $\tau=\tau/\sqrt{k}$. After running each algorithm, the gap value in \eqref{s-gap} is calculated for $Z=\dom f \times \dom h$ using MOSEK in CVX. Moreover, we select the mini-batch size for the dual variable as $u=n$ since the cost of calculating $\grad_\by\Phi$ is low while we consider different mini-batch sizes for the primal variable. 
\begin{table}[h!]
    \centering
    \caption{Performance comparison among RB-PDA under different block partitioning strategies, SMD, and SMP for Experiment 1 (logistic loss). Tables from top to bottom: deterministic and stochastic settings, and left to right: \texttt{w7a}, \texttt{madelon}, and \texttt{mushroom} datasets.}
    \scriptsize
    \begin{tabular}{|c|c|c|c|c|}
    \hline
    \multicolumn{5}{|c|}{\texttt{w7a} -- deterministic} \\ \hline
    Method & $M$ & $N$ & Gap & \# iter. \\ \hline
    RB-PDA & $1$ & $1$ & 1.1e-2 & 10412\\ \hline
    RB-PDA & $1$ & $37$  & 1.1e-2 & 11760\\ \hline
    RB-PDA & $1$ & $634$  & 1.5e-1 & 11806\\ \hline
    RB-PDA & \bf 3 & \bf 1  & \bf 1.2e-3 & \bf 20597\\ \hline
    RB-PDA & $10$ & $1$  & 2.9e-2 & 24863\\ \hline
    RB-PDA & $3$ & $37$  & 1.3e-2 & 25956\\ \hline
    RB-PDA & $10$ & $37$ & 6.5e-2 & 33368\\ \hline
    SMD & $1$ & $1$  & 5.9e0 & 95835\\ \hline
    SMP & $1$ & $1$ & 6.0e0 & 47564\\ \hline
    \end{tabular}
    \hspace{3mm}
    \begin{tabular}{|c|c|c|c|c|}
    \hline
    \multicolumn{5}{|c|}{\texttt{madelon} -- deterministic} \\ \hline
    Method & $M$ & $N$ & Gap & \# iter. \\ \hline
    RB-PDA & 1 & 1 & 2.3e-1 & 3794 \\ \hline
    RB-PDA & \bf 10 & \bf 1 & \bf 5.7e-2 & \bf 18357 \\ \hline
    RB-PDA & 1 & 10 & 1.2e0 & 4177 \\ \hline
    RB-PDA & 10 & 10 & 2.4e-1 & 31129 \\ \hline
    RB-PDA & 1 & 100 & 7.9e0 & 4193 \\ \hline
    RB-PDA & 10 & 100 & 8.9e-1 & 32699 \\ \hline
    SMD & 1 & 1 & 6.4e-1 & 18745 \\ \hline
    SMP & 1 & 1 & 9.0e-1 & 10608 \\ \hline
    \end{tabular}
    \hspace{3mm}
    \begin{tabular}{|c|c|c|c|c|}
    \hline
    \multicolumn{5}{|c|}{\texttt{mushroom} -- deterministic} \\ \hline
    Method & $M$ & $N$ & Gap & \# iter. \\ \hline
    RB-PDA & 1 & 1 & 5.4e-3 & 1177 \\ \hline
    RB-PDA & \bf 8 & \bf 1 & \bf 1.2e-3 & \bf 5986 \\ \hline
    RB-PDA & 1 & 12 & 8.3e-3 & 7949 \\ \hline
    RB-PDA & 8 & 12 & 9.4e-3 & 9931 \\ \hline
    SMD & 1 & 1 & 9.4e-1 & 23745 \\ \hline
    SMP & 1 & 1 & 1.1e0 & 12608 \\ \hline
    \end{tabular}
    
    \vspace{3mm}
    
    \begin{tabular}{|c|c|c|c|c|c|}
    \hline
    \multicolumn{5}{|c|}{\texttt{w7a} -- stochastic} \\ \hline
    Method & $M$ & $N$ & Gap & \# iter. \\ \hline
    RB-PDA & $1$ & $1$ &  4e-3 & 55821\\ \hline
    RB-PDA & \bf 1 & \bf 37 &  \bf 1.5e-3 & \bf 145115\\ \hline
    RB-PDA & $1$ & $634$ & 7.5e-2 & 155808\\ \hline
    RB-PDA & $3$ & $1$ & 8.4e-3 & 56196\\ \hline
    RB-PDA & $10$ & $1$ & 1.4e-2 & 56396\\ \hline
    RB-PDA & \bf 3 & \bf 37 &  \bf 3.7e-3 & \bf 148395\\ \hline
    RB-PDA & $10$ & $37$ &  9.9e-3 & 148633\\ \hline
    SMD & $1$ & $1$ &  5.2e0 & 301182\\ \hline
    SMP & $1$ & $1$ &  6.3e0 & 152452\\ \hline
    \end{tabular}
    \hspace{3mm}
    \begin{tabular}{|c|c|c|c|c|c|}
    \hline
    \multicolumn{5}{|c|}{\texttt{madelon} -- stochastic} \\ \hline
    Method & $M$ & $N$ & Gap & \# iter. \\ \hline
    RB-PDA & 1 & 1 & 6.9e-2 & 15448 \\ \hline
    RB-PDA & \bf 10 & \bf 1 & \bf 1.6e-2 & \bf 16517 \\ \hline
    RB-PDA & 1 & 10 & 1.2e-1 & 23327 \\ \hline
    RB-PDA & 10 & 10 & 1.2e-1 & 25970 \\ \hline
    RB-PDA & 1 & 100 & 4.2e-1 & 22790 \\ \hline
    RB-PDA & 10 & 100 & 6.6e-1 & 26878 \\ \hline
    BMD & 1 & 1 & 6.1e-1 & 49403 \\ \hline
    SMP & 1 & 1 & 2.1e0 & 24906 \\ \hline
    \end{tabular}
    \hspace{3mm}
    \begin{tabular}{|c|c|c|c|c|c|}
    \hline
    \multicolumn{5}{|c|}{\texttt{mushroom} -- stochastic} \\ \hline
    Method & $M$ & $N$ & Gap & \# iter. \\ \hline
    RB-PDA & 1 & 1 & 7.4e-2 & 6419 \\ \hline
    RB-PDA & 8 & 1 & 2.6e-1 & 6738 \\ \hline
    RB-PDA & \bf 1 & \bf 12 & \bf 7.3e-2 & \bf 10749 \\ \hline
    RB-PDA & 8 & 12 & 2.5e-1 & 11088 \\ \hline
    SMD & 1 & 1 & 9.6e-1 & 23745 \\ \hline
    SMP & 1 & 1 & 1.5e0 & 12608 \\ \hline
    \end{tabular}

    \vspace*{0mm}
    \label{tab:numerics}
\end{table}

\begin{table}[h!]
    \centering
    \caption{Performance comparison among RB-PDA under different block partitioning strategies, SMD, and SMP for Experiment 2 (hinge loss). Tables from top to bottom: deterministic and stochastic settings, and left to right: \texttt{w7a}, \texttt{madelon}, and \texttt{mushroom} datasets.}
    \scriptsize
    \begin{tabular}{|c|c|c|c|c|}
    \hline
    \multicolumn{5}{|c|}{\texttt{w7a} -- deterministic} \\ \hline
    Method & $M$ & $N$ & Gap & \# iter. \\ \hline
    RB-PDA & 1 & 1 & 5.1e-1 & 11147 \\ \hline
    RB-PDA & \bf 10 & \bf 1 & \bf 3.1e-1 & \bf 45470 \\ \hline
    RB-PDA & 1 & 37 & 4.0e0 & 12553 \\ \hline
    RB-PDA & 10 & 37 & 1.5e0 & 91981 \\ \hline
    RB-PDA & 1 & 634 & 4.4e1 & 12508 \\ \hline
    RB-PDA & 10 & 634 & 1.7e1 & 96081 \\ \hline
    SMD & $1$ & $1$  & 5.1e0 & 95835\\ \hline
    SMP & $1$ & $1$ & 5.9e0 & 47564\\ \hline
    \end{tabular}
    \hspace{3mm}
    \begin{tabular}{|c|c|c|c|c|}
    \hline
    \multicolumn{5}{|c|}{\texttt{madelon} -- deterministic} \\ \hline
    Method & $M$ & $N$ & Gap & \# iter. \\ \hline
    RB-PDA & 1 & 1 & 5.5e-1 & 3445 \\ \hline
    RB-PDA & 10 & 1 & 1.8e-1 & 16880 \\ \hline
    RB-PDA & 1 & 10 & 5.5e-1 & 3819 \\ \hline
    RB-PDA & \bf 10 & \bf 10 & \bf 8.5e-2 & \bf 27866 \\ \hline
    RB-PDA & 1 & 100 & 1.5e+1 & 3853 \\ \hline
    RB-PDA & 10 & 100 & 1.8e-1 & 30743 \\ \hline
    SMD & 1 & 1 & 9.5e-1 & 40055  \\ \hline
    SMP & 1 & 1 & 1.5e0 & 18351  \\ \hline
    \end{tabular}
    \hspace{3mm}
    \begin{tabular}{|c|c|c|c|c|}
    \hline
    \multicolumn{5}{|c|}{\texttt{mushroom} -- deterministic} \\ \hline
    Method & $M$ & $N$ & Gap & \# iter. \\ \hline
    RB-PDA & 1 & 1 & 5.4e-3 & 1177 \\ \hline
    RB-PDA & \bf 8 & \bf 1 & \bf 1.2e-3 & \bf 5986 \\ \hline
    RB-PDA & 1 & 12 & 8.3e-3 & 7949 \\ \hline
    RB-PDA & 8 & 12 & 9.4e-3 & 9931 \\ \hline
    SMD & 1 & 1 & 9.4e-1 & 21745 \\ \hline
    SMP & 1 & 1 & 1.2e0 & 10608 \\ \hline
    \end{tabular}
    
    \vspace{3mm}
    
    \begin{tabular}{|c|c|c|c|c|c|}
    \hline
    \multicolumn{5}{|c|}{\texttt{w7a} -- stochastic} \\ \hline
    Method & $M$ & $N$ & Gap & \# iter. \\ \hline
    RB-PDA & 1 & 1 & 3.5e-1 & 54752 \\ \hline
    RB-PDA & \bf 10 & \bf 1 & \bf 2.5e-1 & \bf 53001 \\ \hline
    RB-PDA & 1 & 37 & 4.0e0 & 132815 \\ \hline
    RB-PDA & 10 & 37 & 6.8e0 & 135122 \\ \hline
    RB-PDA & 1 & 634 & 1.4e2 & 134156 \\ \hline
    RB-PDA & 10 & 634 & 9.3e1 & 136903 \\ \hline
    SMD & $1$ & $1$ &  5.2e0 & 301182\\ \hline
    SMP & $1$ & $1$ &  6.3e0 & 152452\\ \hline
    \end{tabular}
    \hspace{3mm}
    \begin{tabular}{|c|c|c|c|c|c|}
    \hline
    \multicolumn{5}{|c|}{\texttt{madelon}  -- stochastic} \\ \hline
    Method & $M$ & $N$ & Gap & \# iter. \\ \hline
    RB-PDA & 1 & 1 & 7.7e-1 & 13302 \\ \hline
    RB-PDA & 10 & 1 & 2.1e-1 & 14292 \\ \hline
    RB-PDA & 1 & 10 & 7.2e-1 & 20348 \\ \hline
    RB-PDA & \bf 10 & \bf 10 & \bf 2.0e-1 & \bf 21939 \\ \hline
    RB-PDA & 1 & 100 & 5.1e0 & 24439 \\ \hline
    RB-PDA & 10 & 100 & 2.4e0 & 23199 \\ \hline
    BMD & 1 & 1 & 5.2e-1 & 50069 \\ \hline
    SMP & 1 & 1 & 4.4e-1 & 25356 \\ \hline
    \end{tabular}
    \hspace{3mm}
    \begin{tabular}{|c|c|c|c|c|c|}
    \hline
    \multicolumn{5}{|c|}{\texttt{mushroom}  -- stochastic} \\ \hline
    Method & $M$ & $N$ & Gap & \# iter. \\ \hline
    RB-PDA & 1 & 1 & 7.4e-2 & 6419 \\ \hline
    RB-PDA & 8 & 1 & 2.6e-1 & 6738 \\ \hline
    RB-PDA & \bf 1 & \bf 12 & \bf 7.3e-2 & \bf 10749 \\ \hline
    RB-PDA & 8 & 12 & 2.5e-01 & 11088 \\ \hline
    SMD & 1 & 1 & 9.6e-1 & 23745 \\ \hline
    SMP & 1 & 1 & 1.5e0 & 12608 \\ \hline
    \end{tabular}

    \vspace*{0mm}
    \label{tab:numerics2}
\end{table}

\rev{\noindent\textbf{Experiment 1.} In the first experiment, we use the logistic loss $L_j(\mathbf{x}) = \log(1 + \exp(-b_j a_j^\top \mathbf{x}))$ to evaluate the performance of RB-PDA across different primal-dual block configurations $(M, N)$ on both datasets. The results in Table \ref{tab:numerics} present the performance of methods in terms of the gap function within a fixed time budget. The results demonstrate that RB-PDA consistently achieves smaller optimality gaps than SMD and SMP within a fixed time budget that is same for all the methods. On \texttt{w7a}, moderate block sizes such as $(M,N)=(3,1)$ in deterministic setting and $(1,37)$ in stochastic setting yield the smallest gaps ($1.2\times 10^{-3}$ and $1.5\times 10^{-3}$, respectively), whereas extremely fine-grained partitions ($N=634$) degrade performance. On \texttt{madelon} and \texttt{mushroom}, block sizes $(M,N)=(8,1)$ and $(M,N)=(10,1)$ achieve similarly low gaps ($\approx 5.6\times 10^{-2}\text{ and } 1.2\times 10^{-3}$) in the deterministic setting, repsectively, far outperforming SMD and SMP whose gaps remain in the order of 1. Moreover, a more extreme partitioning $(M,N)=(8,12)$ in \texttt{mushroom} has a comparable performance with a cheaper per-iteration cost reflected in a larger number of iterations performed in the fixed time budget. Similar performance characteristics can also be observed in the stochastic setting. These findings confirm that carefully chosen block partitions enhance convergence speed while excessive partitioning can be detrimental to the practical performance.}

\rev{\noindent\textbf{Experiment 2.} In the second experiment, we replace the logistic loss with the squared hinge loss $L_j(\mathbf{x}) = \big(1 - b_j a_j^\top \mathbf{x}\big)_+^2$ and repeat the tests on the same datasets and block configurations. The results as presented in Table \ref{tab:numerics2} show that RB-PDA retains its performance advantage regardless of the loss function. On \texttt{w7a}, deterministic RB-PDA with $(M,N)=(1,10)$ achieves a gap of $3.1\times 10^{-1}$, considerably better than the gaps of SMD and SMP, and similar trends are observed in the stochastic setting. On \texttt{madelon} and \texttt{mushroom}, RB-PDA with $(M,N)=(10,10)$ and $(M,N)=(8,1)$ in a deterministic setting achieves the smallest gap ($8.5\times 10^{-2}$ and $1.2\times 10^{-3}$), respectively, outperforming other configurations and surpassing SMD/SMP. Overall, RB-PDA demonstrates robustness to changes in loss functions, consistently delivering lower gaps and confirming the efficacy of block-coordinate schemes for large-scale binary classification tasks.}

\section{Conclusions} \label{sec:con}
\nsa{We studied a randomized block coordinate primal-dual 
method to solve the large-scale 
SP problems.} {The method can contend with non-bilinear non-separable maps $\Phi(\bx,\by)$, possibly with a stochastic objective function as well as multiple primal and/or dual blocks.}
\us{At} each iteration, a pair of primal and dual blocks is randomly selected and updated. In the stochastic setting, we utilize a {\em mini-batch} of 
\saa{stochastic partial} gradients in each iteration to approximate the partial gradient maps. 
Under this setup, \rev{for deterministic setting, it is shown that the iterates converge to a saddle point almost surely with an ergodic convergence rate {of} ${\cO}(1/k)$ in terms of the expected gap function; moreover, for the stochastic setting, an ergodic convergence rate {of} $\tilde{\cO}(1/\sqrt{k})$ has been shown. These results are valid whether the problem domain is bounded or unbounded.}
\aj{In addition, the oracle complexity for both deterministic and stochastic settings and their dependency on the number of primal and dual blocks are obtained.} 
Empirical experiments suggest that the scheme compares well with existing techniques.

\rev{A potential future direction would be on double-loop stochastic-coordinate methods for solving for large-scale non-convex saddle-point problems and non-convex constrained optimization. Indeed, for the nonconvex-(strongly) concave setting when there is a single coordinate block for both primal and dual variables, in \citep{zhang2022sapd+} we proposed an ineaxact proximal-point method SAPD+ which repeatedly calls the stochastic primal-dual method we proposed in~\citep{zhang2024robust} for inexactly solving stochastic strongly convex-strongly concave saddle point problems -- the main proof technique here is to carefully control the change in expected duality gap evaluated at the inexact solutions to each saddle-point subproblem. The technique we designed in \citep{zhang2022sapd+} can be combined with the block-coordinate method RB-PDA we propose in this paper.}



\vspace*{0mm}





\section*{Appendix}\label{sec:append}

 \section{Proofs of Supporting Lemmas}
 \label{sec:proof-supporting}
 \subsection{Proof of Lemma \ref{lem:inner-w}.} \label{sec:proof-prelim}
 Since $\bv^{k+1}$-update is separable in 
 $i\in\cM$, one can apply Lemma \ref{lem_app:prox} for each coordinate 
 to obtain a bound for $\fprod{\bdelta^k,\bx-\bv^{k+1}}$. Then we have that\vspace*{0mm}
 {\small
 \begin{align*}
 \fprod{\bdelta^k,\bx-\bv^k}&=\fprod{\bdelta^k,\bx-\bv^{k+1}}+\fprod{\bdelta^k,\bv^{k+1}-\bv^k}\nonumber\\
 &\leq \bD_\cX^{\cA}(\bx,\bv^k)-\bD_\cX^{\cA}(\bx,\bv^{k+1})-\bD_\cX^{\cA}(\bv^{k+1},\bv^k)+\fprod{\bdelta^k,\bv^{k+1}-\bv^k}\nonumber\\
 &\leq \bD_\cX^{\cA}(\bx,\bv^k)-\bD_\cX^{\cA}(\bx,\bv^{k+1})+\frac{1}{2}\norm{\bdelta^k}_{\nsa{*,\cA^{-1}}}^2,
 \end{align*}}%
 where in the last inequality we used 
 \nsa{$\fprod{\delta^k_i,v_i^{k+1}-v^k_i}\leq\tfrac{a_i}{2}\norm{v_i^{k+1}-v^k_i}_{\cX_i}^2+\tfrac{1}{2a_i}\norm{\delta_i^k}_{\cX_i^*}^2$ for $i\in\cM$ together with $\sp{\bD_\cX^\cA(\bv^{k+1},\bv^k)}\geq \tfrac{1}{2}\norm{\bv^{k+1}-\bv^k}_\cA^2$.}
 \qed
 \subsection{Proof of Lemma \ref{lem:one-step}.}\label{sec:proof-lem-one-step}
 For $k\geq0$, 
 Lemma~\ref{lem_app:prox} 
 \nsa{applied to 
  the $\tilde{y}_j$-subproblem in \eqref{eq:tilde-problem-y} and the $\tilde{x}_i$-subproblem in \eqref{eq:tilde-problem-x},
 implies that for \sa{any $y_j\in\cY_j$, $j\in\cN$, and any $x_i\in\cX_i$, $i\in\cM$}:} 
 {\small
 \begin{align}
  h_j(\tilde y_j^{k+1})-h_j(y_j)-\fprod{\tilde s_j^k+\ey{Nw_{y_j}^k},~ \tilde y_j^{k+1}-y_j} &\leq \frac{1}{\sigma^k_j}\Big[\bD_{\cY_j}(y_j,y_j^k) -\bD_{\cY_j}(y_j,\tilde{y}_j^{k+1})-\bD_{\cY_j}(\tilde{y}_j^{k+1},y_j^k)\Big], \label{eq:sc-h}\\
 f_i(\tilde{x}_i^{k+1})- f_i(x_i)+\fprod{\tilde{r}_i^k+\ey{M w_{x_i}^k},\tilde{x}_i^{k+1}-x_i} &\leq \frac{1}{\tau^k_i} \Big[\bD_{\cX_i}(x_i,x_i^k) -\bD_{\cX_i}(x_i,\tilde{x}_i^{k+1})-\bD_{\cX_i}(\tilde{x}_i^{k+1},x_i^k)\Big].  \label{eq:sc-f}
 \end{align}}%
 We define two auxiliary sequences \sa{by summing the Bregman terms in \eqref{eq:sc-h} over $j\in\cN$ and those in \eqref{eq:sc-f} over $i\in\cM$, i.e., for $k\geq 0$,}
 \sp{\small
 \begin{align*}
     A_1^{k+1}\triangleq \bD_{\cY}^{\bS^k}(\by,\by^k) -\bD_{\cY}^{\bS^k}(\by,\tilde{\by}^{k+1})-\bD_{\cY}^{\bS^k}(\tilde{\by}^{k+1},\by^k),\quad
     A_2^{k+1}\triangleq \bD_{\cX}^{\bT^k}(\bx,\bx^k) -\bD_{\cX}^{\bT^k}(\bx,\tilde{\bx}^{k+1})-\bD_{\cX}^{\bT^k}(\tilde{\bx}^{k+1},\bx^k).
 \end{align*}}%
 \sa{For $k\geq 0$, using} coordinate-wise Lipschitz continuity of \sa{$\grad_\by\Phi(\bx^k,\cdot)$}, we obtain that \sa{for $j\in\cN$,}
 {\small
 \begin{align}\label{eq:lip-yy}
 \MoveEqLeft\fprod{\grad_{y_j}\Phi(\bx^k,~\by^k),~\tilde y_j^{k+1}-y_j}=\fprod{\grad_{y_j}\Phi(\bx^k,\by^k),~\tilde y_j^{k+1}-y_j^k}+\fprod{\grad_{y_j}\Phi(\bx^k,\by^k),~y_j^k-y_j} \\ \nonumber
 &\leq \Phi(\bx^k,\by^k+V_j(\tilde y_j^{k+1}-y_j^k))-\Phi(\bx^k,\by^k)+\frac{L_{y_jy_j}}{2}\norm{\tilde y_j^{k+1}-y_j^k}_{\cY_j}^2+\fprod{\grad_{y_j}\Phi(\bx^k,\by^k),~y_j^k-y_j}.
 \end{align}}%
 {Summing \eqref{eq:sc-h} over $j\in\cN$, and multiplying by $\frac{1}{N}$ gives 
 \sa{$h(\tilde{\by}^{k+1})-h(\by)-\tfrac{1}{N}\fprod{\ey{\tilde\bs^k+N\bw_\by^k},~ \tilde{\by}^{k+1}-\by}$ $\leq \tfrac{1}{N} A_1^{k+1}$}.
 Now, using \eqref{eq:lip-yy} and the fact that $\tilde s_j^k=N\grad_{y_{j}}\Phi(\bx^k,\by^k)+\sp{\tilde q_{y_j}^k}$, we obtain}
 {\small
 \begin{align}
 &h(\bty^{k+1})-h(\by)- \sum_{j\in\cN}\Phi(\bx^k,\by^k+V_j(\tilde y_j^{k+1}-y_j^k)) \label{eq:y-ineq}\\
 & \leq-N\Phi(\bx^k,\by^k)+\fprod{\grad_\by\Phi(\bx^k,\by^k),~\by^k-\by}+\tfrac{1}{N}\fprod{\sp{\tilde\bq_{\by}^k}+\ey{N\bw_\by^k},~\bty^{k+1}-\by}\nonumber +\tfrac{1}{2}\norm{\bty^{k+1}-\by^k}_{\bL_{\by\by}}^2+\tfrac{1}{N}A_1^{k+1} \nonumber\\ 
 &\leq-(N-1)\Phi(\bx^k,\by^k)-\Phi(\bx^k,\by)+\tfrac{1}{N}\fprod{\sp{\tilde\bq_{\by}^k}+\ey{N\bw_\by^k},~\bty^{k+1}-\by}\nonumber +\tfrac{1}{2}\norm{\bty^{k+1}-\by^k}_{\bL_{\by\by}}^2+\tfrac{1}{N}A_1^{k+1},\nonumber
 \end{align}%
 where \us{$\bL_{\by \by}$ is defined in \aj{\eqref{def L}}} and in the last inequality, we invoke concavity of $\Phi(\bx^k,\cdot)$ \sa{as $\bx^k\in\dom f$}. \sa{Also note that for $k\geq 0$ and $i\in\cM$,} 
 the coordinate-wise Lipschitz continuity of \sa{$\grad_\bx\Phi(\cdot,\by^{k+1})$} implies 
 \begin{align}
 \fprod{\grad_{x_i}\Phi(\bx^k,\by^{k+1}),\tilde x_i^{k+1}-x_i} & \geq \Phi(\bx^k+U_i(\tilde x_i^{k+1}-x_i^k),\by^{k+1})-\Phi(\bx^k,\by^{k+1})\label{eq:lip-g} \\
  &\quad -\frac{L_{x_ix_i}}{2}\norm{\tilde x_i^{k+1}- x_i^k}^2_{\cX_i}+\fprod{\grad_{x_i}\Phi(\bx^k,\by^{k+1}),~x_i^{k}-x_i}.\nonumber
 \end{align}}%
 Similar to the derivation of \eqref{eq:y-ineq}, using \eqref{eq:sc-f} and \eqref{eq:lip-g}, 
 leads to
 {\small
 \begin{align}
 \MoveEqLeft f(\btx^{k+1})+ \sum_{i\in\cM}\Phi(\bx^k+U_i(\tilde x_i^{k+1}-x_i^k),\by^{k+1})- M\Phi(\bx^k,\by^{k+1}) \label{eq:convex-Lip-f}\\ 
 &\nonumber \leq f(\bx) -\fprod{\grad_\bx\Phi(\bx^k,\by^{k+1}),\bx^k-\bx}+\frac{1}{2}\norm{\btx^{k+1}- \bx^k}^2_{\bL_{\bx\bx}}+\frac{1}{M}A_2^{k+1}-\frac{1}{M}\fprod{\sp{\tilde\bq_\bx^k}+M\bw_\bx^k,~\btx^{k+1}-\bx} \nonumber\\ 
 &\nonumber\leq f(\bx)+\Phi(\bx,\by^{k+1})-\Phi(\bx^k,\by^{k+1})+\frac{1}{2}\norm{\btx^{k+1}- \bx^k}^2_{\bL_{\bx\bx}}+\frac{1}{M}A_2^{k+1}-\frac{1}{M}\fprod{\sp{\tilde\bq_\bx^k}+M\bw_\bx^k,~\btx^{k+1}-\bx}, \nonumber
 \end{align}}%
 where \sa{$\bL_{\bx \bx}$ is defined in \eqref{def L},} and the last inequality \sa{follows from the convexity of $\Phi(\cdot,\by^{k+1})$.} 
 \rev{Rearranging the terms in \eqref{eq:y-ineq} and adding and subtracting $N(h(\by^{k+1})-h(\by))$, $(N-1)h(\by^k)$ and $N\Phi(\bx^k,\by^{k+1})$, we get 
 {\small
 \begin{align}\nonumber
 N(h(\by^{k+1})-h(\by))&\leq N(h(\by^{k+1})-h(\by))+h(\by)-h(\tilde \by^{k+1})+(N-1)h(\by^k)-(N-1)h(\by^k)+N\Phi(\bx^k,\by^{k+1})\nonumber\\ &\quad-N\Phi(\bx^k,\by^{k+1})+\sum_{j\in\cN}\Phi(\bx^k,\by^k+V_j(\tilde y_j^{k+1}-y_j^k))-(N-1)h(\by^k)+Nh(\by^{k+1})-h(\bty^{k+1})\nonumber\\ 
 &\quad -\Phi(\bx^k,\by)+\tfrac{1}{N}A_1^{k+1}+\tfrac{1}{2}\norm{\bty^{k+1}-\by^k}^2_{\bL_{\by\by}}+\tfrac{1}{N}\fprod{{\tilde\bq_\by^k}+{N\bw_\by^k},~\bty^{k+1}-\by}\nonumber\\
 &\leq (N-1)\left(h(\by^k)-h(\by)\right)+N\Phi(\bx^k,\by^{k+1})-(N-1)\Phi(\bx^k,\by^k)-\Phi(\bx^k,\by)\label{eq:expect-h}\nonumber\\
 &\quad+\tfrac{1}{N}A_1^{k+1}+\tfrac{1}{2}\norm{\bty^{k+1}-\by^k}^2_{\bL_{\by\by}}+\tfrac{1}{N}\fprod{{\tilde\bq_\by^k}+{N\bw_\by^k},~\bty^{k+1}-\by}+E_1,
 \end{align}}%
 {where $E_1 \triangleq \sum_{j\in\cN}\Phi(\bx^k,\by^k+V_j(\tilde y_j^{k+1}-y_j^k))-N\Phi(\bx^k,\by^{k+1})-(N-1)h(\by^k)+Nh(\by^{k+1})-h(\bty^{k+1})$.}}
 {Similarly, by adding and subtracting $M(f(\bx^{k+1})-f(\bx))$ and $M\Phi(\bx^{k+1},\by^{k+1})$ to \eqref{eq:convex-Lip-f}, 
 we get}
 {\small
 \begin{eqnarray}
 \lefteqn{M(f(\bx^{k+1})-f(\bx)+\Phi(\bx^{k+1},\by^{k+1})) -\Phi(\bx,\by^{k+1})+\tfrac{1}{M}\fprod{\sp{\tilde\bq_\bx^k}+M\bw_\bx^k,\btx^{k+1}-\bx}} \label{eq:expect-f} \\
 &&\leq\left(M-1\right)\big(f(\bx^k)-f(\bx)+\Phi(\bx^k,\by^{k+1})\big)+
 \tfrac{1}{2}\norm{\btx^{k+1}-\bx^k}_{\bL_{\bx\bx}}^2+\tfrac{1}{M}A_2^{k+1}+E_2,\nonumber
 \end{eqnarray}}%
 where $E_2\triangleq  Mf(\bx^{k+1})-(M-1)f(\bx^k)-f(\btx^{k+1})-\sum_{i\in\cM}\Phi(\bx^k+U_i(\tilde{x}_i^{k+1}-x_i^k),\by^{k+1})+M\Phi(\bx^{k+1},\by^{k+1})$.
 Summing 
 \eqref{eq:expect-h} and \eqref{eq:expect-f}, adding \eyh{$\Phi(\bx^{k+1},\by)$} to both sides, and rearranging terms lead to 
 {\small
 \begin{eqnarray}
 \lefteqn{f(\bx^{k+1}) -f(\bx) +\Phi(\bx^{k+1},\by)+h(\by^{k+1})-h(\by) -\Phi(\bx,\by^{k+1})}\label{eq:sum-f-h}\\ 
 & & \leq M(\Phi(\bx^{k+1},\by)-\Phi(\bx^{k+1},\by^{k+1})) +M(\Phi(\bx^k,\by^{k+1})-\Phi(\bx^k,\by)) \nonumber \\
 & & \hbox{ } +\tfrac{1}{2}\norm{\btx^{k+1}-\bx^k}_{\bL_{\bx\bx}}^2+\tfrac{1}{2}\norm{\bty^{k+1}-\by^k}_{\bL_{\by\by}}^2+ \left(M-1\right)\big(f(\bx^k)-f(\bx)+\Phi(\bx^k,\by)\big) \nonumber \\
 & & \hbox{ } -\left(M-1\right) 
 \nsa{\left(f(\bx^{k+1}) -f(\bx) +\Phi(\bx^{k+1},\by)\right)}+(N-1)\left(h(\by^k)-h(\by)-\Phi(\bx^k,\by^k)\right)\nonumber \\
 & & \hbox{ } -(N-1)
 \nsa{\left(h(\by^{k+1})-h(\by)-\Phi(\bx^k,\by^{k+1})\right)}+\tfrac{1}{N}A_1^{k+1}+\tfrac{1}{M}A_2^{k+1}\nonumber \\
 & & \hbox{ } +\tfrac{1}{N}\fprod{\sp{\tilde\bq_\by^k}+\ey{N\bw_\by^k},\bty^{k+1}-\by}-\tfrac{1}{M}\fprod{\sp{\tilde\bq_\bx^k}+M\bw_\bx^k,\btx^{k+1}-\bx}+E_1+E_2\nonumber \\
 & & \leq M\big(G_\Phi^{k+1}(\bz)-G_\Phi^k(\bz)\big)+M\fprod{\grad_\by\Phi(\bx^k,\by^k),\by^{k+1}-\by^k}\nonumber\\
 & & \hbox{ } +\tfrac{1}{2}\norm{\btx^{k+1}-\bx^k}_{\bL_{\bx\bx}}^2+\tfrac{1}{2}\norm{\bty^{k+1}-\by^k}_{\bL_{\by\by}}^2 +(M-1)\big(G_f^k(\bz)-G_f^{k+1}(\bz)\big)\nonumber \\
 & & \hbox{ } +(N-1)(G_h^k(\bz)-G_h^{k+1}(\bz))+(N-1)\fprod{\grad_\bx\Phi(\bx^k,\by^{k+1}),\bx^k-\bx^{k+1}}\nonumber\\
 & & \hbox{ } +\tfrac{1}{N}A_1^{k+1}+\tfrac{1}{M}A_2^{k+1}+\tfrac{1}{N}\fprod{\sp{\tilde\bq_\by^k}+\ey{N\bw_\by^k},\bty^{k+1}-\by}-\tfrac{1}{M}\fprod{\sp{\tilde\bq_\bx^k}+M\bw_\bx^k,\btx^{k+1}-\bx}+E_1+E_2, \nonumber
 \end{eqnarray}}%
 where in the last inequality, concavity of $\Phi(\bx^k,\cdot)$, convexity of $\Phi(\cdot,\by^{k+1})$, and \rev{the definitions of $G_f^k(\bz)$, $G_h^k(\bz)$, and $G_\Phi^k(\bz)$ as given in Definition \ref{def:pos-set} are used.}
 Moreover, using 
 \eqref{eq:G-connection}, the definitions of $G_f^k$ and $\cL(\cdot,\cdot)$, \sa{
the concavity of $\Phi(\bx^k,\cdot)$, and \rev{the fact that we chose $\theta^k\geq 1$ in Algorithm \ref{alg:RB-PDA}}, we obtain:}\vspace*{0mm}
 {\small
 \begin{eqnarray}
 \lefteqn{(1-\theta^k)[(M-1)G_f^k(\bz)+(N-1)G_h^k(\bz)-M G_\Phi^k(\bz)]}  \label{sum-ghf}\\
 & & = (1-\theta^k)[(M-1)\big(\cL(\bx^k,\by)-\cL(\bx,\by)\big)+(N-1)G_h^k(\bz)-M G_\Phi^k(\bz)] \nonumber\\
 & & \leq (1-\theta^k)\Big[(M-1)\big(\cL(\bx^k,\by)-\cL(\bx,\by)\big)+(N-1)\big(\cL(\bx,\by)-\cL(\bx,\by^k)\big)\Big]\nonumber\\
 & & \hbox{ }  +(1-\theta^k)\Big[(N-1)\fprod{\grad_\bx\Phi(\bx^k,\by^k),\bx-\bx^k}-M\fprod{\grad_\by\Phi(\bx^k,\by^k),\by-\by^k}\Big].\nonumber
 \end{eqnarray}}%
 Therefore, utilizing \eqref{sum-ghf} within \eqref{eq:sum-f-h} and rearranging the terms leads to
 {\small
 \begin{align}
 & \cL(\bx^{k+1},\by)-\cL(\bx,\by^{k+1})
+(\theta^k-1)\Big[(M-1)\big(\cL(\bx^k,\by)-\cL(\bx,\by)\big)+(N-1)\big(\cL(\bx,\by)-\cL(\bx,\by^k)\big)\Big] \label{eq:sum-f-h-bound}\\
 & \leq M\fprod{\grad_\by\Phi(\bx^k,\by^k),\by^{k+1}-\by^k}+\tfrac{1}{2}\norm{\btx^{k+1}-\bx^k}_{\bL_{\bx\bx}}^2+\tfrac{1}{2}\norm{\bty^{k+1}-\by^k}_{\bL_{\by\by}}^2 \nonumber\\
 & \hbox{ } +(M-1)\big(\theta^kG_f^k(\bz)-G_f^{k+1}(\bz)\big)+(N-1)\big(\theta^kG_h^k(\bz)-G_h^{k+1}(\bz)\big)+M\big(G_\Phi^{k+1}(\bz)-\theta^k G_\Phi^k(\bz)\big)\nonumber\\
 & \hbox{ } +(N-1)\fprod{\grad_\bx\Phi(\bx^k,\by^{k+1}),\bx^k-\bx^{k+1}}+M(\theta^k-1)\fprod{\grad_\by\Phi(\bx^k,\by^k),\by-\by^k} \nonumber\\
 & \hbox{ } +(N-1)(1-\theta^k)\fprod{\grad_\bx\Phi(\bx^k,\by^k),\bx-\bx^k}  +\tfrac{1}{N}A_1^{k+1}+\tfrac{1}{M}A_2^{k+1}+E_1+E_2 \nonumber\\ 
 & \hbox{ } +\tfrac{1}{N}\fprod{\sp{\tilde\bq_\by^k},\bty^{k+1}-\by}+\nsa{\tfrac{1}{M}\fprod{\sp{\tilde\bq_\bx^k},\bx-\btx^{k+1}}+\ey{\fprod{\bw_\by^k,\bty^{k+1}-\by}+\fprod{\bw_\bx^k,\bx-\btx^{k+1}}}}. \nonumber
 \end{align}}%
Recall that 
 $\sp{\tilde\bq_{\by}^k} \triangleq M N\theta^k(\grad_\by\Phi(\bx^k,\by^k)-\grad_\by\Phi(\bx^{k-1},\by^{k-1}))$ for all $k\geq 0$. Now we rearrange the inner products containing $\grad_\by\Phi$ on the right hand side of \eqref{eq:sum-f-h-bound}: 
 {\small
 \begin{eqnarray}
 \label{eq:arrange-inner-y}
 \lefteqn{M\fprod{\grad_\by\Phi(\bx^{k},\by^{k}),\by^{k+1}-\by^{k}}-M(1-\theta^k)\fprod{\grad_\by\Phi(\bx^k,\by^k),\by-\by^k}+\tfrac{1}{N}\fprod{\sp{\tilde\bq_\by^k},\bty^{k+1}-\by}}  \\ 
 & & = M\fprod{\grad_\by\Phi(\bx^{k},\by^{k}),\by^{k+1}-\by)}+M\theta^k\fprod{\grad_\by\Phi(\bx^k,\by^k),\by-\by^k} +\fprod{\sp{\tilde\bq_\by^k},\by^{k+1}-\by^k}\nonumber \\
 && \quad {+\tfrac{1}{N}\fprod{\sp{\tilde\bq_{\by}^k},\by^{k}-\by}+\tfrac{1}{N}\fprod{\sp{\tilde\bq_\by^k},\bty^{k+1}\nsa{-}N\by^{k+1}\nsa{+}(N-1)\by^k}}\nonumber\\
 & & = M\fprod{\grad_\by\Phi(\bx^k,\by^k),\by^{k+1}-\by}-M\theta^k\fprod{\grad_\by\Phi(\bx^{k-1},\by^{k-1}),\by^k-\by}+\fprod{\sp{\tilde\bq_{\by}^k},\by^{k+1}-\by^k}\nonumber\\ 
 && \quad  {+\tfrac{1}{N}\fprod{\sp{\tilde\bq_\by^k},\bty^{k+1}-N\by^{k+1}+(N-1)\by^k}}. \nonumber
 \end{eqnarray}}%
 \nsa{Similarly, rearranging the inner products involving $\grad_\bx\Phi$ on the right hand side of \eqref{eq:sum-f-h-bound} and using $\sp{\tilde\bq_\bx^k}\triangleq M(N-1)\theta^k(\grad_\bx\Phi(\bx^{k},\by^{k})-\grad_\bx\Phi(\bx^{k-1},\by^{k-1}))$, we get}
 {\small
 \begin{align}
 \label{eq:arrange-inner-x}
 &(N-1)\fprod{\grad_\bx\Phi(\bx^k,\by^{k+1}),\bx^k-\bx^{k+1}}+(N-1)(1-\theta^k)\fprod{\grad_\bx\Phi(\bx^k,\by^k),\bx-\bx^k} +\tfrac{1}{M}\fprod{\sp{\tilde\bq_\bx^k},\bx-\btx^{k+1}}\nonumber\\
 & =(N-1)\fprod{\grad_\bx\Phi(\bx^k,\by^{k+1}),\bx^k-\bx^{k+1}}+(N-1)(1-\theta^k)\fprod{\grad_\bx\Phi(\bx^k,\by^k),\bx-\bx^k}\nonumber\\
 & \quad +\tfrac{1}{M}\fprod{\sp{\tilde\bq_\bx^k},\bx-\bx^k}+\fprod{\sp{\tilde\bq_\bx^k},\bx^k-\bx^{k+1}}+\tfrac{1}{M}\fprod{\sp{\tilde\bq_\bx^k},M\bx^{k+1}-(M-1)\bx^k-\btx^{k+1}}\nonumber\\
 & =(N-1)\fprod{\grad_\bx\Phi(\bx^k,\by^k),\bx-\bx^{k+1}}-(N-1)\theta^k\fprod{\grad_\bx\Phi(\bx^{k-1},\by^{k-1}),\bx-\bx^k} \nonumber\\
 & \quad +(N-1)\fprod{\grad_\bx\Phi(\bx^k,\by^{k+1})-\grad_\bx\Phi(\bx^k,\by^{k}),~\bx^k-\bx^{k+1}} +\fprod{\sp{\tilde\bq_\bx^k},\bx^k-\bx^{k+1}}\nonumber\\
 & \quad {+\tfrac{1}{M}\fprod{\sp{\tilde\bq_\bx^k},M\bx^{k+1}-(M-1)\bx^k-\btx^{k+1}}}. 
 \end{align}}%
 \nsa{Next}, substituting the results derived in \eqref{eq:arrange-inner-y} and \eqref{eq:arrange-inner-x} back into \eqref{eq:sum-f-h-bound}, 
 we get the following bound:
 {\small
 \begin{align}
 &\cL(\bx^{k+1},\by)-\cL(\bx,\by^{k+1})+(\theta^k-1)\big[(M-1)\big(\cL(\bx^k,\by)-\cL(\bx,\by)\big)+(N-1)\big(\cL(\bx,\by)-\cL(\bx,\by^k)\big)\big] \label{eq:pre_bound} \\
 &  \leq (M-1)\big(\theta^kG_f^k(\bz)-G_f^{k+1}(\bz)\big)+(N-1)\big(\theta^kG^k_h(\bz)-G_h^{k+1}(\bz)\big)+\underbrace{\fprod{\sp{\tilde\bq_{\by}^k},~\by^{k+1}-\by^k}}_{(a)}\nonumber\\
 &  \hbox{ } + \underbrace{\fprod{\sp{\tilde\bq_{\bx}^k},~\bx^{k}-\bx^{k+1}}}_{(b)}+\underbrace{(N-1)\fprod{\grad_\bx\Phi(\bx^k,\by^{k+1})-\grad_\bx\Phi(\bx^k,\by^k),~\bx^k-\bx^{k+1}}}_{(c)} \sa{+}\underbrace{\fprod{\bw_\bx^k,~\sa{\bx-\btx^{k+1}}}}_{(d)}+\underbrace{\ey{\fprod{\bw_\by^k,\bty^{k+1}-\by}}}_{(e)}\nonumber\\
 &  \hbox{ } + M\fprod{\grad_\by\Phi(\bx^k,\by^k),~\by^{k+1}-\by}-M\theta^k\fprod{\grad_\by\Phi(\bx^{k-1},\by^{k-1}),\by^k-\by} \nonumber\\
 &  \hbox{ } +(N-1)\fprod{\grad_\bx\Phi(\bx^k,\by^{k}),~\bx-\bx^{k+1}}-(N-1)\theta^k\fprod{\grad_\bx\Phi(\bx^{k-1},\by^{k-1}),\bx-\bx^k}\nonumber\\
 &  \hbox{ } +\tfrac{1}{2}\norm{\btx^{k+1}-\bx^k}^2_{\bL_{\bx\bx}} +\tfrac{1}{N}A_1^{k+1}+\tfrac{1}{M}A_2^{k+1}+\tfrac{1}{2}\norm{\bty^{k+1}-\by^k}^2_{\bL_{\by\by}}+ M \big(G_\Phi^{k+1}(\bz)-\theta^kG_\Phi^k(\bz)\big)\nonumber \\
 & \ \hbox{ } +\frac{1}{N}\fprod{\sp{\tilde\bq_\by^k},\bty^{k+1}-N\by^{k+1}+(N-1)\by^k}{+\frac{1}{M}\fprod{\sp{\tilde\bq_\bx^k},M\bx^{k+1}-(M-1)\bx^k-\btx^{k+1}}}+E_1+E_2. \nonumber
 \end{align}}%
 \nsa{We now} provide \us{bounds} for $(a)$, $(b)$, $(c)$ \nsa{and $(d)$} using the fact that for any \sa{$\by'\in\cY^*$}, $\by\in\cY$ and 
 $\bar{\theta}>0$, \nsa{$|\fprod{\by',\by}|\leq\frac{1}{2}\bar{\theta}\norm{\by'}_{\cY^*}^2+\frac{1}{2}\norm{\by}_{\cY}^2/\bar{\theta}$.} 
 In particular, {together with the Lipschitz continuity (Assumption~\ref{assum-lip}), for each term we use the mentioned inequality twice for different $\bar{\theta}$: for part (a), $\bar{\theta}=\lambda_1$ and $\bar{\theta}=\lambda_2$; for part (b), $\bar{\theta}=\gamma_1$ and $\bar{\theta}=\gamma_2$; and for part (c), we only use $\bar{\theta}=\gamma_2$ to derive}
 {\small
 \begin{align} \label{eq:bound-q-y}
 (a)&= M N\theta^k\sa{\fprod{\grad_{y_{j_k}}\Phi(\bx^k,\by^k)-\grad_{y_{j_k}}\Phi(\bx^k,\by^{k-1}),y_{j_k}^{k+1}-y_{j_k}^k}}
 \\ \nonumber
 &\quad+M N\theta^k\fprod{\grad_{y_{j_k}}\Phi(\bx^k,\by^{k-1})-\grad_{y_{j_k}}\Phi(\bx^{k-1},\by^{k-1}),y_{j_k}^{k+1}-y_{j_k}^k} \nonumber \\
 &\leq\frac{M N\theta^k}{2}\Big( \lambda_1L^2_{y_{j_k}y_{j_{k-1}}}\norm{\by^k-\by^{k-1}}_{\cY}^2+{\lambda_2}L^2_{y_{j_k}x_{i_{k-1}}}\norm{\bx^k-\bx^{k-1}}_{\cX}^2 +\norm{\by^{k+1}-\by^k}_{(\lambda_1^{-1}+\lambda_2^{-1})\id_n}^2\Big),\nonumber\\
 &=\frac{M N\theta^k}{2}\Big(\norm{\by^k-\by^{k-1}}_{ \lambda_1\bC_{\by}^2}^2+\norm{\bx^k-\bx^{k-1}}_{{\lambda_2}\bL_{\by\bx}^2}^2 +\norm{\by^{k+1}-\by^k}_{(\lambda_1^{-1}+\lambda_2^{-1})\id_n}^2\Big)+E_3,\nonumber\\
 \text{where }&{\small E_3\triangleq \frac{M N\theta^k}{2}\Big( \norm{\by^k-\by^{k-1}}_{\lambda_1(L^2_{y_{j_k}y_{j_{k-1}}}\bI_n-\bC^2_{\by})}^2+\norm{\bx^k-\bx^{k-1}}_{\lambda_2(L^2_{y_{j_k}x_{i_{k-1}}}\bI_m-\bL^2_{\by\bx})}^2\Big)},\nonumber\\
  \label{eq:bound-q-x}(b)&= M(N-1)\theta^k\sa{\fprod{\grad_{x_{i_k}}\Phi(\bx^k,\by^{k})-\grad_{x_{i_k}}\Phi(\bx^k,\by^{k-1}),x_{i_k}^{k}-x_{i_k}^{k+1}}} \\ \nonumber
 & \quad +M(N-1)\theta^k\fprod{\grad_{x_{i_k}}\Phi(\bx^k,\by^{k-1})-\grad_{x_{i_k}}\Phi(\bx^{k-1},\by^{k-1}),x_{i_k}^{k}-x_{i_k}^{k+1}} \nonumber \\
 &\leq\frac{M(N-1)\theta^k}{2}\Big(\gamma_2 L^2_{x_{i_k}y_{j_{k-1}}}\norm{\by^{k}-\by^{k-1}}_{\cY}^2+\gamma_1L^2_{x_{i_k}x_{i_{k-1}}}\norm{\bx^k-\bx^{k-1}}_{\cX}^2 +\norm{\bx^{k+1}-\bx^k}_{(\gamma_1^{-1}+\gamma_2^{-1})\id_m}^2\Big),\nonumber\\
 &=\frac{M(N-1)\theta^k}{2}\Big(\norm{\by^{k}-\by^{k-1}}_{\gamma_2 \bL_{\bx\by}^2}^2+\norm{\bx^k-\bx^{k-1}}_{\gamma_1\bC_\bx^2}^2 +\norm{\bx^{k+1}-\bx^k}_{(\gamma_1^{-1}+\gamma_2^{-1})\id_m}^2\Big)+E_4,\nonumber\\
 \text{where }&\erfan{E_4\triangleq \frac{M (N-1)\theta^k}{2}\big( \norm{\by^k-\by^{k-1}}_{\gamma_2(L^2_{x_{i_k}y_{j_{k-1}}}\bI_n-\bL^2_{\bx\by})}^2+\norm{\bx^k-\bx^{k-1}}_{\gamma_1(L^2_{x_{i_k}x_{i_{k-1}}}\bI_m-\bC_\bx^2)}^2\big),}\nonumber\\
 (c)&\leq \frac{N-1}{2}\Big(\gamma_2L^2_{x_{i_k}y_{j_k}}\norm{\by^{k+1}-\by^k}^2_{\cY}+\norm{\bx^{k+1}-\bx^k}^2_{\gamma_2^{-1}\id_m}\Big)\\
 &=\frac{N-1}{2}\Big(\norm{\by^{k+1}-\by^k}^2_{\gamma_2\bL^2_{\bx\by}}+\norm{\bx^{k+1}-\bx^k}^2_{\gamma_2^{-1}\id_m}\Big)+E_5,\nonumber
 \end{align}}%
 where $E_5\triangleq \frac{(N-1)}{2}\norm{\by^{k+1}-\by^k}_{\gamma_2(L^2_{x_{i_k}y_{j_k}}\bI_n-\bL^2_{\bx\by})}^2$.
 \ey{Recall $\{\bu_\bx^k,\bu_\by^k\}$ sequence defined in \eqref{eq:uk}. Now, using Lemma \ref{lem:inner-w} with $\cA=\frac{1}{M}\cA^k$ and $\delta^k=\bw_\bx^k$ together with Young's inequality and \eqref{eq:uk}, we can upper bound $(d)$. Similarly, by setting $\cA=\frac{1}{N}\tilde\cA^k$ and $\delta^k=-\bw_\by^k$ together with Young's inequality and \eqref{eq:uk}, we can upper bound $(e)$ using Lemma \ref{lem:inner-w} as follows:
 {\small
 \begin{align}
 (d)=&\fprod{\bw_\bx^k,\bx^k-\btx^{k+1}}+\fprod{\bw_\bx^k,\bx-\bu_\bx^k}+\fprod{\bw_\bx^k,\bu_\bx^k-\bx^k} \nonumber\\
 \leq& \frac{1}{2M}\norm{x_{i_k}^k-\tilde x_{i_k}^{k+1}}^2_{\cA^k} +M\norm{\bw_\bx^k}^2_{\nsa{*,(\cA^k)^{-1}}}+\frac{1}{M}\left(\bD_\cX^{\cA^k}(\bx,\bu_\bx^k)-\bD_\cX^{\cA^k}(\bx,\bu_\bx^{k+1})\right)+\fprod{\bw_\bx^k,\bu_\bx^k-\bx^k},\\
 (e)=&\fprod{\bw_\by^k,\bty^{k+1}-\by^k}+\fprod{\bw_\by^k,\bu_\by^k-\by}+\fprod{\bw_\by^k,\by^k-\bu_\by^k}\nonumber\\ 
 \leq & \frac{1}{2N}\norm{\by^k-\bty^{k+1}}^2_{\tilde\cA^k}+N\norm{\bw_\by^k}^2_{*,(\tilde\cA^k)^{-1}}+\frac{1}{N}\left(\bD_\cY^{\tilde\cA^k}(\by,\bu_\by^k)-\bD_\cY^{\tilde\cA^k}(\by,\bu_\by^{k+1})\right) +\fprod{\bw_\by^k,\by^k-\bu_\by^k},
 \end{align}}}%
 for any \sa{diagonal} $\cA^k\in\mathbb{S}_{++}^m$ and $\tilde\cA^k\in\mS_{++}^n$. 
 Hence, \sa{within \eqref{eq:pre_bound}, replacing the above} bounds on $(a)$, $(b)$,  $(c)$, $(d)$, and $(e)$ 
 we obtain the following bound:
 {\small
 \begin{align}\nonumber
 &\cL(\bx^{k+1},\by)-\cL(\bx,\by^{k+1})+(\theta^k-1)\big[(M-1)(\cL(\bx^k,\by)-\cL(\bx,\by))+(N-1)(\cL(\bx,\by)-\cL(\bx,\by^k))\big] \nonumber \\ 
 &  \leq (M-1)(\theta^kG_f^k(\bz)-G_f^{k+1}(\bz)) +(N-1)(\theta^k G_h^k(\bz) -G_h^{k+1}(\bz)) +M(G_\Phi^{k+1}(\bz)-\theta^kG_\Phi^k(\bz))\nonumber \\
 &  \mbox{ } + M\Big(\fprod{\grad_\by\Phi(\bx^k,\by^k),\by^{k+1}-\by}-\theta^k\fprod{\grad_\by\Phi(\bx^{k-1},\by^{k-1}),\by^k-\by}\Big)  \nonumber \\
 & \mbox{ } +(N-1)\Big(\fprod{\grad_\bx\Phi(\bx^k,\by^{k}),\bx-\bx^{k+1}}-\theta^k\fprod{\grad_\bx\Phi(\bx^{k-1},\by^{k-1}),\bx-\bx^k}\Big)\nonumber\\
 & \mbox{ }+\tfrac{1}{2}\norm{\btx^{k+1}-\bx^k}^2_{\bL_{\bx\bx}} +\tfrac{1}{N}A_1^{k+1}+\tfrac{1}{M}A_2^{k+1}+\tfrac{1}{2}\norm{\bty^{k+1}-\by^k}^2_{\bL_{\by\by}}+\tfrac{1}{2M}\norm{\bx^k-\btx^{k+1}}^2_{\cA^k} +M\norm{\bw_\bx^k}^2_{\nsa{*,(\cA^k)^{-1}}}\nonumber\\
 & \mbox{ }+\frac{1}{M}\big(\bD_\cX^{\cA^k}(\bx,\bu_\bx^k)-\bD_\cX^{\cA^k}(\bx,\bu_\bx^{k+1})\big)+\fprod{\bw_\bx^k,\bu_\bx^k-\bx^k}+\frac{1}{N}\left(\bD_\cY^{\tilde\cA^k}(\by,\bu_\by^k)-\bD_\cY^{\tilde\cA^k}(\by,\bu_\by^{k+1})\right)+\fprod{\bw_\by^k,\by^k-\bu_\by^k}\nonumber\\
 &\mbox{ }+\frac{M\theta^k}{2}\norm{\bx^k-\bx^{k-1}}^2_{N\lambda_2\bL_{\by\bx}^2+(N-1)\gamma_1\bC_{\bx}^2}+\frac{1}{2N}\norm{\by^k-\bty^{k+1}}^2_{\tilde\cA^k}+N\norm{\bw_\by^k}^2_{*,(\tilde\cA^k)^{-1}}\nonumber\\
 & \mbox{ }+\frac{M\theta^k}{2}\norm{\by^k-\by^{k-1}}_{(N-1)\gamma_2\bL^2_{\bx\by}+N \lambda_1\bC_{\by}^2}^2 \sa{+}\frac{1}{2}\norm{\bx^{k+1}-\bx^k}^2_{\substack{[M(N-1)\theta^k(\gamma_1^{-1}+\gamma_2^{-1})+(N-1)\gamma_2^{-1}]\id_m}}\nonumber\\
 &  \mbox{ }\sa{+}\frac{1}{2}\norm{\by^{k+1}-\by^k}^2_{\substack{M N\theta^k(\lambda_1^{-1}+\lambda_2^{-1})\id_n +(N-1){\gamma_2}\bL_{\bx\by}^2}}+\frac{1}{N}\fprod{\erfan{\tilde\bq}_\by^k,\bty^{k+1}-N\by^{k+1}+(N-1)\by^k}\nonumber\\ 
 &\mbox{ }+\frac{1}{M}\fprod{\erfan{\tilde\bq}_\bx^k,M\bx^{k+1}-(M-1)\bx^k-\btx^{k+1}}+E_1+E_2+\erfan{E_3+E_4+E_5}.\nonumber
 \end{align}}%
 Finally, the desired result follows from adding and subtracting $\bD_{\cY}^{\bS^k}(\by,\by^k) -\bD_{\cY}^{\bS^k}(\by,{\by}^{k+1})-\bD_{\cY}^{\bS^k-N\bL_{\by\by}-\tilde\cA^k}({\by}^{k+1},\by^k)$ and $\bD_{\cX}^{\bT^k}(\bx,\bx^k)$  $-\bD_{\cX}^{\bT^k}(\bx,{\bx}^{k+1})-{\bD_{\cX}^{\bT^k-M\bL_{\bx\bx}-\cA^k}({\bx}^{k+1},\bx^k)}$, and using strong convexity of the Bregman distances.
 \end{proof}
 
 \subsection{Proof of Lemma \ref{error_sampling}}\label{sec:proof-lem-error}
 We start by proving 
 \eqref{eq:part1}.
 {\small
 \begin{align*}
 {\mE}^k\Big[\norm{\bw_\bx^k}^2_{*,(\cA^k)^{-1}}\Big]&=\frac{1}{\alpha^k}\mE^k\Big[\norm{\bw_\bx^k}^2_*\Big]\\
 &\leq \frac{1}{\alpha^k}\mE^k\Big[2\norm{\mathbf{e}_\bx^k(\bx^k,\by^{k+1})}_*^2+2\erfan{(N-1)^2(\theta^k)^2}\norm{\mathbf{e}_\bx^k(\bx^k,\by^k)-\mathbf{e}_\bx^k(\bx^{k-1},\by^{k-1})}_*^2 \Big]\\
 &\leq \frac{2\sp{M}\delta_x^2}{v\alpha^k}+\erfan{\frac{8M(N-1)^2(\theta^k)^2}{\alpha^k}}\left(\norm{\bx^k-\bx^{k-1}}_{\bC_\bx^2}^2+\norm{\by^k-\by^{k-1}}_{\bL_{\bx\by}^2}^2\right),
 \end{align*}}%
 where in the first inequality we used the fact that $(a+b)^2\leq 2a^2+2b^2$ and the second inequality follows from Assumption \ref{assum:sample} and 
 Lipschitz continuity of $\grad_\bx\Phi$. 
\eqref{eq:part2} follows by a similar discussion as in part $(i)$.

 Next, recall that $\cF_{k-1}=\sigma(\Psi_{k-1})$ for $k\geq 1$ -- see Definition~\ref{def:sigma_algebra}. To show 
 \eqref{eq:part3}, note that $y_j^{k+1}$ is equal to $\tilde{y}_j^{k+1}$ with probability ${1\over N}$ and equal to $y_j^k$ with probability $(1-{1\over N})$ for all $j\in \cN$ (similar argument holds \nsa{also} for $x_i^{k+1}$ for all $i\in \cM$); \nsa{hence, it follows that} 
 for any \sp{separable function $\psi^y(\by)=\sum_{j\in\cN}\psi^y_j(y_j)$ such that $\psi^y_j:\cY_j\to \reals$ for $j\in\cN$, we have $\mE[\psi^y(\by^{k+1}) \mid \sigma(\psi^y_{k-1}\cup\cB^k)]=\frac{1}{N}\psi^y(\tilde\by^{k+1})+\frac{(N-1)}{N}\psi^y(\by^k)$} which implies that $\sp{\mE^k[\psi^y(\tilde\by^{k+1})]}=N\mE^k[\psi^y(\by^{k+1})]-(N-1)\psi^y(\by^k)$. 
 Similarly, for 
 $\psi^x(\bx)=\sum_{i\in\cM}\psi^x_i(x_i)$ such that $\psi^x_i:\cX_i\to\reals$,  $\sp{\mE[\psi^x(\bx^{k+1})\mid \sigma(\Psi_{k-1}\cup\{j_k,\cB^k,\cS^k\})]}=\frac{1}{M}\psi^x(\tilde\bx^{k+1})+(1-\frac{1}{M})\psi^x(\bx^k)$. Thus, 
 {\small
 \begin{align}
    \mE^k[\psi^x(\tilde\bx^{k+1})]=M\mE^k[\psi^x(\bx^{k+1})]-(M-1)\psi^x(\bx^k),\quad \mE^k[\psi^y(\tilde\by^{k+1})]=N\mE^k[\psi^y(\by^{k+1})]-(N-1)\psi^y(\by^k).\label{eq:Ek-x-relation} 
 \end{align}}%
 Therefore, \sp{using some particular separable functions}, we conclude 
 that for any $\by\in\cY$ and $\bx\in\cX$, we get:
 \begin{subequations}
 \label{eq:tilde-expectation}
 {\small
 \begin{align}        &\sp{\mE^k[\bD^{{\bS^k}}_\cY(\by,\bty^{k+1})]}=N\mE^k[\bD^{\sa{\bS^k}}_\cY(\by,\by^{k+1})]-(N-1)\bD^{\sa{\bS^k}}_\cY(\by,\by^k),\label{eq:EDy}\\ &\sa{\mE^k[}\bD^{\sa{\bT^k}}_\cX(\bx,\btx^{k+1})]=M\mE^k[\bD^{\sa{\bT^k}}_\cX(\bx,\bx^{k+1})]-(M-1)\bD^{\sa{\bT^k}}_\cX(\bx,\bx^k).\label{eq:EDx}\\
        &\mE^k[f(\btx^{k+1})]=M\mE^k[f(\bx^{k+1})]-(M-1)f(\bx^k),\label{eq:Ek-x}\\
        &\sp{\mE^k[h(\bty^{k+1})]}=N\mE^k[h(\by^{k+1})]-(N-1)h(\by^k).\label{eq:Ek-y}\\
        &N\mE^k[\Phi(\bx^k,\by^{k+1})]=\sum_{j\in\cN}\sp{\mE^k}[\Phi(\bx^k,\by^k+V_j(\tilde{y}_j^{k+1}-y_j^k))], \label{eq:phi-y}\\
     &M\mE^k[\Phi(\bx^{k+1},\by^{k+1})]
 =\sum_{i\in\cM}\mE^k[\Phi(\bx^k+U_i(\tilde{x}_i^{k+1}-x_i^k),\by^{k+1})]. \label{eq:phi-x-y} 
 \end{align}}%
 \end{subequations}
 \erfan{Moreover, one can easily verify that 
 $\sp{\mE^k[L^2_{x_{i_k}x_{i_{k-1}}}]}=\frac{1}{M}\sum_{i=1}^M L^2_{x_{i}x_{i_{k-1}}}=C^2_{x_{i_{k-1}}}$; hence, $
 \sp{\mE^k[L^2_{x_{i_k}x_{i_{k-1}}}\norm{\bx^k-\bx^{k-1}}_\cX^2]}=\norm{\bx^k-\bx^{k-1}}^2_{\bC_\bx^2}$, \rev{where $\{C_{x_i}\}_{i\in\cM}$ are defined after Remark \ref{rem:L-constants} and $C_{\bx}$ is defined in \eqref{def L}.} Similarly, 
 $\sp{\mE^k[L^2_{y_{j_k}x_{i_{k-1}}}}\norm{\bx^k-\bx^{k-1}}_\cX^2]=\norm{\bx^k-\bx^{k-1}}^2_{\bL_{\by\bx}^2}$. Therefore, the definition of 
$\bU^k$ in Lemma \ref{lem:one-step} implies that $\mE^k[\norm{\bx^k-\bx^{k-1}}^2_{\bU^k}]=0$. \sp{Furthermore, $\mE^k[\norm{\by^{k+1}-\by^{k}}^2_{\bV_1^k}]=(N-1)\gamma_2\mE^k[(L^2_{x_{i_k}y_{j_k}}-L_{\bx y_{j_k}}^2)\norm{y^{k+1}_{j_k}-y^k_{j_k}}_{\cY_{j_k}}^2]$; hence, $\mE[\norm{\by^{k+1}-\by^{k}}^2_{\bV_1^k}\mid~\sigma(\Psi_{k-1}\cup\{j_k,\cB^k\})]=\mathbf{0}$, implying that $\mE^k[\norm{\by^{k+1}-\by^{k}}^2_{\bV_1^k}]=0$, where we used $\mE^k[L^2_{x_{i_k}y_{j_k}}-L_{\bx y_{j_k}}^2]=0$.} Using a similar argument we also conclude that 
$\mE^k[\norm{\by^k-\by^{k-1}}^2_{\bV_2^k}]=0$ holds.}
 {Now, putting these results together and using the relations in \eqref{eq:EDy}-\eqref{eq:phi-x-y}, we conclude that \nsa{$\mE^k[\Xi_\bx^k(\bx)]=\mE^k[\Xi_\by^k(\by)]=0$,} for any $k\geq 0$, \nsa{$\bx\in\dom f$ and $\by\in\dom h$}.
 
 Finally, {for $\mathbf{e}_\bx^k(\bx,\by)$ and \sp{$\mathbf{e}_\by^k(\bx,\by)$} 
 in Definition~\ref{def:ek}, \sp{using Assumption~\ref{assum:sample} we 
 have $\mE^k[\mathbf{e}^k_\by(\bx^k,\by^k)]=\mE^k[\mathbf{e}^k_\by(\bx^{k-1},\by^{k-1})]=\mathbf{0}$ and $\mE^k[\mathbf{e}^k_\bx(\bx^k,\by^k)]=\mE^k[\mathbf{e}^k_\bx(\bx^{k-1},\by^{k-1})]=\mathbf{0}$. Moreover,} $\mE^k[\sp{\mathbf{e}_\bx^k}(\bx^k,\by^{k+1})]=\mE^k[\mE[\mathbf{e}^k(\bx^k,\by^{k+1})\mid \sigma\big(\Psi_{k-1}\cup\{j_k,~i_k,~\sp{\cB^k}\}\big)]]=\mathbf{0}$, 
 for $k\geq 0$;} hence, 
 we get $\mE^k[\sp{\bw_\bx^k}]=\mE^k[\sp{\bw_\by^k}]=\mathbf{0}$ and 
 \sp{$\mE^k[\fprod{\bw_\bx^k,\bu_\bx^k-\bx^k}]=0$ and $\mE^k[\fprod{\bw_\by^k,\by^k-\bu_\by^k}]=0$ since both $\bu_\bx^k$ and $\bu_\by^k$ are $\cF_{k-1}$-measurable.}}
 \end{proof}
\section{Proof of the Meta-Result in Theorem~\ref{thm:meta}}

\subsection{Proof of Lemma \ref{lem:result-summable}.}\label{sec:proof-lem-summable}
Fix an arbitrary $\varsigma>0$. 
 Let \sp{$\bz^*=(\bx^*,\by^*)\in Z^*$ be an arbitrary} saddle point of $\cL$ defined in~\eqref{eq:original-problem}. 
 \erfan{Using the definitions of $\tilde\bq_\bx^k= M(N-1)\theta^k(\grad_\bx\Phi(\bx^k,\by^{k})-\grad_\bx\Phi(\bx^{k-1},\by^{k-1}))$, $\tilde\bq_\by^k= MN\theta^k(\grad_\by\Phi(\bx^k,\by^k)-\grad_\by\Phi(\bx^{k-1},\by^{k-1}))$, the following inner products can be bounded similar to \eqref{eq:bound-q-y} and \eqref{eq:bound-q-x} for any $k\geq 0$,}
\begin{subequations}
\label{eq:qk-bound}
{\small
 \begin{flalign}
     &\tfrac{1}{N}\Big|\fprod{\erfan{\tilde\bq}_\by^k,\by-\by^k}\Big|\leq M\sa{\theta^k}\Big(\bD_{\cY}^{(\lambda_1^{-1}+\lambda_2^{-1})\bI_n}(\by,\by^k)+\bD_{\cY}^{\lambda_1N\bC_{\by}^2}(\by^k,\by^{k-1})+\bD_{\cX}^{\lambda_2N\bL_{\by\bx}^2}(\bx^k,\bx^{k-1})\Big),
     \\  
     &\tfrac{1}{M}\Big|\fprod{\erfan{\tilde\bq}_\bx^k,\bx^k-\bx}\Big|\leq (N-1)\theta^k\Big(\bD_{\cX}^{(\gamma_1^{-1}+\gamma_2^{-1})\bI_m}(\bx,\bx^k)+\bD_{\cX}^{\gamma_1M\bC_{\bx}^2}(\bx^k,\bx^{k-1})+\bD_{\cY}^{\gamma_2M\bL_{\bx\by}^2}(\by^k,\by^{k-1})\Big). 
 \end{flalign}}%
 \end{subequations}
\rev{Recall that $G_f^{k+1}(\bz^*)=\cL(\bx^{k+1},\by^*)-\cL(\bx^*,\by^*)$, then using \eqref{eq:G-connection} in Lemma \ref{lem:pos-theta}, we obtain
{\small
 \begin{align}
 &t^k\Big[{\cL(\bx^{k+1},\by^*)-\cL(\bx^*,\by^{k+1})}+R^{k+1}(\bz^*)+(M-1)G_f^{k+1}(\bz^*)+(N-1)G_h^{k+1}(\bz^*)- M G_\Phi^{k+1}(\bz^*)\Big]\nonumber \\
 &\geq t^k[\cL(\bx^{k+1},\by^*)-\cL(\bx^*,\by^{k+1})+R^{k+1}(\bz^*)]+(M-1)t^k(\cL(\bx^{k+1},\by^*)-\cL(\bx^*,\by^*)) \nonumber\\
 &\quad + (N-1)t^k(\cL(\bx^*,\by^*)-\cL(\bx^*,\by^{k+1})+\fprod{\nabla_\bx\Phi(\bx^k,\by^k),\bx-\bx^k})-Mt^k\fprod{\nabla_\by\Phi(\bx^k,\by^k),\by^*-\by^k}\nonumber\\ 
 &\geq t^k[\cL(\bx^{k+1},\by^*)-\cL(\bx^*,\by^{k+1})]+(M-1)t^k(\cL(\bx^{k+1},\by^*)-\cL(\bx^*,\by^*))\nonumber\\
 &\quad +(N-1)t^k(\cL(\bx^*,\by^*)-\cL(\bx^*,\by^{k+1}))+t^k\Big(\bD_{\cX}^{\bT^{k}}(\bx^*,\bx^{k+1})+\bD_{\cY}^{\bS^{k}}(\by^*,\by^{k+1})\nonumber \\
 &\quad {+ \bD_{\cX}^{\bM_3^k-\cA^k}(\bx^{k+1},\bx^k)+\bD_{\cY}^{\bM_4^k-\tilde\cA^k}(\by^{k+1},\by^k)}\Big)+\frac{t^{k+1}}{M}\fprod{{\tilde\bq_\bx^{k+1}},~\bx^*-\bx^{k+1}}+\frac{t^{k+1}}{N}\fprod{{\tilde\bq_\by^{k+1}},~\by^{k+1}-\by^*}\nonumber\\
 &\geq {\bD^{t^{k}\bM_3^{k}}_\cX(\bx^*,\bx^{k+1})+\bD^{t^{k}\bM_4^{k}}_\cY(\by^*,\by^{k+1})}\geq {\varsigma\bD^{\bT^0}_\cX(\bx^*,\bx^{k+1})+\varsigma\bD^{\bS^0}_\cY(\by^*,\by^{k+1}),} \label{Qk-b}
 \end{align}}%
 where 
 in the second inequality we used the definition of $R^{k+1}(\cdot)$ in \eqref{eq:Rk-def}, along with the relation $t^k = t^{k+1} \theta^{k+1}$ (rom Assumption \ref{assum:step}) and definitions of $\tilde\bq_\bx^k$ and $\tilde\bq_\by^k$, moreover, in the third inequality} 
 we used \sp{\eqref{eq:qk-bound},}  
 $\cL({\bx}^{k},\by^*)-\cL(\bx^*,\by^*)\geq 0$, $\cL({\bx}^*,\by^*)-\cL(\bx^*,\by^{k})\geq 0$,  $t^k(\bM_3^k-\cA^k)\succeq t^{k+1}\bM_1^{k+1}$ and $t^k(\bM_4^k-\tilde\cA^k)\succeq t^{k+1}\bM_2^{k+1}$, which follow from \eqref{eq:step-size-lip} since $\cD^{k+1}\succeq 0$ and $\tilde\cD^{k+1}\succeq 0$ for both $\bar\delta>0$ and $\bar\delta=0$,
 \sp{and we also used $\bT^k-(N-1)(\gamma_1^{-1}+\gamma_2^{-1})\id\succeq\bM_3^k$ and $\bS^k-M(\lambda_1^{-1}+\lambda_2^{-1})\id\succeq\bM_4^k$, which follows from the definitions of $\bM_3^k$ and $\bM_4^k$,} and for the last inequality we used \eqref{eq:step-size-lip} \sp{together with $\bM_1^{k+1}\succeq 0$ and $\bM_2^{k+1}\succeq 0$}.

 \rev{
 From \eqref{eq:Rk-def} and \eqref{eq:Qk-def}, we can provide a lower bound for $t^k R^{k+1}(\bz)$ in terms of $t^{k+1}Q^{k+1}(\bz)$ by using the relations between the coefficients of each term based on \eqref{eq:step-size-condition} in Assumption~\ref{assum:step} 
 as follows,}
 {\small
\begin{align}
 &\sp{t^k R^{k+1}(\bz) \geq t^{k+1} Q^{k+1}(\bz)+\bD^{t^{k+1}\cD^{k+1}+\varsigma\bT^0}_\cX(\bx^{k+1},\bx^k)+\bD_\cY^{t^{k+1}\tilde\cD^{k+1}+\varsigma\bS^0}(\by^{k+1},\by^k)}, \quad\forall~\bz\in\cZ,\label{Qk-a} 
 \end{align}}%
\sp{with $\cD^{k+1}=\mathbf{0}$ and $\tilde\cD^{k+1}=\mathbf{0}$ when $\bar\delta=0$.} 
\rev{Next, rearranging the terms in \eqref{eq:one-step-main} evaluated at $(\bx,\by)=(\bx^*,\by^*)$, multiplying both sides by $t^k$, using $t^k=t^{k+1}\theta^{k+1}$ as in Assumption \ref{assum:step}, taking conditional expectation, and using Lemma \ref{error_sampling} we obtain
{\small
 \begin{align*}
 &{t^{k}}\mE^k\Big[{\cL(\bx^{k+1},\by^*)-\cL(\bx^*,\by^{k+1})}+R^{k+1}(\bz^*)+(M-1)G_f^{k+1}(\bz^*)+(N-1)G_h^{k+1}(\bz^*)- M G_\Phi^{k+1}(\bz^*)\Big]  \\
 &\leq t^kQ^k(\bz^*)+t^{k-1}\big[(M-1)G_f^k(\bz^*)+(N-1)G_h^k(\bz^*)- M G_\Phi^k(\bz^*)\big]+B^k_\delta+\bD_\cX^{t^k\cD^k}(\bx^k,\bx^{k-1})+\bD_\cY^{t^k\tilde\cD^k}(\by^k,\by^{k-1}), \nonumber
 \end{align*}}%
 where ${B^k_\delta}\triangleq 0$ 
 {when $\bar\delta=0$;} ${B^k_\delta}\triangleq 2{t^k}\left(\frac{{M^2}\delta_x^2}{v\alpha^k}+\frac{{N^2}\delta_y^2}{u\tilde\alpha^k}\right)$ 
 {when $\bar\delta>0$} and we also dropped the non-negative terms involving $\cL({\bx}^{k},\by^*)-\cL(\bx^*,\by^*)\geq 0$ and $\cL({\bx}^*,\by^*)-\cL(\bx^*,\by^{k})\geq 0$ in the left-hand side.}  
 \rev{Next, using \eqref{Qk-a} as well as adding and subtracting $t^{k-1}(\cL(\bx^{k},\by^*)-\cL(\bx^*,\by^{k}))$ to the right hand side, we obtain} 
 {\small
 \begin{align}\label{eq:bound-iterates}
 &\sa{t^{k}}\mE^k\Big[\eyz{\cL(\bx^{k+1},\by^*)-\cL(\bx^*,\by^{k+1})}+R^{k+1}(\bz^*)+(M-1)G_f^{k+1}(\bz^*)+(N-1)G_h^{k+1}(\bz^*)- M G_\Phi^{k+1}(\bz^*)\Big] \nonumber \\
 &\leq   \sa{t^{k-1}}\big[\eyz{\cL(\bx^{k},\by^*)-\cL(\bx^*,\by^{k})}+R^k(\bz^*)+(M-1)G_f^k(\bz^*)+(N-1)G_h^k(\bz^*)- M G_\Phi^k(\bz^*)\big]\nonumber\\
 &\quad +\sp{B^k_\delta-\varsigma \bD^{\bT^0}_\cX(\bx^k,\bx^{k-1})-\varsigma\bD_\cY^{\bS^0}(\by^k,\by^{k-1})}-\eyz{t^{k-1}\big(\cL(\bx^{k},\by^*)-\cL(\bx^*,\by^{k})\big)}.
 \end{align}}%
 \sp{Let $b_k\triangleq t^{k-1}(\cL(\bx^{k},\by^*)-\cL(\bx^*,\by^{k}))+\varsigma(\bD_\cX^{\bT^0}(\bx^{k},\bx^{k-1})+\bD_\cY^{\bS^0}(\by^{k},\by^{k-1})) \us{ \ \geq \ 0}$ and $c_k\triangleq \sp{B^k_\delta} \us{ \ \geq \ 0}$ for $k\geq 0$, and $a_k\triangleq {t^{k-1}}\Big(\cL(\bx^{k},\by^*)-\cL(\bx^*,\by^{k})+R^k(\bz^*)+(M-1)G_f^k(\bz^*)+(N-1)G_h^k(\bz^*)- M G_\Phi^k(\bz^*)\Big)$ for $k\geq 1$ and
 {\small
 \begin{align*}
     a_0
     &\triangleq \cL(\bx^{0},\by^*)-\cL(\bx^*,\by^{0})+Q^0(\bz^*)+(M-1)G_f^0(\bz^*)+(N-1)G_h^0(\bz^*)- M G_\Phi^0(\bz^*)\\
     &=\bD_\cX^{\bT^0+\frac{1}{M}\cA^0}(\bx^*,\bx^0)+M(\cL(\bx^{0},\by^*)-\cL(\bx^*,\by^*))+(N-1)(\Phi(\bx^*,\by^0)-\Phi(\bx^0,\by^0)-\fprod{\grad_\bx\Phi(\bx^0,\by^0),\bx^*-\bx^0})\\
     &\quad +\bD_\cY^{\bS^0+\frac{1}{N}\tilde\cA^0}(\by^*,\by^0)+N(\cL(\bx^*,\by^*)-\cL(\bx^*,\by^{0}))+M(\Phi(\bx^0,\by^0)+\fprod{\grad_\by\Phi(\bx^0,\by^0),\by^*-\by^0}-\Phi(\bx^0,\by^0)),
 \end{align*}}%
 where we used $\bu_\bx^0=\bx^0=\bx^{-1}$, $\bu_\by^0=\by^0=\by^{-1}$ and the definitions of $Q^0(\cdot)$, $G_f^0(\cdot)$, $G_h^0(\cdot)$ and $G^0_\Phi(\cdot)$. Hence,}
 {\small
 \begin{equation}\label{eq:a0-bound}
 \begin{aligned}
     0\leq a_0 
     &\leq \bD_\cX^{\bT^0+\frac{1}{M}\cA^0+(N-1)M\bL_{\bx\bx}}(\bx^*,\bx^0)+\bD_\cY^{\bS^0+\frac{1}{N}\tilde\cA^0+MN\bL_{\by\by}}(\by^*,\by^0)\\
     &\quad+M(\cL(\bx^{0},\by^*)-\cL(\bx^*,\by^*))+N(\cL(\bx^*,\by^*)-\cL(\bx^*,\by^{0})),
 \end{aligned}
  \end{equation}}%
 where the lower bound follows from $\Phi$ being convex-concave, and the upper bound follows from 
 \sp{\cite[(2.10) in Lemma~2 with $\alpha=0$]{nesterov2012efficiency}. Thus, we have 
 \begin{align}
    \label{eq:a0-bound-simple}
     a_0\leq B(\bz^*)\triangleq\Delta(\bz^*)+ C_0(\bz^*)+\cL(\bx^0,\by^*)-\cL(\bx^*,\by^0),
 \end{align} 
 where $\Delta(\cdot)$ and $C_0(\cdot)$ are defined in~\eqref{eq:delta-function} and \eqref{eq:C0-function}, respectively.} \sp{Moreover, according to \eqref{Qk-b} we have $a_k\geq 0$ for $k\geq 1$ as well, and $\sum_{k=0}^\infty c_k<\infty$ since ${\sum_{k=0}^\infty t^k(\frac{1}{\alpha^k}+\frac{1}{\tilde\alpha^k})}<+\infty$ from Assumption~\ref{assum:step}; therefore,} from Lemma~\ref{lem:supermartingale}, $\lim_{k\rightarrow \infty} a_k$ exists 
 and $\sum_{k=1}^\infty b_k <\infty$ almost surely. 
 Therefore, 
 $\sum_{k=0}^\infty \bD_\cX(\bx^{k},\bx^{k-1})+\bD_\cY(\by^{k},\by^{k-1})\leq \sum_{k=1}^\infty b_k <\infty$ almost surely. \sp{Note that 
 taking expectation of both sides of \eqref{eq:bound-iterates}, we get $\mE[a_{k+1}]\leq\mE[a_k]-\mE[b_k]+c_k$ for $k\geq 0$; hence, 
 $0\leq \varsigma\sum_{k=0}^{K} \mE[\bD_\cX^{\bT^0}(\bx^{k+1},\bx^{k})+\bD_\cY^{\bS^0}(\by^{k+1},\by^{k})]\leq \sum_{k=0}^{K}\mE[b_k]\leq \mE[a_0]-\mE[a_{K+1}]+\sum_{k=0}^Kc_k\leq \mE[a_0]+\sum_{k=0}^\infty c_k$ for all $K\geq 0$.} Thus, as $K\to \infty$,  
 using the monotone convergence theorem, we can interchange the sum and expectation operations to obtain the desired result.

\rev{Next, we show that for any $p\in (0,1)$ there exists $\Delta_p>0$ with $\Delta_p=\mathcal{O}(1/p)$ such that $\mathbb P\!\left(\sup_{k\ge0}\{\| \bx^k\|^2_\cX+\| \by^k\|^2_\cY\}\le \Delta_p\right)\ge 1-p$. Using \eqref{eq:bound-iterates} and the definitions of $a_k,b_k,c_k\geq 0$ given above, it holds that $\{a_k\}$ and $\{b_k\}$ are adapted sequences, $\{c_k\}$ is deterministic, and
\begin{equation}\label{eq:drift}
\mathbb E^k[a_{k+1}]\;\le\; a_k-b_k+c_k\quad\text{a.s.}\quad\forall~k\geq 0.
\end{equation}
Define $\bar a_k\triangleq a_k+\sum_{i=k}^{\infty}c_i$ for $k\geq 0$. Note that $\{\bar a_k\}$ is a nonnegative, adapted sequence; moreover, $a_0\geq 0$ is a deterministic and $0\leq \sum_{k=0}^\infty c_k<\infty$ --hence, $0\leq \bar{a}_0<+\infty$. Therefore, we can conclude that
\begin{equation}\label{eq:abar-seq}
\mathbb E^k[\bar a_{k+1}]
= \mathbb E^k[a_{k+1}]+\sum_{i=k+1}^{\infty}c_i
\overset{\eqref{eq:drift}}{\le} a_k-b_k+c_k+\sum_{i=k+1}^{\infty}c_i
= \bar a_k-b_k\leq \bar a_k,\quad\forall~k\geq 0.
\end{equation}
Hence $\{\bar a_k\}_{k\geq 0}$ is a nonnegative super-martingale. 
Now, fix $\gamma>0$ and $\bar k\in\mathbb Z_+$. Define the bounded stopping time 
$T_\gamma\triangleq \inf\{k\ge0:\ \bar a_k\ge\gamma\}\wedge \bar k$, where $a\wedge b\triangleq \min\{a,b\}$. 
Since $\{T_\gamma = k\}\in\mathcal F_k$, $T_\gamma$ is a stopping time for the super-martingale sequence $\{\bar a_k\}_{k\geq 0}$; therefore, we can conclude that the stopped process $\{\bar{a}_{k\wedge T_\gamma}\}_{k\geq 0}$ is also a nonnegative super-martingale. Thus, optional stopping theorem implies that 
\begin{equation}\label{eq:OST}
\mathbb E[\bar a_{k\wedge T_\gamma}] \leq \mathbb E[\bar a_0] = \bar a_0,\quad \forall~k\geq 0.
\end{equation}
Note that $\bar a_{T_\gamma}=\bar a_{\bar k\wedge T_\gamma}$; hence, from \eqref{eq:OST} it follows that  
\begin{align}
\label{eq:OST-2}
    \mathbb E[\bar a_{T_\gamma}]=\mathbb E[\bar a_{\bar k\wedge T_\gamma}] \leq \bar a_0.
\end{align}
Next, since $\bar a_{T_\gamma}\geq 0$, by Markov’s inequality and \eqref{eq:OST-2}, we have 
\[
\mathbb P\!\left(\sup_{0\le k\le \bar k}\bar a_k\ge\gamma\right) = 1- \mathbb P\!\left(\sup_{0\le k\le \bar k}\bar a_k<\gamma\right)=1 - \mathbb P(\bar a_{T_\gamma}<\gamma)
= \mathbb P(\bar a_{T_\gamma}\ge\gamma)\ \le\ \frac{\mathbb E[\bar a_{T_\gamma}]}{\gamma}\ \le\ \frac{\bar a_0}{\gamma}.
\]
Let $A(\bar k)\triangleq \left\{\sup_{0\le k\le \bar k}\bar a_k\ge\gamma\right\}$ for $\bar k\in\integers_+$. Note that $A(k')\subset A(k'')$ for $k',k''\in\integers_+$ such that $k''\geq k'$; hence, $\lim_{\bar k\to\infty}A(\bar k)=\bigcup_{k\geq 0}A(k)$.
Moreover, using continuity from below property of probability, we obtain $\lim_{\bar k\to \infty}\mathbb P\!\left(A(\bar k)\right)=\mathbb P\!\left(\lim_{\bar k\to \infty} A(\bar k)\right)$. Therefore, we can conclude that $\mathbb P\!\left(\sup_{k\ge0}\bar a_k\ge\gamma\right)\ \le\ \lim_{\bar k\to\infty}\mathbb P\!\left(\sup_{0\le k\le \bar k}\bar a_k\ge\gamma\right)\ \le\ \frac{\bar a_0}{\gamma}$ for any $\gamma>0$. 
Hence, selecting $\gamma=\bar a_0/p$ for any $p\in(0,1)$ leads to
$\mathbb P\!\left(\sup_{k\ge0}\bar a_k\le \frac{\bar a_0}{p}\right)\ \ge\ 1-p$. 
Since $a_k\le\bar a_k$ for all $k\geq 0$ by definition of $\bar a_k$, we also get
\begin{equation}\label{eq:ak-high-prob-bound}
\mathbb P\!\left(\sup_{k\ge0}a_k\le \frac{\bar a_0}{p}\right)\ \ge\
\mathbb P\!\left(\sup_{k\ge0}\bar a_k\le \frac{\bar a_0}{p}\right)\ \ge\ 1-p .
\end{equation}
Finally, note that \eqref{Qk-b} implies $a_k\geq \varsigma\bD^{\bT^0}_\cX(\bx^*,\bx^{k})+\varsigma\bD^{\bS^0}_\cY(\by^*,\by^{k})$ for $k\geq 0$; moreover, recall that $\sum_{k=0}^\infty c_k=\sum_{k=0}^\infty B_\delta^k<\infty$ and from \eqref{eq:a0-bound-simple} we also have $a_0\leq B(\bz^*)$. Therefore, we can conclude that 
\begin{equation*}
\mathbb P\!\left(\sup_{k\ge0}\left\{\bD^{\bT^0}_\cX(\bx^*,\bx^{k+1})+\bD^{\bS^0}_\cY(\by^*,\by^{k+1})\right\}\leq \frac{B(\bz^*)+\sum_{k=0}^{\infty}B_\delta^k}{\varsigma p}\right) \geq\ 1-p ,
\end{equation*}
and the desired result follows from 1-strong convexity of Bregman distance functions.} 
 \end{proof}
\subsection{Proof of Lemma~\ref{lem:sup-bound}}
\label{sec:proof-sup-bound}
 \rev{From the definitions of $\Xi_x^k$ and $\Xi_y^k$ in \eqref{eq:error-x} and \eqref{eq:error-y}, and using \eqref{eq:Ek-x-relation} and
 \eqref{eq:tilde-expectation} 
 together with $\mE^k[\norm{\bx^k-\bx^{k-1}}^2_{\bU^k}]=\mE^k[\norm{\by^{k+1}-\by^{k}}^2_{\bV_1^k}]=\mE^k[\norm{\by^k-\by^{k-1}}^2_{\bV_2^k}]=0$, which follows from the proof of Lemma \ref{error_sampling}, we get
 {\small
 \begin{align}\label{eq:bound-Xi}
    &\mE\Big[\sup_{\bz\in Z} \Big\{C_0(\bz)+ \bD_\cX^{\bT^0+(N-1)M\bL_{\bx\bx}+\frac{1}{M}\cA^0}(\bx,\bx^0)+\bD_\cY^{\bS^0+NM\bL_{\by\by}+\frac{1}{N}\tilde\cA^0}(\by,\by^0)+\sum_{k=0}^{K-1}t^k\Xi_x^k(\bx)+\sum_{k=0}^{K-1}t^k\Xi_y^k(\by)\Big\}\Big]\nonumber \\
    &\leq \mE\Big[\sup_{\bz\in Z} \Big\{C_0(\bz)+ \bD_\cX^{\bT^0+(N-1)M\bL_{\bx\bx}+\frac{1}{M}\cA^0}(\bx,\bx^0)+\bD_\cY^{\bS^0+NM\bL_{\by\by}+\frac{1}{N}\tilde\cA^0}(\by,\by^0) \nonumber\\
    &\quad + \sum_{k=0}^{K-1}t^k\left(\tfrac{1}{M}\big(\bD^{\bT^k}_\cX(\bx,\bx^k)-\bD^{\bT^k}_\cX(\bx,\tilde\bx^{k+1})\big) -\big(\bD^{\bT^k}_\cX(\bx,\bx^k)-\bD^{\bT^k}_\cX(\bx,\bx^{k+1})\big)\right)\nonumber\\
    &\quad +\sum_{k=0}^{K-1}t^k\left(\tfrac{1}{N}\big(\bD^{\bS^k}_\cY(\by,\by^k)-\bD^{\bS^k}_\cY(\by,\tilde\by^{k+1})\big) -\big(\bD^{\bS^k}_\cY(\by,\by^k)-\bD^{\bS^k}_\cY(\by,\by^{k+1})\big)\right) \Big\}\Big].
 \end{align}}%
 Next, we provide some upper bounds for the terms in the last two lines of the r.h.s of the above inequality.}
 For $k\geq 0$, let $\tilde\Gamma_\bx^{k+1}\triangleq \bT^k(\grad\varphi_\cX(\tilde\bx^{k+1})-\grad\varphi_\cX(\bx^k))$ and $\Gamma_\bx^{k+1}\triangleq\bT^k(\grad\varphi_\cX(\bx^{k+1})-\grad\varphi_\cX(\bx^k))$. Then from the definition of Bregman function, one can observe that $\bD^{\bT^k}_\cX(\bx,\bx^k)-\bD^{\bT^k}_\cX(\bx,\tilde\bx^{k+1})-\bD^{\bT^k}_\cX(\tilde\bx^{k+1},\bx^k)=\fprod{\nsa{\tilde\Gamma_\bx^{k+1}},~\bx-\tilde\bx^{k+1}}$ and $\bD^{\bT^k}_\cX(\bx,\bx^k)-\bD^{\bT^k}_\cX(\bx,\bx^{k+1})-\bD^{\bT^k}_\cX(\bx^{k+1},\bx^k)=\fprod{\nsa{\Gamma_\bx^{k+1}},~\bx-\tilde\bx^{k+1}}$. Thus,
 {\small
 \begin{equation}\label{eq:bregman-inner}
 \begin{aligned}
 \MoveEqLeft\tfrac{1}{M}\Big(\bD^{\bT^k}_\cX(\bx,\bx^k)-\bD^{\bT^k}_\cX(\bx,\tilde\bx^{k+1})\Big) -\Big(\bD^{\bT^k}_\cX(\bx,\bx^k)-\bD^{\bT^k}_\cX(\bx,\bx^{k+1})\Big) \\
 &= \tfrac{1}{M}\Big(\bD^{\bT^k}_\cX(\tilde\bx^{k+1},\bx^k)+\fprod{\tilde\Gamma_\bx^{k+1},\bx-\tilde\bx^{k+1}}\Big) -\Big(\bD^{\bT^k}_\cX(\bx^{k+1},\bx^k)+\fprod{\Gamma_\bx^{k+1},\bx-\bx^{k+1}}\Big),\\ 
 &=\erfan{-\tfrac{1}{M}\bD_\cX^{\bT^k}(\bx^k,\tilde\bx^{k+1})+\bD_\cX^{\bT^k}(\bx^k,\bx^{k+1})+\fprod{\tfrac{1}{M}\tilde\Gamma_\bx^{k+1}-\Gamma_\bx^{k+1},\bx-\bx^k}}, 
 \end{aligned}
 \end{equation}}%
 where we used 
 $\fprod{\tilde\Gamma_\bx^{k+1},\bx^k-\tilde\bx^{k+1}}=-\bD_\cX^{\bT^k}(\bx^k,\tilde\bx^{k+1})-\bD_\cX^{\bT^k}(\tilde\bx^{k+1},\bx^k)$ \nsa{and $\fprod{\Gamma_\bx^{k+1},\bx^k-\bx^{k+1}}=-\bD_\cX^{\bT^k}(\bx^k,\bx^{k+1})-\bD_\cX^{\bT^k}(\bx^{k+1},\bx^k)$.}
 \rev{Next, we provide an upper bound for the inner product
 on the r.h.s. of 
 \eqref{eq:bregman-inner} by invoking Lemma \ref{lem:inner-w}.} Specifically, for the auxiliary sequence $\{\bv_\bx^k\}$ defined
 such that \sp{$\bv_\bx^0=\bx^0\in\cX$, $\bv_\bx^{k+1}\triangleq\argmin_{\bx\in\cX}\{-\fprod{\bdelta_\bx^k,\bx}+\bD_\cX^{{\bT^k}}(\bx,\bv_\bx^k)\}$ and \sp{$\bdelta_\bx^k\triangleq \tfrac{1}{M}\tilde{\Gamma}_\bx^{k+1}-\Gamma_\bx^{k+1}$} for $k\geq 0$, Lemma \ref{lem:inner-w} implies that}
 \begin{align}\label{eq:inner-Gamma}
     \fprod{\tfrac{1}{M}\tilde\Gamma_\bx^{k+1}-\Gamma_\bx^{k+1},\bx-\sp{\bv_\bx^k}}\leq \bD_\cX^{\bT^k}(\bx,\bv_\bx^k)-\bD_\cX^{\bT^k}(\bx,\bv_\bx^{k+1})+\tfrac{1}{2}\norm{
     \sp{\bdelta_\bx^k}}^2_{*,(\bT^k)^{-1}}.
 \end{align}
 Therefore, using \eqref{eq:inner-Gamma} within \eqref{eq:bregman-inner} we obtain
 \begin{equation}\label{eq:bregman-inner-telescope-x}
 {\small
 \begin{aligned}
     &\tfrac{1}{M}\Big(\bD^{\bT^k}_\cX(\bx,\bx^k)-\bD^{\bT^k}_\cX(\bx,\tilde\bx^{k+1})\Big) -\Big(\bD^{\bT^k}_\cX(\bx,\bx^k)-\bD^{\bT^k}_\cX(\bx,\bx^{k+1})\Big) \\
     &\leq -\tfrac{1}{M}\bD_\cX^{\bT^k}(\bx^k,\tilde\bx^{k+1})+\bD_\cX^{\bT^k}(\bx^k,\bx^{k+1})+\tfrac{1}{2}\norm{
     \sp{\bdelta_\bx^k}}^2_{*,(\bT^k)^{-1}}+\fprod{\sp{\bdelta_\bx^k}
     ,\sp{\bv_\bx^k-\bx^k}}+\bD_\cX^{\bT^k}(\bx,\bv_\bx^{k}) -\bD_\cX^{\bT^k}(\bx,\bv_\bx^{k+1}).
 \end{aligned}}%
 \end{equation}
 Similarly, for $k\geq 0$, defining $\tilde\Gamma_\by^{k+1}\triangleq \bS^k(\grad\varphi_\cY(\tilde\by^{k+1})-\grad\varphi_\cY(\by^k))$, \sp{$\Gamma_\by^{k+1}\triangleq \bS^k(\grad\varphi_\cY(\by^{k+1})-\grad\varphi_\cY(\by^k))$ and $\bdelta_\by^k\triangleq \tfrac{1}{N}\tilde{\Gamma}_\by^{k+1}-\Gamma_\by^{k+1}$,} we can show that 
 \begin{equation}\label{eq:bregman-inner-telescope-y}
 {\small
 \begin{aligned}
     &\tfrac{1}{N}\Big(\bD^{\bS^k}_\cY(\by,\by^k)-\bD^{\bS^k}_\cY(\by,\tilde\by^{k+1})\Big) -\Big(\bD^{\bS^k}_\cY(\by,\by^k)-\bD^{\bS^k}_\cY(\by,\by^{k+1})\Big) \\
     &\leq -\tfrac{1}{N}\bD_\cY^{\bS^k}(\by^k,\tilde\by^{k+1})+\bD_\cY^{\bS^k}(\by^k,\by^{k+1})+\tfrac{1}{2}\norm{
     \sp{\bdelta_\by^k}}^2_{*,(\bS^k)^{-1}}+\fprod{\sp{\bdelta_\by^k}
     ,\sp{\bv_\by^k-\by^k}} + \bD_\cY^{\bS^k}(\by,\bv_\by^{k}) - \bD_\cY^{\bS^k}(\by,\bv_\by^{k+1}).
     \end{aligned}}%
 \end{equation}
\sp{Next, define $\psi(\bx)\triangleq \|\tfrac{1}{M}\tilde{\Gamma}_\bx^{k+1}-\bT^k(\grad\varphi_\cX(\bx)-\grad\varphi_\cX(\bx^k))\|^2_{*,(\bT^k)^{-1}}$ and note that $\psi(\bx^{k+1})=\|\bdelta_\bx^k\|^2_{*,(\bT^k)^{-1}}$, $\psi(\tilde\bx^{k+1})=\|(1-\tfrac{1}{M})\tilde{\Gamma}_\bx^{k+1}\|^2_{*,(\bT^k)^{-1}}$ and $\psi(\bx^{k})=\|\tfrac{1}{M}\tilde{\Gamma}_\bx^{k+1}\|^2_{*,(\bT^k)^{-1}}$. Since $\psi(\bx)$ is separable in $\bx=[x_i]_{i=1}^M$, \eqref{eq:Ek-x-relation} implies that $\mE^k[\psi(\bx^{k+1})]=\frac{1}{M}\mE^k[\psi(\tilde\bx^{k+1})]+(1-\frac{1}{M})\psi(\bx^k)$, i.e.,}
{\small
 \begin{equation*}
     \mE^k\left[\|\sp{\bdelta_\bx^k}
     \|^2_{*,(\bT^k)^{-1}}\right]=\tfrac{1}{M}\mE^k\big[\|{(1-\tfrac{1}{M})\tilde{\Gamma}_\bx^{k+1}}\|^2_{*,(\bT^k)^{-1}}\big]+(1-\tfrac{1}{M})\mE^k\big[\|{\tfrac{1}{M}\tilde{\Gamma}_\bx^{k+1}}\|^2_{*,(\bT^k)^{-1}}\big],
 \end{equation*}}%
 which implies that
 \begin{equation}\label{eq:variance-Gamma-x}
 {\small
 \begin{aligned}
     \mE^k\left[\norm{\sp{\bdelta_\bx^k}
     }^2_{*,(\bT^k)^{-1}}\right]=\sp{\tfrac{M-1}{M^2}\mE^k\big[\|{\tilde{\Gamma}_\bx^{k+1}}\|^2_{*,(\bT^k)^{-1}}\big]}\leq \tfrac{M-1}{M^2}\mE^k\left[\norm{\grad\varphi_\cX(\tilde\bx^{k+1})-\grad\varphi_\cX(\bx^k)}_{*,\bT^k}^2\right].
 \end{aligned}}%
 \end{equation}
 \sp{Using similar arguments, we also get}
 {\small
 \begin{align}\label{eq:variance-Gamma-y}
     \mE^k\left[\norm{\sp{\bdelta_\by^k}
     }^2_{*,(\bS^k)^{-1}}\right]\leq \tfrac{N-1}{N^2}\mE^k\left[\norm{\grad\varphi_\cY(\tilde\by^{k+1})-\grad\varphi_\cY(\by^k)}_{*,\bS^k}^2\right].
 \end{align}}%
 Finally, also note that $\mE^k[\fprod{\sp{\bdelta_\by^k}
     ,\sp{\bv_\by^k-\by^k}}]=0$ since $\bv_\by^k-\by^k$ is $\Psi_{k-1}$-measurable and $\mE^k[\bdelta_\by^k]=\mathbf{0}$; similarly, we also have $\mE^k[\fprod{\sp{\bdelta_\bx^k}
     ,\sp{\bv_\bx^k-\bx^k}}]=0$. \rev{Thus, 
     using 
 \eqref{eq:bregman-inner-telescope-x}, \eqref{eq:bregman-inner-telescope-y}, \eqref{eq:variance-Gamma-x}, \eqref{eq:variance-Gamma-y} within \eqref{eq:bound-Xi} we obtain
 {\small
 \begin{align*}
 &\mE\Big[\sup_{\bz\in Z} \Big\{C_0(\bz)+ \bD_\cX^{\bT^0+(N-1)M\bL_{\bx\bx}+\frac{1}{M}\cA^0}(\bx,\bx^0)+\bD_\cY^{\bS^0+NM\bL_{\by\by}+\frac{1}{N}\tilde\cA^0}(\by,\by^0)+\sum_{k=0}^{K-1}t^k\Xi_x^k(\bx)+\sum_{k=0}^{K-1}t^k\Xi_y^k(\by)\Big\}\Big]\\
 &\leq \mE\Big[\sup_{\bz\in Z} \Big\{C_0(\bz)+ \bD_\cX^{\bT^0+(N-1)M\bL_{\bx\bx}+\frac{1}{M}\cA^0}(\bx,\bx^0)+\bD_\cY^{\bS^0+NM\bL_{\by\by}+\frac{1}{N}\tilde\cA^0}(\by,\by^0) \nonumber\\
    &\quad +\sum_{k=0}^{K-1}t^k\left(\bD_\cX^{\bT^k}(\bx,\bv_\bx^{k}) -\bD_\cX^{\bT^k}(\bx,\bv_\bx^{k+1})+\bD_\cY^{\bS^k}(\by,\bv_\by^{k}) - \bD_\cY^{\bS^k}(\by,\bv_\by^{k+1})\right)\Big\}\nonumber\\
    &\quad +\sum_{k=0}^{K-1}t^k\left(\tfrac{M-1}{2M^2}\norm{\grad\varphi_\cX(\tilde\bx^{k+1})-\grad\varphi_\cX(\bx^k)}_{*,\bT^k}^2+\tfrac{N-1}{2N^2}\norm{\grad\varphi_\cY(\tilde\by^{k+1})-\grad\varphi_\cY(\by^k)}_{*,\bS^k}^2\right) \Big]\nonumber\\
    &\leq \mE\Big[\sup_{\bz\in Z} \Big\{C_0(\bz)+\bD_\cX^{2\bT^0+(N-1)M\bL_{\bx\bx}+\frac{1}{M}\cA^0}(\bx,\bx^0) + \bD_\cY^{2\bS^0+NM\bL_{\by\by}+\frac{1}{N}\tilde\cA^0}(\by,\by^0)\Big\}\nonumber\\
 &\quad +\underbrace{\sum_{k=0}^{\infty}\mE\left[t^k\Big(\tfrac{M-1}{\sp{2}M^2}\norm{\grad\varphi_\cX(\tilde\bx^{k+1})-\grad\varphi_\cX(\bx^k)}_{*,\bT^k}^2+\tfrac{N-1}{\sp{2}N^2}\norm{\grad\varphi_\cY(\tilde\by^{k+1})-\grad\varphi_\cY(\by^k)}_{*,\bS^k}^2\Big)\right]}_{\triangleq B^\infty},\nonumber
 \end{align*}}}%
 where in the last inequality \sp{we used the fact that $\sum_{k=0}^{K-1}t^k\Big(\bD_\cX^{\bT^k}(\bx,\bv_\bx^{k})-\bD_\cX^{\bT^k}(\bx,\bv_\bx^{k+1})\Big)\leq \bD_\cX^{\bT^0}(\bx,\bx^0)$ and $\sum_{k=0}^{K-1}t^k\Big(\bD_\cY^{\bS^k}(\by,\bv_\by^{k})-\bD_\cY^{\bS^k}(\by,\bv_\by^{k+1})\Big)\leq \bD_\cY^{\bS^0}(\by,\by^0)$ since $\bv_\bx^0=\bx^0$, $\bv_\by^0=\by^0$, $t_0=1$ and $t^k\bT^{k}=t^0\bT^0$, $t^k\bS^k=t^0\bS^0$ for $k\geq 0$. Note that from \rev{Lipschitz continuity of $\grad\varphi_\cX$ and $\grad\varphi_\cY$}, and using the fact that $t^k\bT^k=t^0\bT^0=\bT^0$ and $t^k\bS^k=t^0\bS^0=\bS^0$ for $k\geq 0$, we get
 {\small
 \begin{equation*}
     B^\infty\leq \sum_{k=0}^{+\infty} \mE\Big[\tfrac{M-1}{M^2}L^2_{\varphi_\cX}\bD_\cX^{\bT^0}(\btx^{k+1},\bx^k)+\tfrac{N-1}{N^2}L^2_{\varphi_\cY}\bD_\cY^{\bS^0}(\bty^{k+1},\by^k)\Big]=C_1(\bz^*)+C_2,
 \end{equation*}}%
 where the equality follows from \eqref{eq:Ek-x-relation} and Lemma \ref{lem:result-summable}.}
\end{proof}
\subsection{Proof of Theorem~\ref{thm:meta}}
\label{sec:meta-proof}
\rev{Fix an arbitrary saddle point  of \eqref{eq:original-problem}, $\bz^*\in Z^*$, and $p\in (0,1)$. Lemma \ref{lem:result-summable} shows that $\mathbb P\left(\sup_{k\ge 0}\{\| \bx^k-\bx^*\|^2_{\bT^0}+\| \by^k-\by^*\|^2_{\bS^0}\}\leq \Delta_p\right)\geq 1-p$; therefore, for $Z\subset\cZ$ given in Definition~\ref{def:restricted-domain}, it immediately holds that $\{\bz^k\}\subset Z$ with probability at least $1-p$.}

\ey{Recall that $\bw^k_\bx= \mathbf{e}^k_\bx(\bx^k,\by^{k+1})+(N-1)\theta^k({\mathbf{e}}^k_\bx(\bx^k,\by^{k})-{\mathbf{e}}^{k}_\bx(\bx^{k-1},\by^{k-1}))$ and $\bw^k_\by\triangleq \mathbf{e}^k_\by(\bx^k,\by^{k+1})+M\theta^k({\mathbf{e}}^k_\by(\bx^k,\by^{k})-{\mathbf{e}}^{k}_\by(\bx^{k-1},\by^{k-1}))$ denote the error of sampled gradients \sp{for all $k\geq 0$}, where $\mathbf{e}^k_\bx(\bx,\by)$ and $\mathbf{e}^k_\by(\bx,\by)$ are given in Definition~\ref{def:error}. 
When $\bar\delta>0$, $\bw^k_\bx$ and $\bw^k_\by$ are the errors of using a fixed-batch sampled gradients with sizes \sp{$u\geq 1$ and $v\geq 1$}, respectively; \sp{note that when $\bar\delta=0$, we can set $u=v=1$ and $\bw^k_\bx=\mathbf{0}$ and $\bw^k_\by=\mathbf{0}$ for $k\geq 0$.}} 
Let $\Upsilon_x^k\triangleq \frac{8\erfan{(M(N-1)\theta^k)^2}}{\alpha^k}\big(\norm{\bx^k-\bx^{k-1}}_{\bC_\bx^2}^2+\norm{\by^k-\by^{k-1}}_{\bL_{\bx\by}^2}^2\big)$ and $\Upsilon_y^k\triangleq \frac{8\erfan{(NM\theta^k)^2}}{\tilde\alpha^k}\big(\norm{\bx^k-\bx^{k-1}}_{\bL_{\by\bx}^2}^2+\norm{\by^k-\by^{k-1}}_{\bC_{\by}^2}^2\big)$ for $k\geq 0$. \ey{It 
 follows \nsa{from 
 \eqref{eq:part1} and \eqref{eq:part2} in Lemma~\ref{error_sampling}  
 that} 
 {\small
 \begin{align}
 {\mE}^k\Big[M\norm{\bw_\bx^k}^2_{\nsa{*,(\cA^k)^{-1}}}-\Upsilon_x^k\Big]\leq 
 \frac{2\erfan{M^2}\delta_x^2}{v\alpha^k},\qquad{\mE}^k\Big[N\norm{\bw_\by^k}^2_{\nsa{*,(\tilde\cA^k)^{-1}}}-\Upsilon_y^k\Big]\leq \frac{2\erfan{N^2}\delta_y^2}{u\tilde\alpha^k},\quad\forall~k\geq 0.\label{eq:w-lip}
 \end{align}}}%
 Now, we are ready to prove the rate statement. 
 Since Assumption~\ref{assum:step} holds, \sp{it implies that \eqref{Qk-a} holds, i.e., $t^{k+1} Q^{k+1}(\bz)-t^k R^{k+1}(\bz) \leq -\bD^{t^{k+1}\cD^{k+1}}_\cX(\bx^{k+1},\bx^k)-\bD_\cY^{t^{k+1}\tilde\cD^{k+1}}(\by^{k+1},\by^k)$} for all $\bz\in\cZ$ and $k\geq 0$, where $Q^k(\cdot)$ and $R^{k+1}(\cdot)$ are defined in \eqref{eq:Qk-def} and \eqref{eq:Rk-def}, respectively.  
 Therefore, multiplying both sides of \eqref{eq:one-step-main} by \sa{$t^{k}$} \sp{after adding and subtracting $\Upsilon_x^k$ and $\Upsilon_y^k$ to the right-hand side, and summing over $k=0$ to $K-1$, we conclude that for any $K\geq 1$ and $\bz=(\bx,\by)\in\cZ$, 
 {\small
 \begin{align}
 \label{eq:simple-rate-bound}
 \MoveEqLeft\sum_{k=0}^{K-1}\sa{t^{k}}\left(\cL(\bx^{k+1},\by)-\cL(\bx,\by^{k+1})\right)+\sum_{k=0}^{K-1}\sa{t^{k}}(\theta^k-1)\left((M-1)(\cL(\bx^k,\by)-\cL(\bx,\by))+(N-1)(\cL(\bx,\by)-\cL(\bx,\by^k))\right) \nonumber\\
 &\leq (M-1)(t^0G_f^0(\bz)-\sa{t^{K-1}}G_f^{K}(\bz))+(N-1)(t^0G_h^0(\bz)-\sa{t^{K-1}}G_h^K(\bz))+M(t^{K-1}G^K_\Phi(\bz)-t^0G^0_\Phi(\bz)) \nonumber\\
 &\quad +t^0Q^0(\bz)-\sa{t^{K-1}}R^K(\bz) +\sum_{k=0}^{K-1}{t^{k}}\big(M\norm{\bw_\bx^k}^2_{*,(\cA^k)^{-1}}-\erfan{\Upsilon_x^k}\big)+\sum_{k=0}^{K-1}{t^{k}}\big(N\norm{\bw_\by^k}^2_{*,(\tilde\cA^k)^{-1}}\erfan{-\Upsilon_y^k}\big)\nonumber\\
 &\quad +\sum_{k=0}^{K-1}t^k\Big(\fprod{\bw_\bx^k,\nsa{\bu_\bx^k-\bx^k}}+\fprod{\bw_\by^k,\nsa{\by^k-\bu_\by^k}}+{\Xi_x^k(\bx)}+{\Xi_y^k(\by)}\Big),
 \end{align}}%
 where we used that $\Upsilon_x^k+\Upsilon_y^k\leq \bD_\cX^{\cD^k}(\bx^k,\bx^{k-1})+\bD_\cY^{\tilde\cD^k}(\by^k,\by^{k-1})$ for $k\geq 1$ and $\Upsilon_x^0+\Upsilon_y^0=0$.} Moreover,
\sp{
{\small
\begin{align*}
    \MoveEqLeft t^0Q^0(\bz)-t^{K-1}R^K(\bz)\label{eq:Q-R-bound}\\
    & \leq \bD_\cX^{\bT^0+\frac{1}{M}\cA^0}(\bx,\bx^0)+\bD_\cY^{\bS^0+\frac{1}{N}\tilde\cA^0}(\by,\by^0)-M\fprod{\grad_\by\Phi(\bx^0,\by^0),\by^0-\by}-(N-1)\fprod{\grad_\bx\Phi(\bx^0,\by^0),\bx-\bx^0}\nonumber\\
 &\quad -t^{K-1}\Big(\bD_\cX^{\bT^{K-1}}(\bx,\bx^K)+\bD_\cY^{\bS^{K-1}}(\by,\by^K)+\bD_\cX^{\bM_3^{K-1}-\cA^{K-1}}(\bx^K,\bx^{K-1})+\bD_\cY^{\bM_4^{K-1}-\tilde\cA^{K-1}}(\by^K,\by^{K-1})\Big)\nonumber\\
 &\quad +Mt^{K-1}\fprod{\grad_\by\Phi(\bx^{K-1},\by^{K-1}),\by^{K}-\by}+(N-1)t^{K-1}\fprod{\grad_\bx\Phi(\bx^{K-1},\by^{K-1}),\bx-\bx^K},\nonumber
\end{align*}}}%
\sp{which follows from dropping the non-negative terms $\bD_\cX^{\cA^{K-1}}(\bx,\bu_\bx^K)$ and $\bD_\cY^{\tilde\cA^{K-1}}(\by,\bu_\by^K)$, and using the fact that \sa{$t^0=\theta^0=1$} and $(\bx^{-1},\by^{-1})=(\bx^0,\by^0)$.} \sp{Moreover, using the definitions of $\tilde\bq_\bx^K,\tilde\bq_\by^K$ given in the proof of Lemma \ref{lem:one-step} together with Lemma~\ref{lem:pos-theta}, we get}
 {\small
 \begin{align*}
     G_\Phi^K(\bz)+\fprod{\grad_\by\Phi(\bx^{K-1},\by^{K-1}),\by^{K}-\by}&\leq \frac{1}{\nsa{M}N\sa{\theta^K}}\fprod{\sp{\tilde\bq_\by^K},\by-\by^K},\quad G_f^K(\bz)= \cL(\bx^K,\by)-\cL(\bx,\by),\\
     \fprod{\grad_\bx\Phi(\bx^{K-1},\by^{K-1}),\bx-\bx^K}-G_h^K(\bz)&\leq \frac{1}{\nsa{(N-1)}M\sa{\theta^K}}\fprod{\sp{\tilde\bq_\bx^K},\bx^K-\bx}+\cL(\bx,\by^K)-\cL(\bx,\by).
 \end{align*}}%
 Finally, \sp{\cite[(2.10) in Lemma~2 with $\alpha=0$]{nesterov2012efficiency} implies that
 }
 \begin{equation*}
 \label{eq:Lip-global}
 {\small
 \begin{aligned}
 &(M-1)G_f^0(\bz)+(N-1)G_h^0(\bz)-M G^0_\Phi(\bz)-M\fprod{\grad_\by\Phi(\bx^0,\by^0),\by^0-\by}-(N-1)\fprod{\grad_\bx\Phi(\bx^0,\by^0),\bx-\bx^0}\\
 & \leq (M-1)(\cL(\bx^0,\by)-\cL(\bx,\by))+(N-1)(\cL(\bx,\by)-\cL(\bx,\by^0))+\tfrac{(N-1)M}{2}\norm{\bx^0-\bx}^2_{\bL_{\bx\bx}}+\tfrac{NM}{2}\norm{\by^0-\by}^2_{\bL_{\by\by}}.
 \end{aligned}}%
 \end{equation*}
\sp{Due to \sp{\eqref{eq:qk-bound} with $k=K$} and the fact that $t^{K-1}=T^K\theta^K$, by using the three inequalities above within \eqref{eq:simple-rate-bound}, invoking {Jensen}'s inequality on the left side for the averaged iterate $\bar\bz^k$ defined in~\eqref{eq:avg-iterate}, and dropping the non-positive terms from the right side implied by \eqref{eq:step-size-lip}, for all $K\geq 1$, we get}
 \begin{equation}
 \label{eq:lagrange-simplified-single}
 {\small
 \begin{aligned}
 \MoveEqLeft(T_K+M-1)(\cL(\bar\bx^K,\by)-\cL(\bx,\by))+(T_K+N-1)(\cL(\bx,\by)-\cL(\bx,\bar\by^K))\\
 &\leq \sum_{k=0}^{K-1}\sa{t^{k}}\left[\cL(\bx^{k+1},\by)-\cL(\bx,\by^{k+1})+
 (\theta^k-1)\left((M-1)(\cL(\bx^k,\by)-\cL(\bx,\by))+(N-1)(\cL(\bx,\by)-\cL(\bx,\by^k))\right)\right]\\
 &\quad +\sa{t^{K-1}}\left((M-1)\big(\cL(\bx^K,\by)-\cL(\bx,\by)\big)+(N-1)\big(\cL(\bx,\by)-\cL(\bx,\by^K)\big)\right)\\
 & \leq (M-1)\big(\cL(\bx^0,\by)-\cL(\bx,\by)\big)+(N-1)\big(\cL(\bx,\by)-\cL(\bx,\by^0)\big)\\
 &\quad+\bD_\cX^{\bT^0+\frac{1}{M}\cA^0+(N-1)M\bL_{\bx\bx}}(\bx,\bx^0)+\bD_\cY^{\bS^0+\frac{1}{N}\tilde\cA^0+NM\bL_{\by\by}}(\by,\by^0)\\
&\quad-{t^{K-1}}\Big(\bD_\cX^{\bT^{K-1}-(N-1){(\gamma_1^{-1}+\gamma_2^{-1})\id_m}}(\bx,\bx^K)+\bD_\cY^{\bS^{K-1}-M({\lambda_1^{-1}+\lambda_2^{-1})\id_n}}(\by,\by^K)\Big)\\
 &\quad +\sum_{k=0}^{K-1}{t^{k}}\left(M\norm{\bw_\bx^k}^2_{*,(\cA^k)^{-1}}-\Upsilon_x^k\right) +\sum_{k=0}^{K-1}{t^{k}}\left(N\norm{\bw_\by^k}^2_{*,(\tilde\cA^k)^{-1}}-\Upsilon_y^k\right)\\
 &\quad+\sum_{k=0}^{K-1}t^k\Big(\fprod{\bw_\bx^k,\nsa{\bu_\bx^k-\bx^k}}\ey{+\fprod{\bw_\by^k,\nsa{\by^k-\bu_\by^k}}+{\Xi_x^k(\bx)}+{\Xi_y^k(\by)}\Big)},
 \end{aligned}}%
 \end{equation}
 where $T_K\triangleq\sum_{k=0}^{K-1} t^{k}$ for $K\geq 1$. \sp{Since Assumption \ref{assum:step} holds, we have $\bM_3^{K-1}\succeq 0$, $\bM_4^{K-1}\succeq 0$ and $\theta^{K-1}\geq 1$, which together imply that both $\bT^{K-1}-(N-1)({\gamma_1^{-1}+\gamma_2^{-1}})\id_m$ and $\bS^{K-1}-M({\lambda_1^{-1}+\lambda_2^{-1}})\id_n$ are diagonal matrices with positive diagonal entries; hence, dropping the nonpositive terms in \eqref{eq:lagrange-simplified-single}, first taking supremum over $\bz\in Z$, and then the expectation of both sides gives the desired result for the case $\bar\delta>0$ because of \eqref{eq:w-lip} and 
 \eqref{eq:part3} of Lemma~\ref{error_sampling}, and Lemma~\ref{lem:sup-bound}.} 
 Finally, recall that to obtain \eqref{eq:lagrange-simplified-single}, we sum \eqref{eq:one-step-main} after multiplying its both sides by $t^{k}$, and to derive \eqref{eq:lagrange-simplified-single} we use \eqref{eq:pre_bound} in the proof of Lemma~\ref{lem:one-step}. It should be emphasized that when $\bar\delta=0$, i.e., for the deterministic case, we have $\bw_\bx^k=\bw_\by^k={\bf 0}$ for all $k\geq 0$; therefore, $(d)$ and $(e)$ terms appearing in \eqref{eq:pre_bound} are both $0$, and one does not need to bound these inner products using Lemma~\ref{lem:inner-w}. Thus, anything in \eqref{eq:one-step} involving $\bw_\bx^k$ and $\bw_\by^k$ are all $0$ including $\norm{\bw_\bx^k}^2_{*,(\cA^k)^{-1}}=0$ and $\norm{\bw_\by^k}^2_{*,(\tilde\cA^k)^{-1}}=0$. Therefore, for the case $\bar\delta=0$, whenever Assumption~\ref{assum:step} holds with $\cD^k=\cA^k=\mathbf{0}$ and $\tilde\cD^k=\tilde\cA^k=\mathbf{0}$, \eqref{eq:lagrange-simplified-single} implies that  $\mE\left[\cG_Z(\bar \bz^K,T_K)\right] \leq\frac{\bar \Delta
+\sp{C_1(\bz^*)}}{T_K}$. \qed

\section{Proof of Theorem \ref{thm:bounded-stochastic}}\label{sec:proof-bounded-thm}
We first show that choosing $t^0=1$ and $t^k=\Big(\sqrt{k+1}\cdot\log(k+3)\Big)^{-1}$ for $k\geq 1$, setting $\theta^0=1$ and $\theta^k=t^{k-1}/t^k$ for $k\geq 1$, 
 and selecting $\{\tau_i^k\}_{i\in\cM}$ and $\{\sigma_j^k\}_{j\in\cN}$ as in the statement of the theorem satisfies Assumption~\ref{assum:step} for any 
 \sp{$\gamma_1,\gamma_2,\lambda_1,\lambda_2,\rev{\varsigma}>0$ and for $\{\cA^k\}$ and $\{\tilde\cA^k\}$ chosen as $\cA^k=\sp{\cA^0}/\sa{t^{k}}$ for $k\geq 0$ with $\sp{\cA^0=\alpha^0\id_{m}}$ \sp{for any $\alpha^0>0$}, and \ey{$\tilde\cA^k=\sp{\tilde\cA^0}/\sa{t^{k}}$ for $k\geq 0$ with $\tilde\cA^0=\tilde\alpha^0\id_{n}$ \sa{for any  $\tilde\alpha^0>0$}}.}
 
 \sp{Indeed, we trivially have $t^k=t^{k+1}\sa{\theta^{k+1}}$ for $k\geq 0$. The conditions $t^{k}\bT^k= t^{k+1}\bT^{k+1}$, $t^{k}\bS^k= t^{k+1}\bS^{k+1}$, $t^k\cA^k\succeq t^{k+1}\cA^{k+1}$ and $t^k\tilde\cA^k\succeq t^{k+1}\tilde\cA^{k+1}$ hold \sp{since our selection satisfies $\bT^k=\bT^0/t^k$, $\bS^k=\bS^0/t^k$, $\cA^k=\cA^0/t^k$ and $\tilde\cA^k=\tilde\cA^0/t^k$ for $k\geq 0$. Thus, \eqref{eq:step-size-theta} in Assumption~\ref{assum:step} holds. Moreover,} 
 it holds that $\theta^k\in [1,2]$ for all $k\geq 0$, which} implies that $\bM_3^k\succeq\sp{\bar{\bM}_3^k\triangleq}\bT^k-M\bL_{\bx\bx}-2M(N-1)(\gamma_1^{-1}+\gamma_2^{-1})\id_m-{\gamma_2^{-1}}(N-1)\id_m$ and $\bM_4^k\succeq\sp{\bar{\bM}_4^k\triangleq}\bS^k-N\bL_{\by\by}- 2M N({\lambda_1^{-1}+\lambda_2^{-1}})\id_n -{\gamma_2}(N-1)\bL_{\bx\by}^2$. Therefore, 
 \eqref{eq:step-size-lip} in Assumption~\ref{assum:step} holds, i.e., $\sa{t^{k+1}}(\bM_1^{k+1}\erfan{+\cD^{k+1}}) +\varsigma\bT^0\preceq t^{k}(\sp{\bar\bM_3^k}-\cA^k)$ and \sp{${t^{k+1}}(\bM_2^{k+1}{+\tilde\cD^{k+1}}) +\varsigma\bS^0\preceq t^{k}(\sp{\bar\bM_4^k}-\tilde\cA^k)$} hold for $k\geq 0$, due to $t^k=t^{k+1}\theta^{k+1}$ and $t^k\in (0,1]$ for $k\geq 0$.
 This shows that Assumption \ref{assum:step} holds, hence, the result follows immediately from  Lemma \ref{lem:result-summable} \rev{after bounding $\sum_{k=0}^\infty B_\delta^k=\sum_{k=0}^{\infty}t^k\Big(\frac{2{M^2}\delta_x^2}{v\alpha^k}+\frac{2{N^2}\delta_y^2}{u\tilde\alpha^k}\Big)$. Since $\alpha^k=\alpha^0/t^k$ and $\tilde\alpha^k=\tilde\alpha^0/t^k$ for all $k\geq 0$, it is sufficient to bound 
 {\small
 \begin{equation*}
 \sum_{k=0}^{K-1}(t^k)^2= 1+\sum_{k=1}^{K-1}(t^k)^2
 \leq 1+\sum_{k=1}^{\infty}\frac{1}{(k+1)\log^2(k+3)}\leq 2. \qed
 \end{equation*}}}%
\section{Proof of Theorem \ref{thm:single-sample}} 
\label{sec:stochastic-proof}
 {\bf Part (I):} \sp{Whenever $\bar\delta>0$,} since only a mini-batch sample gradient is used at each iteration \sp{with constant batch sizes of $u^k=u$ and $v^k=v$ for $k\geq 0$,} we need to decrease the primal and dual step-sizes to ensure that the gradient error \us{does not accumulate}.
 Step-sizes and parameters of our algorithm as stated in Theorem \ref{thm:bounded-stochastic} satisfy Assumption~\ref{assum:step} holds; therefore, the desired result follows from Theorem~\ref{thm:meta} \rev{and using the bound on $\sum_{k=0}^\infty B_\delta^k$ provided in the proof of Theorem \ref{thm:bounded-stochastic}.}
 Moreover, it is easy to verify that $T_K=\sum_{k=0}^{K-1} t^k\geq \frac{1}{\log(K+2)}\sum_{k=0}^{K-1}\frac{1}{\sqrt{k+1}}\geq \frac{2\sqrt{K+1}-1}{\log(K+2)}$ for $K\geq 1$.
 \\
 {\bf Part (II):}
  This result immediately follows from \eqref{bound_c_2} and Remark \ref{rem:gap}. \qed
\section{Proof of Theorem \ref{thm:variable-sample}} 
\label{sec:deterministic-proof}
\label{sec:proof-thm-variable}
 {\bf Part (I):} 
 \rev{Consider the deterministic setting, i.e., $\bw_\bx=\mathbf{0}$ and $\bw_\by=\mathbf{0}$, with $\theta^k=t^k=1$, $\bT^k=\sp{\bT}\triangleq \diag([1/\tau_i]_{i\in \cM})$ and $\bS^k=\bS\triangleq \diag(([1/\sigma_j]_{j\in\cN})$ for $k\geq 0$ such that $\tau_i$ for $i\in\cM$ and $\sigma_j$ for $j\in\cN$ 
as stated in \eqref{eq:step-size}. Following the same argument as in the proof of Theorem \ref{thm:bounded-stochastic}, one can verify that Assumption~\ref{assum:step} holds and the result follows from Theorem~\ref{thm:meta} for the case $\delta_x=\delta_y=0$.}\\ 
{\bf Part (II):} This result immediately follows from 
{\bf Part (I)} and Remark \ref{rem:gap}.\\ 
 \missing{\bf Part (III):} The proof is similar to that of \cite[Theorem 2. Part I.]{hamedani23a} and it is provided in Section \ref{sec:proof of conv det} of Appendix for completeness.\qed
\section{Detail of Complexity Results for Different Partitioning Strategies}\label{sec:complexity}
\sp{Given \rev{some $\bar M,\bar N\in\integers_+$, let the primal and dual dimensions $m$ and $n$ be $m=\bar M$ and $n=\bar N$,} and in the rest of this section, consider two cases: \textit{(i)} $\delta_x=\delta_y=\delta>0$, which we call the stochastic setting, and \textit{(ii)} $\delta_x=\delta_y=0$, which we call the deterministic setting -- here, $\delta^2$ denotes the variance bound for each coordinate of the stochastic oracle; hence, the variance bound for $\grad_\bx\Phi(\cdot,\cdot;\xi^x)$ and $\grad_\by\Phi(\cdot,\cdot;\xi^y)$ are \rev{$\bar M\delta^2$ and $\bar N\delta^2$}, respectively. For each setting, we discuss $4$ different primal-dual partitioning strategies.}

 {\it Scenario 1:} {\bf (Full)}. \nsa{\rev{Consider $M=N=1$, i.e., $\cM=\{1\}$ and $\cN=\{1\}$ such that $m_1=\bar M$ and $n_1=\bar N$; hence,} at each iteration we compute $\grad_{\bx}\Phi$ and $\grad_{\by}\Phi$ when $\delta=0$ \sp{or their stochastic estimates} when $\delta>0$. According to \cite{nesterov2012efficiency}, 
 \rev{$L_{xx}=\cO(\bar M)$, $L_{yy}=\cO(\bar N)$ and $L_{yx}=\cO(\sqrt{\bar M\bar N})$, \sp{which are tight (see Remark~\ref{rem:L-constants})}; therefore, $L=\max\{L_x,L_y\}=\cO(\max\{\bar M,\bar N\})$.} 
 \eqref{eq:bound-O-1-s} and \eqref{eq:bound-O-1} imply that the iteration complexity of RB-PDA is $\tilde\cO([\mathfrak{Q}(\cR_x^2+\cR_y^2)+\delta^2]^2/\varepsilon^2)$ when $\delta>0$ and $\cO(\mathfrak{Q}(\cR_x^2+\cR_y^2)/\varepsilon)$ when $\delta=0$, where \rev{$\mathfrak{Q}\triangleq \max\{\bar M, \bar N\}=\max\{m,n\}$}. Note that the variance bounds in \eqref{eq:bound-O-1-s} should be taken as \rev{$\delta_x^2=\bar M\delta^2$ and $\delta_y^2=\bar N\delta^2$}, while all $M$ and $N$ appearing as block numbers in \eqref{eq:bound-O-1-s} and \eqref{eq:bound-O-1} should be set to $1$.} 
 Moreover, each iteration 
 requires \rev{$\cO(\bar M)$ primal and $\cO(\bar N)$} dual oracles. \erfan{Thus, overall primal and dual oracle complexities for the stochastic setting are \sp{$\tilde\cO(\rev{m}[\mathfrak{Q}(\cR_x^2+\cR_y^2)+\delta^2]^2/\varepsilon^2)$} and $\tilde\cO(\rev{n}(\mathfrak{Q}(\cR_x^2+\cR_y^2)+\delta)^2/\varepsilon^2)$, and for the deterministic setting \rev{the corresponding complexities} are $\cO(\rev{m}\mathfrak{Q}(\cR_x^2+\cR_y^2)/\varepsilon)$ and  $\cO(\rev{n}\mathfrak{Q}(\cR_x^2+\cR_y^2)/\varepsilon)$, respectively.} 

 {\it Scenario 2:} {\bf (Block $\bx$)}.
 {
 \rev{Consider $M=\bar M$ and $N=1$, i.e., $\cM=\{1,\ldots, \bar M\}$ and $\cN=\{1\}$ such that $m_i=1$ for $i\in\cM$ and $n_1=\bar N$,} and at each iteration we compute $\grad_{\by}\Phi$ and $\grad_{x_i}\Phi$ for some $i\in\cM$ \sp{or their stochastic estimates}. In this case, $L_{xx}=L_{x_i x_i}=\cO(1)$ for $i\in\cM$, \rev{$L_{yy}=\cO(\bar N)$ and $L_{yx_i}=L_{yx}=\cO(\sqrt{\bar N})$ for all $i\in\cM$; therefore, $L_x=\cO(\sqrt{\bar N})$ and $L_y=\cO(\bar N)$.} \erfan{The results in \eqref{eq:bound-O-1-s} and \eqref{eq:bound-O-1} imply that the iteration complexities can be bounded by \rev{$\tilde\cO([\sqrt{\bar N}\cR_x^2+\bar N\cR_y^2+\delta^2]^2\bar M^2/\varepsilon^2)$ and $\cO((\sqrt{\bar N}\cR_x^2+\bar N\cR_y^2)\bar M/\varepsilon)$} for the stochastic and deterministic settings, respectively.} \sp{Note that the variance bounds in \eqref{eq:bound-O-1-s} should be taken as \rev{$\delta_x^2=\delta^2$ and $\delta_y^2=\bar N\delta^2$}, while all $N$ appearing as the number of dual blocks in \eqref{eq:bound-O-1-s} and \eqref{eq:bound-O-1} should be set to $1$.}
 Moreover, each iteration of this implementation requires $\cO(1)$ primal- and \rev{$\cO(\bar N)$} dual-oracle calls. Thus,
 overall primal and dual oracle complexities for the stochastic case are \rev{$\tilde\cO((\sqrt{n}\cR_x^2+n\cR_y^2+\delta^2)^2m^2/\varepsilon^2)$ and $\tilde\cO((\sqrt{n}\cR_x^2+n\cR_y^2+\delta^2)^2m^2n/\varepsilon^2)$}, and for the deterministic setting are \rev{$\cO((\sqrt{n}\cR_x^2+n\cR_y^2)m/\varepsilon)$ and  $\cO((\sqrt{n}\cR_x^2+n\cR_y^2)mn/\varepsilon)$, respectively.} 

 {\it Scenario 3:} {\bf (Block $\by$)}. \rev{Consider $M=1$ and $N=\bar N$, i.e., $\cM=\{1\}$ and $\cN=\{1,\ldots,\bar N\}$ such that $m_1=\bar M$ and $n_j=1$ for $j\in\cN$,} and at each iteration we compute $\grad_{\bx}\Phi$ and $\grad_{y_j}\Phi$ for some $j\in\cN$ or their stochastic estimates.
 In this case, $L_{yy}=L_{y_j y_j}=\cO(1)$ for $j\in\cN$, \rev{$L_{xx}=\cO(\bar M)$ and $L_{y_jx}=L_{yx}=\cO(\sqrt{\bar M})$ for all $j\in\cN$; therefore, $L_x=\cO(\bar M)$ and $L_y=\cO(\sqrt{\bar M})$.} \erfan{The results in \eqref{eq:bound-O-1-s} and \eqref{eq:bound-O-1} imply that the iteration complexities can be bounded by \rev{$\tilde\cO([\bar M\cR_x^2+\sqrt{\bar M}\cR_y^2+\delta^2]^2\bar N^2/\varepsilon^2)$ and $\cO((\bar M\cR_x^2+\sqrt{\bar M}\cR_y^2)\bar N/\varepsilon)$} for the stochastic and deterministic settings, respectively.} \sp{Note that the variance bounds in \eqref{eq:bound-O-1-s} should be taken as \rev{$\delta_x^2=\bar M\delta^2$} and $\delta_y^2=\delta^2$, while all $M$ appearing as the number of primal blocks in \eqref{eq:bound-O-1-s} and \eqref{eq:bound-O-1} should be set to $1$.}
 Moreover, each iteration of this implementation requires \rev{$\cO(\bar M)$} primal- and $\cO(1)$ dual-oracle calls. Thus,
 overall primal and dual oracle complexities for the stochastic case are \rev{$\tilde\cO((m\cR_x^2+\sqrt{m}\cR_y^2+\delta^2)^2n^2m/\varepsilon^2)$ and $\tilde\cO((m\cR_x^2+\sqrt{m}\cR_y^2+\delta^2)^2n^2/\varepsilon^2)$}, and for the deterministic setting are \rev{$\cO((m\cR_x^2+\sqrt{m}\cR_y^2)mn/\varepsilon)$ and  $\cO((m\cR_x^2+\sqrt{m}\cR_y^2)n/\varepsilon)$}, respectively. 

 {\it Scenario 4:} {\bf (Block $\bx,\by$)}. \rev{Consider $M=\bar M$ and $N=\bar N$, i.e., $\cM=\{1,\ldots,\bar M\}$ and $\cN=\{1,\ldots,\bar N\}$} such that $m_i=1$ and $n_j=1$ for $i\in\cM$ and $j\in\cN$, and at each iteration we compute $\grad_{x_i}\Phi$ and $\grad_{y_j}\Phi$ for some $i\in\cM$ and $j\in\cN$ or their stochastic estimates. Note $L_{xx}$, $L_{yy}$, $L_{xy}$ and $L_{yx}$ are all $\cO(1)$; thus, $L_x$ and $L_y$ are also $\cO(1)$. Moreover,
 at each iteration $\cO(1)$ primal- 
 and $\cO(1)$ dual-oracles are called; 
 therefore, 
 \rev{the computational complexity is directly determined by the iteration complexity. Hence,}
 the results \eqref{eq:bound-O-1-s} and \eqref{eq:bound-O-1} imply that 
 both primal and dual oracle complexities are \rev{$\tilde\cO((mn(\cR_x^2+\cR_y^2)+\max\{\tfrac{m}{n},\tfrac{n}{m}\}\delta^2)^2/\varepsilon^2)$} for the stochastic setting, and \rev{$\cO(mn(\cR_x^2+\cR_y^2)/\varepsilon)$} for the deterministic setting. \sp{Note that the variance bounds in \eqref{eq:bound-O-1-s} should be taken as $\delta_x^2=\delta^2$ and $\delta_y^2=\delta^2$ since $m_i=1$ and $n_j=1$ for $i\in\cM$ and $j\in\cN$.}
\section{Proof of Theorem~\ref{thm:variable-sample} Part (III)}\label{sec:proof of conv det}
Here we show a.s. convergence of the actual iterate sequence. 
 Lemma~\ref{lem:supermartingale} will be \sa{essential to establish the convergence result.} 
 {To fix the notation, \sa{for a random variable} $\bF:\Omega\rightarrow\reals$, $\bF(\omega)$ denotes a particular realization of $\bF$ corresponding to $\omega\in\Omega$ where $\Omega$ denotes the sample space.}


 \sa{Given an arbitrary $\omega\in\Omega$, with probability 1, $\{\bz^k(\omega)\}$ is a bounded sequence;} 
 hence, it has a convergent subsequence \sa{$\{k_q\}_{q\geq 1}\subset\integers_+$ such that} $\bz^{k_q}(\omega)\rightarrow \bz^*(\omega)$ as $q\rightarrow \infty$ for some $\bz^*(\omega)\in\cX\times\cY$ -- note that $k_q$ also depends on $\omega$ which is omitted to simplify the notation. Define $\bz^*=(\bx^*,\by^*)$ such that $\bz^*=[\bz^*(\omega)]_{\omega\in\Omega}$. 
 Moreover, Lemma \ref{lem:result-summable} implies that 
 for arbitrary $\nu>0$ 
 and for any realization $\omega\in\Omega$, there exists $N_1(\omega)$ such that for any $q\geq N_1(\omega)$, we have $\max\{\norm{\bz^{k_q}(\omega)-\bz^{k_q-1}(\omega)},~\norm{\bz^{k_q}(\omega)-\bz^{k_q+1}(\omega)}\}< \frac{\nu}{2}$. {In addition}, convergence of {subsequence} $\{\bz^{k_q}(\omega)\}$ sequence implies that there exists $N_2(\omega)$ such that for any $q\geq N_2(\omega)$, $\norm{\bz^{k_q}(\omega)-\bz^*(\omega)}< \frac{\nu}{2}$. Therefore, for $\omega\in\Omega$, letting $N_3(\omega)\triangleq \max\{N_1(\omega),N_2(\omega)\}$, we conclude that $\norm{\bz^{k_q\pm 1}(\omega)-\bz^*(\omega)}< \nu$, i.e., $\bz^{k_q\pm 1}(\omega) \rightarrow \bz^*(\omega)$ for $\omega\in\Omega\setminus\cN$ as $q\rightarrow \infty$ where $\cN$ is some measure zero set.

 Now we show that for any $\omega\in\Omega$, $\bz^*(\omega)$ is indeed a saddle point of \eqref{eq:original-problem}. 
 In particular, we show the result by considering the optimality conditions for Line~5 and Line~9 of the RB-PDA.
 Fix an arbitrary $\omega\in\Omega$ and consider the subsequence $\{k_q\}_{q\geq 1}$. For all $q\in\mathbb{Z}_+$, one can conclude that 
 {\small
 \begin{align*}
 &\tfrac{1}{\sigma_{j_{k_q}}}\Big(\grad\varphi_{\cY_{j_{k_q}}}\big(y_{j_{k_q}}^{k_q}(\omega)\big)-\grad\varphi_{\cY_{j_{k_q}}}\big(y_{j_{k_q}}^{k_q+1}(\omega)\big)\Big)+\tilde s_{j_{k_q}}^{k_q} \ey{+Nw_{y_{j_{k_q}}}^{k_q}(\omega)}\in\partial h_{j_{k_q}}\big(y_{j_{k_q}}^{k_{q}+1}(\omega)\big), \\ 
 &\tfrac{1}{\tau_{i_{k_q}}} \Big(\grad\varphi_{\cX_{i_{k_q}}}\big(x_{i_{k_q}}^{k_q}(\omega)\big)-\grad\varphi_{\cX_{i_{k_q}}}\big(x_{i_{k_q}}^{k_q+1}(\omega)\big)\Big)
 -\tilde{r}_{i_{k_q}}^{k_q }-M{w}_{x_{i_{k_q}}}^{k_q}(\omega)\in \partial f_{i_{k_q}}\big(x_{i_{k_q}}^{k_q+1}(\omega)\big), 
 \end{align*}}%
 where \ey{$\btr^k= M\grad_\bx\Phi(\bx^k,\by^{k+1})+ M(N-1)(\grad_\bx\Phi(\bx^k,\by^{k})-\grad_\bx\Phi(\bx^{k-1},\by^{k-1}))$, $\tilde\bs^k= \break M\grad_\by\Phi(\bx^k,\by^k)+ MN\theta^k(\grad_\by\Phi(\bx^k,\by^k)-\grad_\by\Phi(\bx^{k-1},\by^{k-1}))$} and $w^k_{y_j},w^k_{x_i}$ are defined in Definition \ref{def:ek}.
 Note that 
 the sequence of randomly chosen block coordinates in RB-PDA, i.e., $\{i_{k_q},j_{k_q}\}_{q\geq 1}$ is a Markov chain containing a single recurrent class. More specifically, the state space is represented by $\cM\times \cN$ and starting from a state $(i,j)\in\cM\times\cN$ the probability of eventually returning to state $(i,j)$ is strictly positive for all $i\in\cM$ and $j\in\cN$. Therefore, for any $i\in\cM$ and $j\in\cN$, we can select a further subsequence 
 $\cK^{ij}\subseteq \{k_q\}_{q\in\mathbb{Z}_+}$ such that $i_\ell=i$ and $j_\ell=j$ for all $\ell\in \cK^{ij}$. Note that $\cK^{ij}$ is an infinite subsequence with probability 1 and $\{\cK^{ij}\}_{(i,j)\in\cM\times\cN}$ is a partition of $\{k_n\}_{n\in\mathbb{Z}_+}$. Thus, for any $i\in\cM$ and $j\in\cN$, one can conclude 
 that for all $\ell\in\cK^{ij}$, 
 {\small
 \begin{subequations}\label{eq:subsublimit}
 \begin{align}
 &\frac{1}{\sigma_{j}}\Big(\grad\varphi_{\cY_{j}}\big(y_{j}^{\ell}(\omega)\big)-\grad\varphi_{\cY_{j}}\big(y_{j}^{\ell+1}(\omega)\big)\Big)+\tilde{s}_{j}^{\ell}\ey{+Nw_{y_j}^\ell(\omega)}\in\partial h_{j}\big(y_{j}^{\ell+1}(\omega)\big), \\
 &\frac{1}{\tau_{i}} \Big(\grad\varphi_{\cX_{i}}\big(x_{i}^{\ell}(\omega)\big)-\grad\varphi_{\cX_{i}}\big(x_{i}^{\ell+1}(\omega)\big)\Big)
 -\tilde{r}_{i}^{\ell }-M{w}_{x_i}^{\ell}(\omega)\in \partial f_{i}\big(x_{i}^{\ell+1}(\omega)\big).
 \end{align}
 \end{subequations}}%
 Since 
 $\cK^{ij}\subseteq \{k_q\}_{q\in\mathbb{Z}_+}$, we have that $\lim_{\ell\in\cK^{ij}}\bz^{\ell}(\omega)=\lim_{\ell\in\cK^{ij}}\bz^{\ell+1}(\omega)=\bz^*(\omega)$. {Moreover, for any $i\in \cM$ \ey{and $j\in\cN$}, \ey{$w_{x_i}^\ell= 0$ and $w_{y_j}^\ell= 0$}.}
 Using the fact that for any $i\in\cM$ and $j\in\cN$, $\grad\varphi_{\cX_i}$ and $\grad\varphi_{\cY_j}$ are continuously differentiable on $\dom f_i$ and $\dom h_j$, respectively, it follows from Theorem 24.4 in~\cite{rockafellar2015convex} that by taking the limit of both sides of \eqref{eq:subsublimit}, and noting that the step-sizes are fixed, we get $\mathbf{0}\in \grad_{x_i}\Phi\big(\bx^*(\omega),\by^*(\omega)\big)+\partial f_i(x_i^*(\omega))$, and $\mathbf{0}\in \partial h_j(y_j^*(\omega))-\grad_{y_j}\Phi(\bx^*(\omega),\by^*(\omega))$, which implies that $\bz^*(\omega)$ is a saddle point of \eqref{eq:original-problem} for any $\omega\in\Omega$. 
 Finally, since \eqref{Qk-b} and \eqref{eq:bound-iterates} are true for any saddle point \sp{$\bz^*\in Z^*$}, 
 from Lemma~\ref{lem:supermartingale}, one can conclude that $a^*=\lim_{k\rightarrow \infty} a_k \geq 0$ exists almost surely.
 On the other hand, using the fact that $\bz^{k_q\pm 1}(\omega)\rightarrow \bz^*(\omega)$ for any $\omega\in\Omega$, we have that $\lim_{q\rightarrow \infty} a_{k_q}(\omega)=0$ for any $\omega\in\Omega$;
 henceforth, $a^*=\lim_{k\rightarrow \infty} a_k=0$ almost surely, and $a_k\geq \varsigma\bD_\cX(\bx^*,\bx^k)+\varsigma\bD_\cY(\by^*,\by^k)$, for some $\varsigma>0$, implies that $\bz^k\rightarrow \bz^*$ almost surely. \qed

\bibliography{papers}
\end{document}